\documentclass[11pt,leqno]{article}
\normalfont
\renewcommand{\baselinestretch}{1.2}
\usepackage{amsfonts}

\newcommand{\QED}{\raisebox{0.5mm}{\fbox{\rule{0mm}{1.5mm}\ }}}

\newcounter{myfn}[page]

\newcommand{\myfootnote}[1]{\setcounter{footnote}{\value{myfn}}%
    \footnote{#1}\stepcounter{myfn}}

\newcommand{\IntroFigOne}{Figure 1.1}

\newcommand{\PosetAndLatticeFig}{Figure 2.1}
\newcommand{\ComponentFig}{Figure 2.2}
\newcommand{\ExampleFigList}{Figures 2.1 and 2.2} 
\newcommand{\CartanMatrixEntriesTable}{Figure 2.3}
\newcommand{\DynkinDiagramSymmetryLemma}{Lemma 2.1}
\newcommand{\DualCharacterLemma}{Lemma 2.2}
\newcommand{\ConnectedLabellingPosetExists}{Lemma 2.3}
\newcommand{\RGFProp}{Proposition 2.4} 
\newcommand{\ThreeElementGridPosetFigure}{Figure 3.1} 
\newcommand{\GridPosets}{Figure 3.2}
\newcommand{\GridPosetsII}{Figure 3.3}
\newcommand{\GridPosetsFigureList}{Figures 3.2 and 3.3}
\newcommand{\MaxPropertyFigList}{Figures 2.1, 3.2, and 3.3} 
\newcommand{\WeightsLemma}{Lemma 3.1}

\newcommand{\CartanTable}{Figure 4.1}
\newcommand{\FundPosets}{Figure 4.2}
\newcommand{\FundLatticeIdealsFigure}{Figure 4.3}
\newcommand{\FundLatticeTableauxFigure}{Figure 4.4}
\newcommand{\IdealTableauFigureList}{Figures 4.3 and 4.4}
\newcommand{\RestrictionChart}{Figure 4.5}
\newcommand{\CaseBFig}{Figure 4.6}
\newcommand{\StructureResultForFundamentals}{Lemma 4.1}
\newcommand{\StructureResult}{Proposition 4.2}
\newcommand{\NewRemark}{Remark 4.3}
\newcommand{\DistinctLemma}{Lemma 4.4}
\newcommand{\TableauProp}{Proposition 4.5}
\newcommand{\OrderRemark}{Remark 4.6}
\newcommand{\TableauWeightProp}{Proposition 4.7}
\newcommand{\FigsForTProp}{Figures 4.4 and 4.5}
\newcommand{\SectFivePropsList}{Propositions 4.7 and 4.2}
\newcommand{\LitTableauProp}{Definition 5.1}
\newcommand{\RankTwoData}{Figure 5.1}
\newcommand{\LitAdmissibleColumns}{Figure 5.2}
\newcommand{\WeightPresBijection}{Proposition 5.2}
\newcommand{\CharacterProposition}{Theorem 5.3}
\newcommand{\CharacterCorollary}{Corollary 5.4}
\newcommand{\BtwoMolevFigureLattice}{Figure 6.1}
\newcommand{\BtwoMolevFigureList}{Figure 6.1}
\newcommand{\NotLieFundPosets}{Figure 6.2}
\newcommand{\NotLieComboLattice}{Figure 6.2}
\newcommand{\ClassificationTheoremOne}{Theorem 6.1}
\newcommand{\ClassificationTheoremTwo}{Theorem 6.2}
\newcommand{\ClassificationTheorems}{Theorems 6.1 and 6.2}

\newcommand{\IntroNum}{1}
\newcommand{\DefsNum}{2}
\newcommand{\GridNum}{3}
\newcommand{\SemiNum}{4}
\newcommand{\CharNum}{5}
\newcommand{\ExamplesNum}{6}

 \newcommand{\relt}{\mathbf{r}}
\newcommand{\selt}{\mathbf{s}} \newcommand{\telt}{\mathbf{t}}
\newcommand{\uelt}{\mathbf{u}} 
 \newcommand{\xelt}{\mathbf{x}}
\newcommand{\yelt}{\mathbf{y}} 

\newcommand{\ecolor}{\mathbf{edgecolor}}
\newcommand{\vcolor}{\mathbf{vertexcolor}}

\newcommand{\color}{\mathbf{color}}
\newcommand{\comp}{\mathbf{comp}}

\newcommand{\Twt}{\mathbf{tableauwt}}

\newcommand{\dichromatic}{two-color }

\newcommand{\TabSet}{\mathcal{S}_{\mathfrak{g}}(\lambda)}
\newcommand{\Lba}{L_{\mathfrak{g}}^{\beta\alpha}(\lambda)}
\newcommand{\Lab}{L_{\mathfrak{g}}^{\alpha\beta}(\lambda)}

\newcommand{\LAtwo}{L_{A_{2}}}

\newcommand{\LBtwo}{L_{C_{2}}}

\newcommand{\LGtwo}{L_{G_{2}}}

\newcommand{\PBtwoba}{P_{C_{2}}^{\beta\alpha}(\lambda)}

\newcommand{\digriddelta}{(P, \leq_{_{P}}, 
\mathbf{chain}: P \longrightarrow [m], 
\color: P \longrightarrow \Delta)}
\newcommand{\digrid}{(P, \leq_{_{P}}, 
\mathbf{chain}: P \longrightarrow [m], 
\color: P \longrightarrow \{\alpha, \beta\})}

\newcommand{\myarrow}[1]{\stackrel{#1}{\rightarrow}}

\newcommand{\VertexForPosets}[2]{
\setlength{\unitlength}{1cm}
\begin{picture}(0,0)
\put(-0.38,0){
\begin{picture}(0,0)
\put(0,0){\circle*{0.15}} 
\put(-0.6,-0.1){\footnotesize $v_{#1}$}
\put(0.2,-0.1){\footnotesize $#2$}
\end{picture}
}
\end{picture}
}

\newcommand{\VertexForLatticeI}[1]{
\setlength{\unitlength}{1.5cm}
\begin{picture}(0,0)
\put(-0.25,0){
\begin{picture}(0,0)
\put(0,0){\circle*{0.1}} 
\put(-0.35,0.0){\footnotesize $\telt_{#1}$}
\end{picture}
}
\end{picture}
}

\newcommand{\NEEdgeLabelForLatticeI}[1]{
\setlength{\unitlength}{1.5cm}
\begin{picture}(0,0)
\put(-0.25,0){
\begin{picture}(0,0)
\put(0.4,0.4){\footnotesize $#1$} 
\end{picture}
}
\end{picture}
}

\newcommand{\NEEdgeLabelForLatticeII}[1]{
\setlength{\unitlength}{1cm}
\begin{picture}(0,0)
\put(-0.25,0){
\begin{picture}(0,0)
\put(0.25,0.4){\footnotesize $#1$} 
\end{picture}
}
\end{picture}
}

\newcommand{\NEEdgeLabelForLatticeIII}[1]{
\begin{picture}(0,0)
\put(-0.25,0){
\begin{picture}(0,0)
\put(0.35,0.35){\footnotesize $#1$} 
\end{picture}
}
\end{picture}
}

\newcommand{\NWEdgeLabelForLatticeI}[1]{
\setlength{\unitlength}{1.5cm}
\begin{picture}(0,0)
\put(-0.25,0){
\begin{picture}(0,0)
\put(-0.525,0.4){\footnotesize $#1$} 
\end{picture}
}
\end{picture}
}

\newcommand{\NWEdgeLabelForLatticeII}[1]{
\setlength{\unitlength}{1cm}
\begin{picture}(0,0)
\put(-0.25,0){
\begin{picture}(0,0)
\put(-0.75,0.4){\footnotesize $#1$} 
\end{picture}
}
\end{picture}
}

\newcommand{\NWEdgeLabelForLatticeIII}[1]{
\begin{picture}(0,0)
\put(-0.25,0){
\begin{picture}(0,0)
\put(-0.575,0.35){\footnotesize $#1$} 
\end{picture}
}
\end{picture}
}

\newcommand{\VerticalEdgeLabelForLatticeI}[1]{
\setlength{\unitlength}{1.5cm}
\begin{picture}(0,0)
\put(-0.25,0){
\begin{picture}(0,0)
\put(-0.05,0.4){\footnotesize $#1$} 
\end{picture}
}
\end{picture}
}

\newcommand{\abtab}[3]{\setlength{\unitlength}{0.35cm}\begin{picture}(2,2)
\put(0,-0.5){\begin{picture}(2,2) \put(0,0){\line(0,1){2}}
\put(1,0){\line(0,1){2}} \put(2,1){\line(0,1){1}}
\put(0,0){\line(1,0){1}} \put(0,1){\line(1,0){2}}
\put(0,2){\line(1,0){2}} \put(0.3,1.2){\footnotesize #1}
\put(0.3,0.2){\footnotesize #2} \put(1.3,1.2){\footnotesize #3}
\end{picture}}
\end{picture}}

\newcommand{\lengthonecolumnA}[1]{
\footnotesize \begin{tabular}{|c|}
\hline $#1$ \\%
\hline 
\end{tabular}
}

\newcommand{\lengthonecolumnB}[1]{
\footnotesize \begin{tabular}{|c|c|}
\hline $#1$ & $#1$ \\%
\hline 
\end{tabular}
}

\newcommand{\lengthonecolumnG}[2]{
\footnotesize \begin{tabular}{|c|c|c|c|c|c|}
\hline $#1$ & $#1$ & $#1$ & $#2$ & $#2$ & $#2$\\%
\hline 
\end{tabular}
}

\newcommand{\lengthtwocolumnA}[2]{
\footnotesize \begin{tabular}{|c|}
\hline $#1$ \\%
\hline $#2$\\
\hline
\end{tabular}
}

\newcommand{\lengthtwocolumnB}[4]{
\footnotesize \begin{tabular}{|c|c|}
\hline $#1$ & $#3$\\%
\hline $#2$ & $#4$\\
\hline
\end{tabular}
}

\newcommand{\lengthtwocolumnG}[8]{
\footnotesize \begin{tabular}{|c|c|c|c|c|c|}
\hline $#1$ & $#1$ & $#3$ & $#5$ & $#7$ & $#7$\\%
\hline $#2$ & $#2$ & $#4$ & $#6$ & $#8$ & $#8$\\
\hline
\end{tabular}
}

\newcommand{\lengthonecolumn}[1]{
\setlength{\unitlength}{0.35cm}\begin{picture}(1.15,1)
\put(0,0){\begin{picture}(1,1) \put(0,0){\line(0,1){1}}
\put(1,0){\line(0,1){1}} \put(0,0){\line(1,0){1}}
\put(0,1){\line(1,0){1}} \put(0.3,0.2){\footnotesize #1}
\end{picture}}
\end{picture}}

\newcommand{\lengthtwocolumn}[2]{
\setlength{\unitlength}{0.35cm}\begin{picture}(1.15,1.5)
\put(0,-0.5){\begin{picture}(1,2) \put(0,0){\line(0,1){2}}
\put(1,0){\line(0,1){2}} \put(0,0){\line(1,0){1}}
\put(0,1){\line(1,0){1}} \put(0,2){\line(1,0){1}}
\put(0.3,1.2){\footnotesize #1} \put(0.3,0.2){\footnotesize #2}
\end{picture}}
\end{picture}}

\newcommand{\smallonebox}[1]{
\setlength{\unitlength}{0.2cm}\begin{picture}(1,1)
\put(0,0){\begin{picture}(1,1) \put(0,0){\line(0,1){1}}
\put(1,0){\line(0,1){1}} \put(0,0){\line(1,0){1}}
\put(0,1){\line(1,0){1}} \put(0.2,0.2){\tiny #1}
\end{picture}}
\end{picture}}

\newcommand{\smalltwobox}[2]{
\setlength{\unitlength}{0.2cm}\begin{picture}(1,1.5)
\put(0,-0.5){\begin{picture}(1,2) \put(0,0){\line(0,1){2}}
\put(1,0){\line(0,1){2}} \put(0,0){\line(1,0){1}}
\put(0,1){\line(1,0){1}} \put(0,2){\line(1,0){1}}
\put(0.2,1.2){\tiny #1} \put(0.2,0.2){\tiny #2}
\end{picture}}
\end{picture}}

\newcommand{\DynkinOne}{\setlength{\unitlength}{1cm} \hspace*{0.1cm}
\begin{picture}(1,0.5)
\put(0,0){\circle*{0.2}} \put(1,0){\circle*{0.2}}
\put(0,0){\line(1,0){1}} \put(-0.05,0.2){\scriptsize $i$}
\put(0.95,0.2){\scriptsize $j$}
\end{picture}
\hspace*{0.2cm}}
\newcommand{\DynkinTwo}{\setlength{\unitlength}{1cm}
\hspace*{0.1cm}
\begin{picture}(1,0.5)
\put(0,0){\circle*{0.2}} \put(1,0){\circle*{0.2}}
\put(0,-0.1){\line(1,0){1}} \put(0,0.1){\line(1,0){1}}
\put(0.35,0){\line(2,1){0.36}} \put(0.35,0){\line(2,-1){0.36}}
\put(-0.05,0.2){\scriptsize $i$} \put(0.95,0.2){\scriptsize $j$}
\end{picture}
\hspace*{0.2cm}}
\newcommand{\DynkinThree}{\setlength{\unitlength}{1cm} \hspace*{0.1cm}
\begin{picture}(1,0.5)
\put(0,0){\circle*{0.2}} \put(1,0){\circle*{0.2}}
\put(0,-0.1){\line(1,0){1}} \put(0,0.1){\line(1,0){1}}
\put(0,0){\line(1,0){1}} \put(0.35,0){\line(2,1){0.36}}
\put(0.35,0){\line(2,-1){0.36}} \put(-0.05,0.2){\scriptsize $i$}
\put(0.95,0.2){\scriptsize $j$}
\end{picture}
\hspace*{0.2cm}}
\newcommand{\DynkinZero}{\setlength{\unitlength}{1cm} \hspace*{0.1cm}
\begin{picture}(1,0.5)
\put(0,0){\circle*{0.2}} \put(1,0){\circle*{0.2}}
\put(-0.05,0.2){\scriptsize $i$} \put(0.95,0.2){\scriptsize $j$}
\end{picture}
\hspace*{0.2cm}}
\newcommand{\AOneAOneAlpha}{
\setlength{\unitlength}{1cm}
\begin{picture}(0.5,0.75)
\put(-0.45,0.2){\footnotesize $v_{1}$}
\put(0,0.25){\circle*{0.15}} 
\put(0.2,0.25){\footnotesize $\alpha$}
\end{picture}
}
\newcommand{\AOneAOneBeta}{
\setlength{\unitlength}{1cm}
\begin{picture}(0.5,0.75)
\put(-0.45,0.2){\footnotesize $v_{1}$}
\put(0,0.25){\circle*{0.15}} 
\put(0.2,0.25){\footnotesize $\beta$}
\end{picture}
}
\newcommand{\ATwoAlpha}{
\setlength{\unitlength}{1cm}
\begin{picture}(2,2)
\put(0,0.25){
\begin{picture}(2,2)
\put(0.5,-0.1){\footnotesize $v_{2}$}
\put(1,0){\circle*{0.15}} 
\put(1.2,-0.1){\footnotesize $\beta$}
\put(-0.5,0.9){\footnotesize $v_{1}$}
\put(0,1){\circle*{0.15}} 
\put(0.2,0.9){\footnotesize $\alpha$}
\put(1,0){\line(-1,1){1}} 
\end{picture}
}
\end{picture}
}
\newcommand{\ATwoBeta}{
\setlength{\unitlength}{1cm}
\begin{picture}(2,2)
\put(0,0.25){
\begin{picture}(2,2)
\put(0.5,-0.1){\footnotesize $v_{2}$}
\put(1,0){\circle*{0.15}} 
\put(1.2,-0.1){\footnotesize $\alpha$}
\put(-0.5,0.9){\footnotesize $v_{1}$}
\put(0,1){\circle*{0.15}} 
\put(0.2,0.9){\footnotesize $\beta$}
\put(1,0){\line(-1,1){1}} 
\end{picture}
}
\end{picture}
}
\newcommand{\BTwoAAlpha}{
\setlength{\unitlength}{1cm}
\begin{picture}(3,4)
\put(0,0.75){
\begin{picture}(3,3)
\put(1.5,-0.1){\footnotesize $v_{3}$}
\put(2,0){\circle*{0.15}} 
\put(2.2,-0.1){\footnotesize $\alpha$}
\put(0.5,0.9){\footnotesize $v_{2}$}
\put(1,1){\circle*{0.15}} 
\put(1.2,0.9){\footnotesize $\beta$}
\put(-0.5,1.9){\footnotesize $v_{1}$}
\put(0,2){\circle*{0.15}} 
\put(0.2,1.9){\footnotesize $\alpha$}
\put(2,0){\line(-1,1){2}} 
\end{picture}
}
\end{picture}
}
\newcommand{\BTwoBBeta}{
\setlength{\unitlength}{1cm}
\begin{picture}(4,4)
\put(1,0.25){
\begin{picture}(3,3.5)
\put(0.5,-0.1){\footnotesize $v_{4}$}
\put(1,0){\circle*{0.15}} 
\put(1.2,-0.1){\footnotesize $\beta$}
\put(-0.5,0.9){\footnotesize $v_{3}$}
\put(0,1){\circle*{0.15}} 
\put(0.2,0.9){\footnotesize $\alpha$}
\put(0.5,1.9){\footnotesize $v_{2}$}
\put(1,2){\circle*{0.15}} 
\put(1.2,1.9){\footnotesize $\alpha$}
\put(-0.5,2.9){\footnotesize $v_{1}$}
\put(0,3){\circle*{0.15}} 
\put(0.2,2.9){\footnotesize $\beta$}
\put(0,1){\line(1,1){1}} 
\put(1,0){\line(-1,1){1}} 
\put(1,2){\line(-1,1){1}} 
\end{picture}
}
\end{picture}
}
\newcommand{\GTwoAAlpha}{
\setlength{\unitlength}{1cm}
\begin{picture}(4,6)
\put(0.25,1.25){
\begin{picture}(4,6)
\put(2.5,-0.1){\footnotesize $v_{6}$}
\put(3,0){\circle*{0.15}} 
\put(3.2,-0.1){\footnotesize $\alpha$}
\put(1.5,0.9){\footnotesize $v_{5}$}
\put(2,1){\circle*{0.15}} 
\put(2.2,0.9){\footnotesize $\beta$}
\put(0.5,1.9){\footnotesize $v_{4}$}
\put(1,2){\circle*{0.15}} 
\put(1.2,1.9){\footnotesize $\alpha$}
\put(1.5,2.9){\footnotesize $v_{3}$}
\put(2,3){\circle*{0.15}} 
\put(2.2,2.9){\footnotesize $\alpha$}
\put(0.5,3.9){\footnotesize $v_{2}$}
\put(1,4){\circle*{0.15}} 
\put(1.2,3.9){\footnotesize $\beta$}
\put(-0.5,4.9){\footnotesize $v_{1}$}
\put(0,5){\circle*{0.15}} 
\put(0.2,4.9){\footnotesize $\alpha$}
\put(1,2){\line(1,1){1}} 
\put(3,0){\line(-1,1){2}} 
\put(2,3){\line(-1,1){2}} 
\end{picture}
}
\end{picture}
}
\newcommand{\GTwoBBeta}{
\setlength{\unitlength}{1cm}
\begin{picture}(4,8)
\put(0.25,0.25){
\begin{picture}(4,8)
\put(1.5,-0.1){\footnotesize $v_{10}$}
\put(2,0){\circle*{0.15}} 
\put(2.2,-0.15){\footnotesize $\beta$}
\put(0.5,0.9){\footnotesize $v_{9}$}
\put(1,1){\circle*{0.15}} 
\put(1.2,0.9){\footnotesize $\alpha$}
\put(1.5,1.9){\footnotesize $v_{8}$}
\put(2,2){\circle*{0.15}} 
\put(2.2,1.9){\footnotesize $\alpha$}
\put(0.5,2.9){\footnotesize $v_{6}$}
\put(1,3){\circle*{0.15}} 
\put(1.2,2.9){\footnotesize $\beta$}
\put(2.5,2.9){\footnotesize $v_{7}$}
\put(3,3){\circle*{0.15}} 
\put(3.2,2.9){\footnotesize $\alpha$}
\put(-0.5,3.9){\footnotesize $v_{4}$}
\put(0,4){\circle*{0.15}} 
\put(0.2,3.9){\footnotesize $\alpha$}
\put(1.5,3.9){\footnotesize $v_{5}$}
\put(2,4){\circle*{0.15}} 
\put(2.2,3.9){\footnotesize $\beta$}
\put(0.5,4.9){\footnotesize $v_{3}$}
\put(1,5){\circle*{0.15}} 
\put(1.2,4.9){\footnotesize $\alpha$}
\put(1.5,5.9){\footnotesize $v_{2}$}
\put(2,6){\circle*{0.15}} 
\put(2.2,5.9){\footnotesize $\alpha$}
\put(0.5,6.9){\footnotesize $v_{1}$}
\put(1,7){\circle*{0.15}} 
\put(1.2,6.9){\footnotesize $\beta$}
\put(1,1){\line(1,1){2}} 
\put(1,3){\line(1,1){1}} 
\put(0,4){\line(1,1){2}} 
\put(2,0){\line(-1,1){1}} 
\put(2,2){\line(-1,1){2}} 
\put(3,3){\line(-1,1){2}} 
\put(2,6){\line(-1,1){1}} 
\end{picture} 
}
\end{picture}
}

\newcommand{\AtwoAlphaIdeals}{
\setlength{\unitlength}{1.1cm}
\begin{picture}(3,4)
\put(0.25,3.5){$\LAtwo(1,0)$}
\put(1.5,1){\circle*{0.125}}
\put(1.5,2){\circle*{0.125}}
\put(1.5,3){\circle*{0.125}}
\put(1.7,0.95){\footnotesize $\emptyset$}
\put(1.6,1.95){\footnotesize $\langle 2 \rangle$}
\put(1.6,2.95){\footnotesize $\langle 1 \rangle$}
\put(1.5,1){\line(0,1){2}}
\put(1.45,1.45){\footnotesize $\beta$}
\put(1.45,2.45){\footnotesize $\alpha$}
\end{picture}
}

\newcommand{\AtwoAlphaTableaux}{
\setlength{\unitlength}{1.1cm}
\begin{picture}(3,4)
\put(0.25,3.5){$\LAtwo(1,0)$}
\put(1.5,1){\circle*{0.125}}
\put(1.5,2){\circle*{0.125}}
\put(1.5,3){\circle*{0.125}}
\put(1.6,0.95){\smallonebox{3}}
\put(1.6,1.95){\smallonebox{2}}
\put(1.6,2.95){\smallonebox{1}}
\put(0.65,0.95){\footnotesize (0,-1)}
\put(0.65,1.95){\footnotesize (-1,1)}
\put(0.7,2.95){\footnotesize (1,0)}
\put(1.5,1){\line(0,1){2}}
\put(1.45,1.45){\footnotesize $\beta$}
\put(1.45,2.45){\footnotesize $\alpha$}
\end{picture}
}

\newcommand{\AtwoBetaIdeals}{
\setlength{\unitlength}{1.1cm}
\begin{picture}(3,4)
\put(0.25,3.5){$\LAtwo(0,1)$}
\put(1.5,1){\circle*{0.125}}
\put(1.5,2){\circle*{0.125}}
\put(1.5,3){\circle*{0.125}}
\put(1.7,0.95){\footnotesize $\emptyset$}
\put(1.6,1.95){\footnotesize $\langle 2 \rangle$}
\put(1.6,2.95){\footnotesize $\langle 1 \rangle$}
\put(1.5,1){\line(0,1){2}}
\put(1.45,1.45){\footnotesize $\alpha$}
\put(1.45,2.45){\footnotesize $\beta$}
\end{picture}
}

\newcommand{\AtwoBetaTableaux}{
\setlength{\unitlength}{1.1cm}
\begin{picture}(3,4)
\put(0.25,3.5){$\LAtwo(0,1)$}
\put(1.5,1){\circle*{0.125}}
\put(1.5,2){\circle*{0.125}}
\put(1.5,3){\circle*{0.125}}
\put(1.6,0.95){\smalltwobox{2}{3}}
\put(1.6,1.95){\smalltwobox{1}{3}}
\put(1.6,2.95){\smalltwobox{1}{2}}
\put(0.65,0.95){\footnotesize (-1,0)}
\put(0.65,1.95){\footnotesize (1,-1)}
\put(0.7,2.95){\footnotesize (0,1)}
\put(1.5,1){\line(0,1){2}}
\put(1.45,1.45){\footnotesize $\alpha$}
\put(1.45,2.45){\footnotesize $\beta$}
\end{picture}
}

\newcommand{\BtwoAlphaIdeals}{
\setlength{\unitlength}{1.1cm}
\begin{picture}(3,5)
\put(0.25,4.5){$\LBtwo(1,0)$}
\put(1.5,1){\circle*{0.125}}
\put(1.5,2){\circle*{0.125}}
\put(1.5,3){\circle*{0.125}}
\put(1.5,4){\circle*{0.125}}
\put(1.7,0.95){\footnotesize $\emptyset$}
\put(1.6,1.95){\footnotesize $\langle 3 \rangle$}
\put(1.6,2.95){\footnotesize $\langle 2 \rangle$}
\put(1.6,3.95){\footnotesize $\langle 1 \rangle$}
\put(1.5,1){\line(0,1){3}}
\put(1.45,1.45){\footnotesize $\alpha$}
\put(1.45,2.45){\footnotesize $\beta$}
\put(1.45,3.45){\footnotesize $\alpha$}
\end{picture}
}

\newcommand{\BtwoAlphaTableaux}{
\setlength{\unitlength}{1.1cm}
\begin{picture}(3,5)
\put(0.25,4.5){$\LBtwo(1,0)$}
\put(1.5,1){\circle*{0.125}}
\put(1.5,2){\circle*{0.125}}
\put(1.5,3){\circle*{0.125}}
\put(1.5,4){\circle*{0.125}}
\put(1.6,0.95){\smallonebox{4}}
\put(1.6,1.95){\smallonebox{3}}
\put(1.6,2.95){\smallonebox{2}}
\put(1.6,3.95){\smallonebox{1}}
\put(0.65,0.95){\footnotesize (-1,0)}
\put(0.65,1.95){\footnotesize (1,-1)}
\put(0.65,2.95){\footnotesize (-1,1)}
\put(0.7,3.95){\footnotesize (1,0)}
\put(1.5,1){\line(0,1){3}}
\put(1.45,1.45){\footnotesize $\alpha$}
\put(1.45,2.45){\footnotesize $\beta$}
\put(1.45,3.45){\footnotesize $\alpha$}
\end{picture}
}

\newcommand{\BtwoBetaIdeals}{
\setlength{\unitlength}{1.1cm}
\begin{picture}(3,6)
\put(0.25,5.5){$\LBtwo(0,1)$}
\put(1.5,1){\circle*{0.125}}
\put(1.5,2){\circle*{0.125}}
\put(1.5,3){\circle*{0.125}}
\put(1.5,4){\circle*{0.125}}
\put(1.5,5){\circle*{0.125}}
\put(1.7,0.95){\footnotesize $\emptyset$}
\put(1.6,1.95){\footnotesize $\langle 4 \rangle$}
\put(1.6,2.95){\footnotesize $\langle 3 \rangle$}
\put(1.6,3.95){\footnotesize $\langle 2 \rangle$}
\put(1.6,4.95){\footnotesize $\langle 1 \rangle$}
\put(1.5,1){\line(0,1){4}}
\put(1.45,1.45){\footnotesize $\beta$}
\put(1.45,2.45){\footnotesize $\alpha$}
\put(1.45,3.45){\footnotesize $\alpha$}
\put(1.45,4.45){\footnotesize $\beta$}
\end{picture}
}

\newcommand{\BtwoBetaTableaux}{
\setlength{\unitlength}{1.1cm}
\begin{picture}(3,6)
\put(0.25,5.5){$\LBtwo(0,1)$}
\put(1.5,1){\circle*{0.125}}
\put(1.5,2){\circle*{0.125}}
\put(1.5,3){\circle*{0.125}}
\put(1.5,4){\circle*{0.125}}
\put(1.5,5){\circle*{0.125}}
\put(1.6,0.95){\smalltwobox{3}{4}}
\put(1.6,1.95){\smalltwobox{2}{4}}
\put(1.6,2.95){\smalltwobox{2}{3}}
\put(1.6,3.95){\smalltwobox{1}{3}}
\put(1.6,4.95){\smalltwobox{1}{2}}
\put(0.65,0.95){\footnotesize (0,-1)}
\put(0.65,1.95){\footnotesize (-2,1)}
\put(0.7,2.95){\footnotesize (0,0)}
\put(0.65,3.95){\footnotesize (2,-1)}
\put(0.7,4.95){\footnotesize (0,1)}
\put(1.5,1){\line(0,1){4}}
\put(1.45,1.45){\footnotesize $\beta$}
\put(1.45,2.45){\footnotesize $\alpha$}
\put(1.45,3.45){\footnotesize $\alpha$}
\put(1.45,4.45){\footnotesize $\beta$}
\end{picture}
}

\newcommand{\GtwoAlphaIdeals}{
\setlength{\unitlength}{1.1cm}
\begin{picture}(3,8)
\put(0.25,7.5){$\LGtwo(1,0)$}
\put(1.5,1){\circle*{0.125}}
\put(1.5,2){\circle*{0.125}}
\put(1.5,3){\circle*{0.125}}
\put(1.5,4){\circle*{0.125}}
\put(1.5,5){\circle*{0.125}}
\put(1.5,6){\circle*{0.125}}
\put(1.5,7){\circle*{0.125}}
\put(1.7,0.95){\footnotesize $\emptyset$}
\put(1.6,1.95){\footnotesize $\langle 6 \rangle$}
\put(1.6,2.95){\footnotesize $\langle 5 \rangle$}
\put(1.6,3.95){\footnotesize $\langle 4 \rangle$}
\put(1.6,4.95){\footnotesize $\langle 3 \rangle$}
\put(1.6,5.95){\footnotesize $\langle 2 \rangle$}
\put(1.6,6.95){\footnotesize $\langle 1 \rangle$}
\put(1.5,1){\line(0,1){6}}
\put(1.45,1.45){\footnotesize $\alpha$}
\put(1.45,2.45){\footnotesize $\beta$}
\put(1.45,3.45){\footnotesize $\alpha$}
\put(1.45,4.45){\footnotesize $\alpha$}
\put(1.45,5.45){\footnotesize $\beta$}
\put(1.45,6.45){\footnotesize $\alpha$}
\end{picture}
}

\newcommand{\GtwoAlphaTableaux}{
\setlength{\unitlength}{1.1cm}
\begin{picture}(3,8)
\put(0.25,7.5){$\LGtwo(1,0)$}
\put(1.5,1){\circle*{0.125}}
\put(1.5,2){\circle*{0.125}}
\put(1.5,3){\circle*{0.125}}
\put(1.5,4){\circle*{0.125}}
\put(1.5,5){\circle*{0.125}}
\put(1.5,6){\circle*{0.125}}
\put(1.5,7){\circle*{0.125}}
\put(1.6,0.95){\smallonebox{7}}
\put(1.6,1.95){\smallonebox{6}}
\put(1.6,2.95){\smallonebox{5}}
\put(1.6,3.95){\smallonebox{4}}
\put(1.6,4.95){\smallonebox{3}}
\put(1.6,5.95){\smallonebox{2}}
\put(1.6,6.95){\smallonebox{1}}
\put(0.65,0.95){\footnotesize (-1,0)}
\put(0.65,1.95){\footnotesize (1,-1)}
\put(0.65,2.95){\footnotesize (-2,1)}
\put(0.7,3.95){\footnotesize (0,0)}
\put(0.65,4.95){\footnotesize (2,-1)}
\put(0.65,5.95){\footnotesize (-1,1)}
\put(0.7,6.95){\footnotesize (1,0)}
\put(1.5,1){\line(0,1){6}}
\put(1.45,1.45){\footnotesize $\alpha$}
\put(1.45,2.45){\footnotesize $\beta$}
\put(1.45,3.45){\footnotesize $\alpha$}
\put(1.45,4.45){\footnotesize $\alpha$}
\put(1.45,5.45){\footnotesize $\beta$}
\put(1.45,6.45){\footnotesize $\alpha$}
\end{picture}
}

\newcommand{\GtwoBetaIdeals}{
\setlength{\unitlength}{1.1cm}
\begin{picture}(5,12)
\put(1.75,11.5){$\LGtwo(0,1)$}
\put(5,6){\circle*{0.125}}
\put(2,5){\circle*{0.125}}
\put(2,7){\circle*{0.125}}
\put(3,1){\circle*{0.125}}
\put(3,2){\circle*{0.125}}
\put(3,3){\circle*{0.125}}
\put(3,4){\circle*{0.125}}
\put(3,6){\circle*{0.125}}
\put(3,8){\circle*{0.125}}
\put(3,9){\circle*{0.125}}
\put(3,10){\circle*{0.125}}
\put(3,11){\circle*{0.125}}
\put(4,5){\circle*{0.125}}
\put(4,7){\circle*{0.125}}
\put(3.15,10.95){\footnotesize $\langle 1 \rangle$}
\put(3.15,9.95){\footnotesize $\langle 2 \rangle$}
\put(3.15,8.95){\footnotesize $\langle 3 \rangle$}
\put(3.15,7.95){\footnotesize $\langle 4,5 \rangle$}
\put(4.15,7){\footnotesize $\langle 4,7 \rangle$}
\put(1.45,6.95){\footnotesize $\langle 5 \rangle$}
\put(2.15,5.95){\footnotesize $\langle 6,7 \rangle$}
\put(5.1,5.95){\footnotesize $\langle 4 \rangle$}
\put(1.45,4.95){\footnotesize $\langle 7 \rangle$}
\put(4.15,4.9){\footnotesize $\langle 6 \rangle$}
\put(3.15,3.95){\footnotesize $\langle 8 \rangle$}
\put(3.1,2.95){\footnotesize $\langle 9 \rangle$}
\put(3.1,1.95){\footnotesize $\langle 10 \rangle$}
\put(3.2,0.95){\footnotesize $\emptyset$}
\put(3,1){\line(0,1){3}}
\put(3,8){\line(0,1){3}}
\put(3,4){\line(-1,1){1}}
\put(3,4){\line(1,1){1}}
\put(2,5){\line(1,1){2}}
\put(4,5){\line(-1,1){2}}
\put(2,7){\line(1,1){1}}
\put(4,7){\line(-1,1){1}}
\put(4,5){\line(1,1){1}}
\put(5,6){\line(-1,1){1}}
\put(2.95,10.45){\footnotesize $\beta$}
\put(2.95,9.45){\footnotesize $\alpha$}
\put(2.95,8.45){\footnotesize $\alpha$}
\put(2.45,7.45){\footnotesize $\alpha$}
\put(3.45,7.45){\footnotesize $\beta$}
\put(2.45,6.45){\footnotesize $\beta$}
\put(3.45,6.45){\footnotesize $\alpha$}
\put(4.45,6.45){\footnotesize $\alpha$}
\put(2.45,5.45){\footnotesize $\beta$}
\put(3.45,5.45){\footnotesize $\alpha$}
\put(4.45,5.45){\footnotesize $\alpha$}
\put(2.45,4.45){\footnotesize $\alpha$}
\put(3.45,4.45){\footnotesize $\beta$}
\put(2.95,3.45){\footnotesize $\alpha$}
\put(2.95,2.45){\footnotesize $\alpha$}
\put(2.95,1.45){\footnotesize $\beta$}
\end{picture}
}

\newcommand{\GtwoBetaTableaux}{
\setlength{\unitlength}{1.1cm}
\begin{picture}(5,12)
\put(1.75,11.5){$\LGtwo(0,1)$}
\put(5,6){\circle*{0.125}}
\put(2,5){\circle*{0.125}}
\put(2,7){\circle*{0.125}}
\put(3,1){\circle*{0.125}}
\put(3,2){\circle*{0.125}}
\put(3,3){\circle*{0.125}}
\put(3,4){\circle*{0.125}}
\put(3,6){\circle*{0.125}}
\put(3,8){\circle*{0.125}}
\put(3,9){\circle*{0.125}}
\put(3,10){\circle*{0.125}}
\put(3,11){\circle*{0.125}}
\put(4,5){\circle*{0.125}}
\put(4,7){\circle*{0.125}}
\put(3.15,10.95){\smalltwobox{1}{2}}
\put(3.15,9.95){\smalltwobox{1}{3}}
\put(3.15,8.95){\smalltwobox{1}{4}}
\put(3.15,7.95){\smalltwobox{1}{5}}
\put(4.15,7){\smalltwobox{1}{6}}
\put(1.45,6.95){\smalltwobox{2}{5}}
\put(2.35,5.95){\smalltwobox{2}{6}}
\put(5.1,5.95){\smalltwobox{1}{7}}
\put(1.45,4.95){\smalltwobox{3}{6}}
\put(4.15,4.9){\smalltwobox{2}{7}}
\put(3.15,3.95){\smalltwobox{3}{7}}
\put(3.1,2.95){\smalltwobox{4}{7}}
\put(3.1,1.95){\smalltwobox{5}{7}}
\put(3.1,0.95){\smalltwobox{6}{7}}
\put(2.25,10.95){\footnotesize (0,1)}
\put(2.2,9.95){\footnotesize (3,-1)}
\put(2.25,8.95){\footnotesize (1,0)}
\put(2.2,7.95){\footnotesize (-1,1)}
\put(3.15,7){\footnotesize (2,-1)}
\put(2.25,6.95){\footnotesize (-3,2)}
\put(3.15,5.95){\footnotesize (0,0)}
\put(4.15,5.95){\footnotesize (0,0)}
\put(2.25,4.95){\footnotesize (3,-2)}
\put(3.15,4.9){\footnotesize (-2,1)}
\put(2.2,3.9){\footnotesize (1,-1)}
\put(2.2,2.95){\footnotesize (-1,0)}
\put(2.2,1.95){\footnotesize (-3,1)}
\put(2.2,0.95){\footnotesize (0,-1)}
\put(3,1){\line(0,1){3}}
\put(3,8){\line(0,1){3}}
\put(3,4){\line(-1,1){1}}
\put(3,4){\line(1,1){1}}
\put(2,5){\line(1,1){2}}
\put(4,5){\line(-1,1){2}}
\put(2,7){\line(1,1){1}}
\put(4,7){\line(-1,1){1}}
\put(4,5){\line(1,1){1}}
\put(5,6){\line(-1,1){1}}
\put(2.95,10.45){\footnotesize $\beta$}
\put(2.95,9.45){\footnotesize $\alpha$}
\put(2.95,8.45){\footnotesize $\alpha$}
\put(2.45,7.45){\footnotesize $\alpha$}
\put(3.45,7.45){\footnotesize $\beta$}
\put(2.45,6.45){\footnotesize $\beta$}
\put(3.45,6.45){\footnotesize $\alpha$}
\put(4.45,6.45){\footnotesize $\alpha$}
\put(2.45,5.45){\footnotesize $\beta$}
\put(3.45,5.45){\footnotesize $\alpha$}
\put(4.45,5.45){\footnotesize $\alpha$}
\put(2.45,4.45){\footnotesize $\alpha$}
\put(3.45,4.45){\footnotesize $\beta$}
\put(2.95,3.45){\footnotesize $\alpha$}
\put(2.95,2.45){\footnotesize $\alpha$}
\put(2.95,1.45){\footnotesize $\beta$}
\end{picture}
}

\begin{document}

\newpage
\setcounter{page}{1} 
\renewcommand{\baselinestretch}{1}

\vspace*{-0.7in}
\hfill {\footnotesize April 28, 2007}

\begin{center}
{\LARGE \bf Distributive lattices defined for 
representations of\\ rank two semisimple Lie algebras}

L.\ Wyatt Alverson II$^{1}$, Robert G.\ Donnelly$^{1}$, Scott J.\ 
Lewis$^{1}$, Marti McClard$^{2}$, Robert Pervine$^{1}$, Robert A.\ 
Proctor$^{3}$, N.\ J.\ Wildberger$^{4}$

\vspace*{-0.05in} 
$^{1}$Department of Mathematics and Statistics, Murray State
University, Murray, KY 42071

\vspace*{-0.05in} 
$^{2}$Department of Mathematics, University of Tennessee, Knoxville, TN 
37996 

\vspace*{-0.05in} 
$^{3}$Department of Mathematics, 
University of North Carolina, Chapel Hill, NC 27599

\vspace*{-0.05in} 
$^{4}$School of Mathematics, University of New South Wales, Sydney, NSW 2052
\end{center}

\begin{abstract}
For a rank two root system and a pair of nonnegative 
integers, using only elementary combinatorics we construct two posets.  
The constructions are uniform across the root systems 
$A_{1} \oplus A_{1}$, $A_{2}$, 
$C_{2}$, and $G_{2}$.  Examples appear in \GridPosetsFigureList.  We then form the 
distributive lattices of order ideals of these posets.  
\CharacterCorollary\  
gives elegant quotient-of-products expressions for the rank generating 
functions of these lattices (thereby providing answers to a 1979 question 
of Stanley). Also, 
\CharacterProposition\ describes how these lattices provide a new 
combinatorial setting for the Weyl characters of representations of rank 
two semisimple Lie algebras.  Most of these lattices are new;  the rest of 
them (or related structures) have arisen in work of Stanley, 
Kashiwara, Nakashima, Littelmann, and Molev.  In a future 
paper, one author shows that the posets constructed here form a Dynkin 
diagram-indexed answer to a combinatorially posed classification 
question.  In a companion paper, some of these 
lattices are used to explicitly construct some  representations 
of rank two semisimple Lie algebras.  This implies that these lattices 
are strongly Sperner.
\end{abstract}



%==================================================================
%\newpage
\vspace{2ex} 
\noindent
{\Large \bf \IntroNum.\ \ Introduction}

One of the earliest combinatorial forays into Lie representation theory 
was Stanley's \cite{StanUnim} in 1979.  Certain polynomials arising from representations 
which had elegant quotient-of-product forms captured his attention.  He 
observed that some of these polynomials were the rank generating functions 
of certain distributive lattices.  In Problem 3 of \cite{StanUnim} he asked if 
further distributive lattices could be found which would be 
associated to more of the polynomials.  Consider the 
poset  ``$\mathbf{2} \times 
\mathbf{3}$'' shown in \IntroFigOne, 
the product of chains with 2 and 3 elements. 
Its lattice $L(2,3) = J(\mathbf{2} \times 
\mathbf{3})$ of order ideals is shown in \IntroFigOne.  
Stanley knew that the rank generating function for the general case 
$L(k,n+1-k) = J\big(\mathbf{k} \times (\mathbf{n+1-k})\big)$ 
satisfies the identity 
\[\sum N_{j}q^{j} = \frac{(1-q^{n+1})(1-q^{n})\cdots(1-q^{n+2-k})}
{(1-q^{k})(1-q^{k-1})\cdots(1-q)},\]
where $N_{j}$ is the number of order ideals in $\mathbf{k} \times 
(\mathbf{n+1-k})$ 
with $j$ elements.  
The right hand side is the ``Gaussian coefficient''  $q$-analog of the binomial 
coefficient ${n+1 \choose k}$.  It is also a shifted version of the principal 
specialization of the Weyl character for the 
$k$th fundamental representation of the %special linear 
Lie algebra $\mathfrak{sl}(n+1,\mathbb{C})$, 
the rank $n$ simple Lie algebra of type $A$. 
These considerations led Stanley to introduce 
the more general distributive lattices $L(\lambda,n+1)$, whose elements are 
semistandard tableaux of shape $\lambda$ with entries from  
$\{1,2,\ldots,n+1\}$.  Similar 
identities hold for the rank generating functions of these lattices. 
Stanley was aware that the polynomial associated to the 
``last'' fundamental representation of the Lie algebra 
$\mathfrak{sp}(2n,\mathbb{C})$ 
specializes to the $(n+1)$st Catalan number $\frac{2}{n+2}{2n+1 
\choose n}$ 
when $q$ is set to 1.  Thus 
the principal specialization of the Weyl character for that 
representation is a $q$-analog to the $(n+1)$st Catalan number.  The 
second author of this paper constructed a poset $P_{n}$ such that the 
distributive lattice $L_{n} = J(P_{n})$ of its order ideals has rank generating 
function $\frac{1-q^{2}}{1-q^{n+2}}{2n+1 \choose n}_{q}$, a shifted 
version of the principal specialization. 
So the total number of order ideals from $P_{n}$ is the $(n+1)$st Catalan 
number. This result now appears as part (ccc) of Exercise 6.19 of 
\cite{StanText2}.  See \IntroFigOne\ for the poset $P_{3}$; it has 14 order 
ideals. 

\vspace*{-0.15in}
%\begin{figure}[ht]
\begin{center}
\hspace*{0.35in}
\setlength{\unitlength}{0.6cm}
\begin{picture}(3,6.5)
\put(-3.5,5.55){\IntroFigOne}
\put(-2,0.85){
\begin{picture}(3,4)
\put(-1.5,2.6){$\mathbf{2} \times \mathbf{3}$}
\multiput(1,0)(1,1){3}{\circle*{0.15}}
\multiput(0,1)(1,1){3}{\circle*{0.15}}
\multiput(1,0)(-1,1){2}{\line(1,1){2}}
\multiput(1,0)(1,1){3}{\line(-1,1){1}}
\end{picture}}
\end{picture}
\hspace*{1in}
\begin{picture}(4,6.5)
\put(-5.25,5){$L(2,3) = J(\mathbf{2} \times 
\mathbf{3})$}
\put(1,0){\circle*{0.15}}
\put(1,1){\circle*{0.15}}
\put(0,2){\circle*{0.15}}
\put(2,2){\circle*{0.15}}
\put(1,3){\circle*{0.15}}
\put(3,3){\circle*{0.15}}
\put(0,4){\circle*{0.15}}
\put(2,4){\circle*{0.15}}
\put(1,5){\circle*{0.15}}
\put(1,6){\circle*{0.15}}
\put(1,0){\line(0,1){1}}
\put(1,5){\line(0,1){1}}
\put(1,1){\line(-1,1){1}}
\put(1,1){\line(1,1){2}}
\put(0,2){\line(1,1){2}}
\put(2,2){\line(-1,1){2}}
\put(3,3){\line(-1,1){2}}
\put(0,4){\line(1,1){1}}
\end{picture}
\hspace*{0.25in}
\begin{picture}(8,6)
\put(0,0.5){
\begin{picture}(5,5)
\put(0,4.5){$P_{3}$}
\multiput(2.5,0.5)(-1,1){3}{\circle*{0.15}}
\multiput(3.5,1.5)(-1,1){3}{\circle*{0.15}}
\multiput(4.5,2.5)(-1,1){3}{\circle*{0.15}}
\qbezier(1.5,1.5)(2.5,1)(3.5,1.5)
\qbezier(0.5,2.5)(2.5,1.5)(4.5,2.5)
\qbezier(1.5,3.5)(2.5,3)(3.5,3.5)
\thicklines
\put(2.6,2){\vector(1,0){0}}
\put(2.6,3.25){\vector(1,0){0}}
\put(2.6,1.25){\vector(1,0){0}}
\thinlines
\put(2.5,0.5){\line(-1,1){2}}
\put(3.5,1.5){\line(-1,1){2}}
\put(4.5,2.5){\line(-1,1){2}}
\put(3.5,1.5){\line(1,1){1}}
\put(0.5,2.5){\line(1,1){1}}
\end{picture}}
\put(5.5,2.75){\large $\cong$}
\put(6.5,0.75){
\begin{picture}(2,4.5)
\put(2,0){\circle*{0.15}}
\put(1,1){\circle*{0.15}}
\put(1,1.5){\circle*{0.15}}
\put(0,2){\circle*{0.15}}
\put(2,2.25){\circle*{0.15}}
\put(0,2.5){\circle*{0.15}}
\put(1,3){\circle*{0.15}}
\put(1,3.5){\circle*{0.15}}
\put(2,4.5){\circle*{0.15}}
\put(2,0){\line(-1,1){2}}
\put(1,1.5){\line(-1,1){1}}
\put(0,2){\line(1,1){1}}
\put(0,2.5){\line(1,1){2}}
\put(1,1.5){\line(4,3){1}}
\put(2,2.25){\line(-4,3){1}}
\put(1,1){\line(0,1){0.5}}
\put(0,2){\line(0,1){0.5}}
\put(1,3){\line(0,1){0.5}}
\end{picture}}
\end{picture}
\end{center}
%\end{figure}

\vspace*{-0.15in}
Here is Stanley's 1979 question: 
\begin{center}
\noindent 
\parbox{5.5in}{Problem 3: Which other of the polynomials [of Theorem 
1] are the rank generating functions for distributive lattices (or 
perhaps just posets) ``naturally associated'' with the root system $R$?}
\end{center}
We supply answers to this question by 
constructing eight two-parameter families of 
distributive lattices.  By the proof of \CharacterCorollary, 
their rank generating functions are the shifted 
principal specializations of the Weyl characters of the 
irreducible finite dimensional 
representations of the four 
rank two semisimple Lie algebras $A_{1} \oplus A_{1}$, $A_{2}$, 
$C_{2}$, and $G_{2}$.  
The answers for $C_{2}$ and $G_{2}$ are largely new. 
Given a rank two semisimple Lie algebra $\mathfrak{g}$ and a pair of non-negative 
integers, we first construct two 
``$\mathfrak{g}$-semistandard posets''.  
The ``$\mathfrak{g}$-semistandard'' distributive lattices are then 
obtained by ordering the order ideals of these posets by inclusion. 
For example, the choices of $G_{2}$ and non-negative 
integer parameters $(2,2)$ specify the last poset in each of 
\GridPosetsFigureList.  According to \CharacterCorollary, both 
of these posets have $\frac{1}{5!}(3 \cdot 3 \cdot 6 \cdot 9 \cdot 12 
\cdot 15) = 729 = 3^{6}$ order ideals.  
The rank generating 
function for both of the corresponding $G_{2}$-semistandard lattices is 
\begin{eqnarray*}
RGF_{G_{2}}(2,2,q) & =\!\!= & \frac{(1-q^{3})(1-q^{3})(1-q^{6})(1-q^{9})(1-q^{12})(1-q^{15})}
{(1-q)(1-q)(1-q^{2})(1-q^{3})(1-q^{4})(1-q^{5})}. 
\end{eqnarray*} 
Hence our lattices $L_{G_{2}}^{\beta\alpha}(2,2)$ and 
$L_{G_{2}}^{\alpha\beta}(2,2)$ 
are two answers to 
Problem 3.

Since the 1970's, the ``zoo'' of finite sets of combinatorial objects which 
are enumerated by quotient-of-products formulas has grown to include 
dozens of species.  
Here \CharacterCorollary\ adds $L_{C_{2}}^{\beta\alpha}(a,b)$, 
$L_{C_{2}}^{\alpha\beta}(a,b)$, $L_{G_{2}}^{\beta\alpha}(a,b)$, and 
$L_{G_{2}}^{\alpha\beta}(a,b)$ to this zoo; they   
are analogs to the lattices  $L(\lambda,n)$.  Our  
$\mathfrak{g}$-semistandard lattices are uniformly 
defined across the four types of rank two semisimple Lie algebras.  
\CharacterCorollary\ 
also notes that the sequence of rank cardinalities 
for any $\mathfrak{g}$-semistandard lattice is symmetric and unimodal.

Only familiarity with the most basic Lie representation theory in 
\cite{Hum} is 
needed to read this paper.  
The central fact needed is that each irreducible finite 
dimensional representation of a semisimple Lie algebra of rank $n$ has 
a unique $n$-variate Weyl character.  

Some of the rank two $\mathfrak{g}$-semistandard lattices constructed 
here (or related objects) have appeared in the work of Stanley, 
Kashiwara, Nakashima, 
Littelmann, Molev, and several of the authors.   However, taken as a whole, each 
of the  $C_{2}$- and  $G_{2}$-families of $\mathfrak{g}$-semistandard 
lattices is new.  
The $A_{2}$-family of semistandard 
lattices here are 
the $n=2$ case of the 
$L(\lambda,n)$ lattices introduced in \cite{StanUnim}.  
A certain infinite subfamily of the $C_{2}$-semistandard lattices 
appeared in \cite{DLP1} as the $n=2$ case of the ``Molev lattices''  
$L_{B}^{Mol}(k,2n)$.  
A certain infinite 
subfamily of the $G_{2}$-semistandard lattices was studied in 
\cite{DLP1}. 

Let  $\mathfrak{g}$  be a semisimple Lie algebra of rank  $n$.  
Various data and 
structures have been associated to each irreducible finite dimensional 
representation of  $\mathfrak{g}$, starting with its highest weight 
and dimension.  Once certain 
subalgebras of  $\mathfrak{g}$  have been fixed, the multiset of weights of a 
representation is determined.  The Weyl character of the representation is 
the generating function for this multiset of weights.  It is a Laurent 
polynomial in  $n$  variables.  The polynomials that caught Stanley's eye 
were shifted versions of the ``principal specializations'' of the Weyl 
characters to the variable  $q$.  A finer version of Stanley's 1979  
question is:  For each Weyl character, find a distributive lattice with 
weighted vertices such that the sum of these weights is the Weyl 
character.  If the lattice elements are assigned weights in a 
reasonable manner, then a shifted version of the principal specialization 
will be the lattice's rank generating function.  An explicit 
combinatorial answer to this question (such as a lattice constructed from 
tableaux) will include a solution to the ``labelling problem'' for the 
character:  the lattice elements will be combinatorial objects 
which can be used as labels for the weights.  The problem considered 
here is a stronger version of this finer version of Stanley's question 
for  $n$ = 2.  The ``stronger'' aspect is described below.

Going further, fixing Chevalley generators for  $\mathfrak{g}$  
and basis vectors for 
the representation space determines the data consisting of the entries of 
the representing matrices for the generators.  At this point in several 
papers (such as \cite{DonSupp}), the second author introduces the ``supporting 
graph'' combinatorial structure.  This is a directed graph whose edges are 
colored by the simple roots of  $\mathfrak{g}$.  
The edges colored by simple root 
$\alpha_{i}$  indicate 
which basis vectors arise with non-zero coefficients when the Chevalley 
generators  $x_{i}$  and  $y_{i}$  of  $\mathfrak{g}$  
act on the various basis vectors.  This 
graph is actually the Hasse diagram of a poset.  Several of the 
authors have been 
able to find distributive lattice supporting graphs for many 
representations 
\cite{DonSupp}, \cite{DLP1}, 
\cite{ADLP}.  
The crystal graph is another combinatorial 
structure associated to a representation.  
For irreducible representations,  
the crystal graph is a supporting graph when the 
weight multiplicities are all one.  Such representations have only one 
supporting graph.  
But otherwise  
the crystal graph has fewer edges than do the most efficient supporting 
graphs; then it cannot support its representation.  

Our original goal for developing  $\mathfrak{g}$-semisimple lattices was 
to supply  
uniformly constructed labels and supporting graphs for explicit 
realizations of all irreducible representations of any rank two 
semisimple Lie algebra 
$\mathfrak{g}$.  Suppose a vertex-weighted 
edge-colored directed graph is proposed to be a supporting graph of a 
representation of  $\mathfrak{g}$:  
In addition to its vertex weighting agreeing with 
the Weyl character, its edge-coloring must also satisfy certain conditions 
specified by the Cartan matrix of  $\mathfrak{g}$.  
(But these conditions alone are not 
sufficient for the graph to be a supporting graph.)  If these 
edge-coloring necessary conditions are also met, the proposed graph is said to 
be a ``splitting poset'' for the representation.  
The edge-coloring conditions are the embodiment of Stanley's request 
that the lattices be {\em natural} with respect to the Lie theory. 
Here is the ``stronger'' 
aspect of our main problem:  Not only do we require that the weighting of 
their elements agree with a given Weyl character, we seek edge-colored 
distributive lattices which are splitting posets.  Our answer to this 
question consists of the  $\mathfrak{g}$-semistandard lattices:  
\StructureResult\  
verifies that the edge colorings satisfy the necessary conditions and 
our main result 
\CharacterProposition\ verifies that the vertex weights agree with the character.  
(The latter verification implies that the order ideals in the  
$\mathfrak{g}$-semistandard posets can serve as new weight labels for these 
representations.)

The necessary edge-color conditions correspond to the Serre relations 
(S1) and (S3) of Proposition 18.1 of \cite{Hum}.  Given a splitting poset for a 
representation of  $\mathfrak{g}$,  if edge coefficients for the actions of the 
generators  $x_{i}$  and  $y_{i}$  of  $\mathfrak{g}$  
can be found that satisfy the relations 
(S2), then a result of Kashiwara's implies that the remaining Serre 
relations  (S$_{ij}^{+}$)  and  (S$_{ij}^{-}$)  are automatically satisfied.  
In certain cases 
the companion paper \cite{ADLP} is able to attain our original goal by 
assigning coefficients satisfying (S2) 
to the edges of the lattices introduced here.  
So \cite{ADLP} presents explicit realizations for the following irreducible 
representations of rank two simple Lie algebras, indexed by their type and 
highest weights:  $A_{2}(a\omega_{1}+b\omega_{2})$,  $C_{2}(a\omega_{1})$,  
$C_{2}(b\omega_{2})$,  $C_{2}(\omega_{1}+b\omega_{2})$,  
$G_{2}(a\omega_{1})$,  
$G_{2}(\omega_{2})$,  for $a, b \geq 0$.  
Since the  $\mathfrak{g}$-semistandard 
lattices are supporting graphs here, as in \cite{PrGZ} they can be 
seen to be ``strongly Sperner''.   
The results of this paper facilitated the new  
$C_{2}(\omega_{1}+b\omega_{2})$ 
constructions and made it possible to now present the 
supporting lattices for all of these representations in a uniform fashion.

It can be shown that the $\mathfrak{g}$-semistandard lattices corresponding to the 
other rank two irreducible representations cannot support their 
corresponding representations.  But to state \CharacterCorollary, one only needs 
to know that the lattice at hand is a splitting poset for an irreducible 
representation.  Hence the beautiful product identities may be written 
down for the rank generating functions of all $\mathfrak{g}$-semistandard lattices.  
The necessary edge-coloring conditions are so strong that the second 
author has been able to prove that the Dynkin diagram-indexed  
$\mathfrak{g}$-semistandard lattices constitute the entire answer to a purely 
combinatorial problem \cite{DonTwoColor}.  See \ClassificationTheorems.

The positioning of splitting posets (in general;  $\mathfrak{g}$-semistandard lattices 
in particular) in the world of combinatorial structures associated to 
representations is vaguely similar in spirit to the positioning of  
crystal graphs:  both the lattices and crystal graphs superimpose 
additional combinatorial structure onto the data contained in the Weyl 
character, but neither can always support the actions of the corresponding 
representations.  In Section \ExamplesNum\ we indicate how some splitting posets may 
hopefully someday be used instead of crystal graphs for some purposes, 
such as computing tensor products.

Many of the definitions, lemmas, and propositions developed in this paper 
are needed in  \cite{ADLP}.  Some of them will also be used 
in \cite{DW} to explicitly construct 
many families of splitting posets for the simple Lie algebras  
$A_{n}$, $B_{n}$, $C_{n}$, 
$D_{n}$, $E_{6}$, $E_{7}$, and $G_{2}$.

Section \DefsNum\ presents definitions and some preliminary and background 
results.  The reader should initially browse this section and then consult 
it as needed.  
Section \GridNum\ further considers ``grid posets'' which were 
introduced in \cite{ADLP} and whose definition is 
purely combinatorial. 
 \WeightsLemma\ is the key decomposition result.  It is proved here and 
used in \cite{ADLP} and \cite{DW}.  Section \SemiNum\ introduces  
$\mathfrak{g}$-semistandard posets,  
$\mathfrak{g}$-semistandard lattices,  and 
$\mathfrak{g}$-semistandard tableaux.  Section \CharNum\ shows 
that the elements of these lattices match up with tableaux presented in 
Littelmann's \cite{Lit1}.  This match-up yields our main results.  Section 
\ExamplesNum\  
contains further remarks and problems.

%==================================================================
%\newpage
\vspace{2ex} 

\noindent
{\Large \bf \DefsNum.\ \ Definitions and preliminary results}

The reference for standard combinatorics material is \cite{Stanley}, 
and the reference for standard representation theory material is \cite{Hum}. 
We use ``$R$'' (and ``$Q$'') 
as a generic name for most of 
the combinatorial 
structures defined in this section: ``edge-colored directed graph,'' 
``vertex-colored directed graph,'' 
``ranked poset,''  
``splitting poset''.  
The letter ``$P$'' is reserved for posets and ``vertex-colored'' 
posets that arise as posets of join irreducibles for 
distributive lattices. The letter ``$L$'' is reserved for 
distributive lattices and ``edge-colored'' distributive 
lattices.  All posets are finite.  We identify a poset with its Hasse 
diagram. 

Let $I$ be any set.  An {\em edge-colored directed graph with
edges colored by the set $I$} is a directed graph $R$ with vertex set 
$\mathcal{V}(R)$ and directed-edge set $\mathcal{E}(R)$ together
with a function $\ecolor_{R}\, :\, \mathcal{E}(R) \longrightarrow I$ 
assigning to each edge of $R$ a {\em color} from the
set $I$.  If an edge $\selt \rightarrow \telt$ in $R$ is assigned 
$i \in I$, we write $\selt \myarrow{i} \telt$.  See 
\PosetAndLatticeFig. 
For $i \in I$, we let 
$\mathcal{E}_{i}(R)$ denote the set of edges in $R$ of color $i$, 
so $\mathcal{E}_{i}(R) = \ecolor_{R}^{-1}(i)$.  
If $J$ is a subset of $I$,
remove all edges from $R$ whose colors are not in $J$; connected
components of the resulting edge-colored directed graph are called
{\em J-components} of $R$. 
For any $\telt$ in $R$ and any $J \subset I$, we let 
$\mathbf{comp}_{J}(\telt)$ denote the $J$-component of $R$ 
containing $\telt$.  
The {\em dual} $R^{*}$ is the edge-colored directed graph whose vertex 
set $\mathcal{V}(R^{*})$ is the set of symbols $\{\telt^{*}\}_{\telt{\in}R}$ 
together with colored edges 
$\mathcal{E}_{i}(R^{*}) := 
\{\telt^{*} \myarrow{i} \selt^{*}\, |\, 
\selt \myarrow{i} \telt \in \mathcal{E}_{i}(R)\}$ 
for each $i \in I$. Let $Q$ be another edge-colored
directed graph with edge colors from $I$. If $R$ and $Q$ have disjoint 
vertex sets, then the {\em disjoint sum} $R
\oplus Q$ is the expected edge-colored directed graph. If $\mathcal{V}(Q) 
\subseteq \mathcal{V}(R)$ and $\mathcal{E}_{i}(Q) \subseteq 
\mathcal{E}_{i}(R)$ for each $i \in
I$, then $Q$ is an {\em edge-colored subgraph} of $R$. 
Let $R \times Q$ denote the expected edge-colored directed graph with vertex 
set $\mathcal{V}(R) \times \mathcal{V}(Q)$.
The notion of isomorphism for edge-colored directed 
graphs is as expected.  (See \cite{ADLP} if any ``expected'' statement 
is unclear.) 
If $R$ is an edge-colored directed graph with edges colored by the 
set $I$, and if $\sigma\, :\, I \longrightarrow I'$ is a mapping of 
sets, then we let $R^{\sigma}$ be the edge-colored directed graph with 
edge color function $\ecolor_{R^{\sigma}} := \sigma \circ 
\ecolor_{R}$. 
We call $R^{\sigma}$ 
a {\em recoloring} of $R$. Observe that $(R^{*})^{\sigma} \cong 
(R^{\sigma})^{*}$. 
We similarly define a {\em vertex-colored directed graph} 
with a function $\vcolor_{R}\, :\, \mathcal{V}(R) \longrightarrow I$ that 
assigns colors to the vertices of $R$.  
In this context, we speak of the {\em dual vertex-colored directed 
graph} $R^{*}$, the {\em disjoint sum} of two vertex-colored directed graphs 
with disjoint vertex sets, 
{\em isomorphism} of vertex-colored directed graphs, {\em recoloring}, 
etc.  See \ExampleFigList.  

For $\selt$ and $\telt$ in a poset $R$, there is a directed 
edge $\selt \rightarrow \telt$ in the Hasse diagram of $R$ 
if and only if 
$\telt$ {\em covers} $\selt$. 
So terminology 
that applies to directed graphs ({\em connected}, 
{\em edge-colored}, {\em dual}, {\em vertex-colored}, etc)  
will also apply to posets.   
The vertex $\selt$ and the edge $\selt \rightarrow \telt$ are  
{\em below} $\telt$, and the
%=============================
% Text before Fig 2.1, 2.2
%=============================

\newpage
%\begin{figure}[thb]
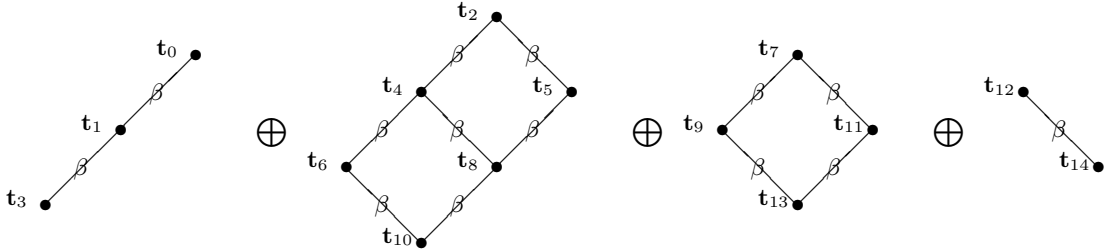
\begin{center}
\PosetAndLatticeFig:  A vertex-colored poset $P$ and an edge-colored 
lattice $L$.

\setlength{\unitlength}{1cm}
\begin{picture}(3,3.5)
\put(0,3){\begin{picture}(3,3.5)
\put(0,3.5){$P$}
%Vertex labels
\put(4,2){\circle*{0.15}} 
\put(3.4,1.9){\footnotesize $v_{6}$}
\put(4.2,1.9){\footnotesize $\beta$} 
\put(1,1){\circle*{0.15}}
\put(0.4,0.9){\footnotesize $v_{5}$} 
\put(1.2,0.9){\footnotesize $\alpha$} 
\put(2,2){\circle*{0.15}} 
\put(1.4,1.9){\footnotesize $v_{4}$} 
\put(2.2,1.9){\footnotesize $\alpha$}
\put(3,3){\circle*{0.15}} 
\put(2.4,2.9){\footnotesize $v_{3}$}
\put(3.2,2.9){\footnotesize $\alpha$} 
\put(0,2){\circle*{0.15}}
\put(-0.6,1.9){\footnotesize $v_{2}$} 
\put(0.2,1.9){\footnotesize $\beta$} 
\put(1,3){\circle*{0.15}} 
\put(0.4,2.9){\footnotesize $v_{1}$} 
\put(1.2,2.9){\footnotesize $\beta$}
%Chains
\put(1,1){\line(1,1){2}} \put(0,2){\line(1,1){1}}
%Other edges
\put(1,1){\line(-1,1){1}} \put(2,2){\line(-1,1){1}}
\put(3,3){\line(1,-1){1}}
\end{picture}
}
\end{picture}
\hspace*{1.5in}
\setlength{\unitlength}{1.5cm}
\begin{picture}(4,6.5)
\put(-0.25,5.5){$L$}
% edges
\put(1,0){\line(-1,1){1}}
\put(1,0){\line(1,1){2}}
\put(0,1){\line(1,1){2}}
\put(2,1){\line(-1,1){1}}
\put(2,1){\line(0,1){1}}
\put(1,2){\line(0,1){1}}
\put(2,2){\line(-1,1){2}}
\put(2,2){\line(1,1){2}}
\put(3,2){\line(-1,1){1}}
\put(3,2){\line(0,1){1}}
\put(1,3){\line(1,1){2}}
\put(2,3){\line(0,1){1}}
\put(3,3){\line(-1,1){2}}
\put(0,4){\line(1,1){2}}
\put(4,4){\line(-1,1){2}}
% colors
% vertices
\put(2,6){\VertexForLatticeI{0}}
\put(1,5){\VertexForLatticeI{1}}
\put(3,5){\VertexForLatticeI{2}}
\put(0,4){\VertexForLatticeI{3}}
\put(2,4){\VertexForLatticeI{4}}
\put(4,4){\VertexForLatticeI{5}}
\put(1,3){\VertexForLatticeI{6}}
\put(2,3){\VertexForLatticeI{7}}
\put(3,3){\VertexForLatticeI{8}}
\put(1,2){\VertexForLatticeI{9}}
\put(2,2){\VertexForLatticeI{10}}
\put(3,2){\VertexForLatticeI{11}}
\put(0,1){\VertexForLatticeI{12}}
\put(2,1){\VertexForLatticeI{13}}
\put(1,0){\VertexForLatticeI{14}}
\put(1,5){\NEEdgeLabelForLatticeI{\beta}}
\put(3,5){\NWEdgeLabelForLatticeI{\alpha}}
\put(0,4){\NEEdgeLabelForLatticeI{\beta}}
\put(2,4){\NWEdgeLabelForLatticeI{\alpha}}
\put(2,4){\NEEdgeLabelForLatticeI{\beta}}
\put(4,4){\NWEdgeLabelForLatticeI{\beta}}
\put(1,3){\NEEdgeLabelForLatticeI{\beta}}
\put(1,3){\NWEdgeLabelForLatticeI{\alpha}}
\put(2,3){\VerticalEdgeLabelForLatticeI{\alpha}}
\put(3,3){\NWEdgeLabelForLatticeI{\beta}}
\put(3,3){\NEEdgeLabelForLatticeI{\beta}}
\put(1,2){\VerticalEdgeLabelForLatticeI{\alpha}}
\put(1.25,2.25){\NEEdgeLabelForLatticeI{\beta}}
\put(2.2,1.8){\NWEdgeLabelForLatticeI{\beta}}
\put(1.8,1.8){\NEEdgeLabelForLatticeI{\beta}}
\put(3,2){\VerticalEdgeLabelForLatticeI{\alpha}}
\put(2.75,2.25){\NWEdgeLabelForLatticeI{\beta}}
\put(0,1){\NEEdgeLabelForLatticeI{\alpha}}
\put(2,1){\VerticalEdgeLabelForLatticeI{\alpha}}
\put(2,1){\NWEdgeLabelForLatticeI{\beta}}
\put(2,1){\NEEdgeLabelForLatticeI{\beta}}
\put(1,0){\NWEdgeLabelForLatticeI{\beta}}
\put(1,0){\NEEdgeLabelForLatticeI{\alpha}}
\end{picture}

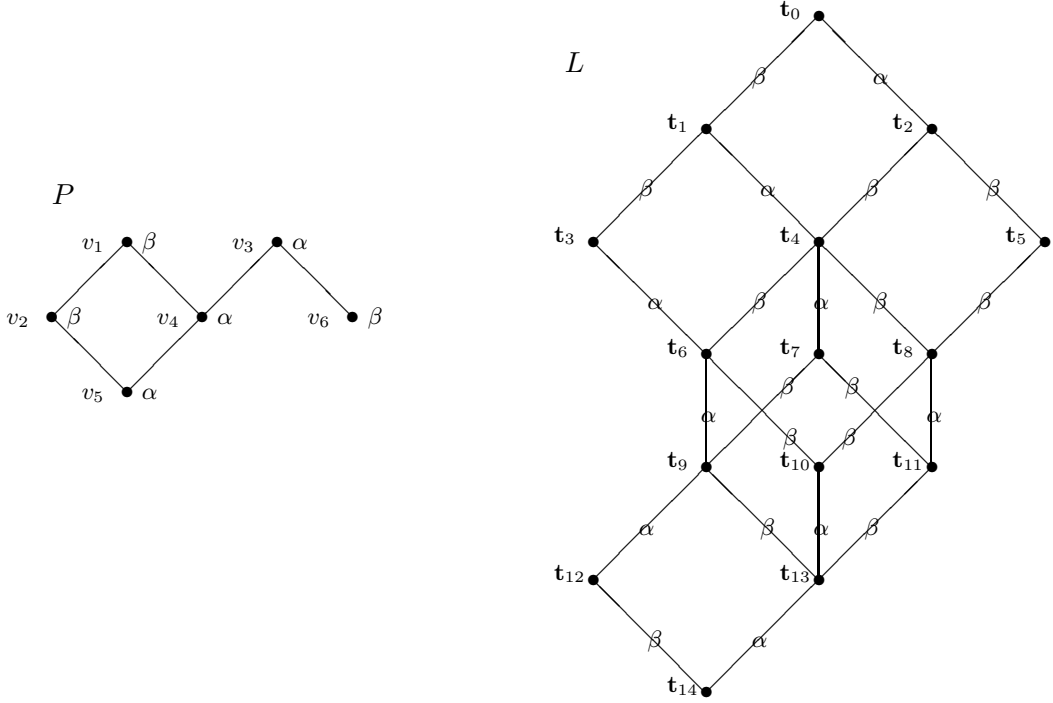
\ComponentFig: The disjoint sum of the $\beta$-components of the edge-colored 
lattice $L$ from \PosetAndLatticeFig. 

\setlength{\unitlength}{1cm}
\hspace*{0.2in}
\begin{picture}(15,3.5)
%Vertex labels
\put(2,2.5){\VertexForLatticeI{0}}
\put(1,1.5){\VertexForLatticeI{1}}
\put(1,1.5){\NEEdgeLabelForLatticeII{\beta}}
\put(0,0.5){\VertexForLatticeI{3}}
\put(0,0.5){\NEEdgeLabelForLatticeII{\beta}}
\put(6,3){\VertexForLatticeI{2}} 
\put(5,2){\VertexForLatticeI{4}}
\put(5,2){\NEEdgeLabelForLatticeII{\beta}}
\put(4,1){\VertexForLatticeI{6}} 
\put(4,1){\NEEdgeLabelForLatticeII{\beta}}
\put(7,2){\VertexForLatticeI{5}}
\put(7,2){\NWEdgeLabelForLatticeII{\beta}}
\put(6,1){\VertexForLatticeI{8}} 
\put(6,1){\NWEdgeLabelForLatticeII{\beta}}
\put(6,1){\NEEdgeLabelForLatticeII{\beta}}
\put(5,0){\VertexForLatticeI{10}}
\put(5,0){\NWEdgeLabelForLatticeII{\beta}}
\put(5,0){\NEEdgeLabelForLatticeII{\beta}}
\put(10,2.5){\VertexForLatticeI{7}}
\put(9,1.5){\VertexForLatticeI{9}}
\put(9,1.5){\NEEdgeLabelForLatticeII{\beta}}
\put(11,1.5){\VertexForLatticeI{11}}
\put(11,1.5){\NWEdgeLabelForLatticeII{\beta}}
\put(10,0.5){\VertexForLatticeI{13}}
\put(10,0.5){\NEEdgeLabelForLatticeII{\beta}}
\put(10,0.5){\NWEdgeLabelForLatticeII{\beta}}
\put(13,2){\VertexForLatticeI{12}}
\put(14,1){\VertexForLatticeI{14}}
\put(14,1){\NWEdgeLabelForLatticeII{\beta}}
%Edges
\put(0,0.5){\line(1,1){2}}
\put(4,1){\line(1,1){2}} 
\put(5,0){\line(1,1){2}}
\put(5,0){\line(-1,1){1}} 
\put(6,1){\line(-1,1){1}}
\put(7,2){\line(-1,1){1}}
\put(9,1.5){\line(1,1){1}} 
\put(10,0.5){\line(1,1){1}}
\put(10,0.5){\line(-1,1){1}} 
\put(11,1.5){\line(-1,1){1}}
\put(14,1){\line(-1,1){1}}
%Sign
\put(2.8,1.375){$\bigoplus$} \put(7.8,1.375){$\bigoplus$}
\put(11.8,1.375){$\bigoplus$}
\end{picture}
\end{center}
%=============================
% Text after Fig 2.1, 2.2
%=============================
vertex $\telt$ and the edge 
$\selt \rightarrow \telt$ 
are {\em above} $\selt$. The vertex $\selt$ is a {\em descendant} of 
$\telt$, and $\telt$ is an {\em ancestor} of $\selt$.  All edge-colored and 
vertex-colored directed graphs in this paper will turn out to be posets.  
For a directed graph $R$, a {\em rank function} is a
surjective function $\rho : R \longrightarrow \{0,\ldots,l\}$
(where $l \geq 0$) with the property that if $\selt \rightarrow
\telt$ in $R$, then $\rho(\selt) + 1 = \rho(\telt)$. If such a rank 
function exists then $R$ is the Hasse diagram for a poset --- a {\em 
ranked} poset.   We call $l$
the {\em length} of $R$ with respect to $\rho$, and the set
$\rho^{-1}(i)$ is the $i${\em th rank} of $R$.  
In an edge-colored ranked
poset $R$, $\mathbf{comp}_{i}(\telt)$ will be a ranked poset for 
each $\telt \in R$ and $i \in I$.  We let
$l_i(\telt)$ denote the length of $\mathbf{comp}_{i}(\telt)$, 
and we let $\rho_i(\telt)$ denote the rank of $\telt$
within this component.  We define the {\em depth} of $\telt$ in its 
$i$-component to be $\delta_{i}(\telt) := l_{i}(\telt) - \rho_{i}(\telt)$.
A ranked poset $R$ with rank function $\rho$ and length $l$ 
is {\em rank symmetric} if
$|\rho^{-1}(i)|=|\rho^{-1}(l-i)|$ for $0 \leq i \leq l$. It 
is {\em rank unimodal} if there is an $m$ such that
$|\rho^{-1}(0)| \leq |\rho^{-1}(1)| \leq \cdots \leq
|\rho^{-1}(m)| \geq |\rho^{-1}(m+1)| \geq \cdots \geq
|\rho^{-1}(l)|$.

The distributive lattice of order ideals of a poset 
$P$, partially ordered by subset containment, will be denoted 
$J(P)$. See \cite{Stanley}. 
A coloring of the vertices of the poset $P$ gives a natural 
coloring of the edges of the distributive lattice 
$L = J(P)$, as follows: Given a function $\vcolor_{P}: 
\mathcal{V}(P) \longrightarrow I$, 
we assign a covering relation $\selt \rightarrow \telt$ in $L$ the 
color $i$ and write $\selt \myarrow{i} \telt$ if 
$\telt \setminus \selt = \{u\}$ and $\vcolor_{P}(u) = i$.  So 
$L$ becomes an edge-colored distributive lattice with edges 
colored by the set $I$; 
we write $L = J_{color}(P)$. 
The edge-colored lattice $L_{G_{2}}(0,1)$ of \FundLatticeIdealsFigure\ is obtained from 
the vertex-colored poset $P_{G_{2}}(0,1)$ of \FundPosets\ in this 
way. 
Note that $J_{color}(P^{*}) \cong  
(J_{color}(P))^{*}$, $J_{color}(P^{\sigma}) \cong  
(J_{color}(P))^{\sigma}$ (recoloring), and $J_{color}(P \oplus Q) 
\cong J_{color}(P) \times J_{color}(Q)$.  
An edge-colored poset $P$ has the {\em diamond coloring property} if 
whenever  
\parbox{1.1cm}{\begin{center}
\setlength{\unitlength}{0.2cm}
\begin{picture}(4,3.5)
\put(2,0){\circle*{0.5}} \put(0,2){\circle*{0.5}}
\put(2,4){\circle*{0.5}} \put(4,2){\circle*{0.5}}
\put(0,2){\line(1,1){2}} \put(2,0){\line(-1,1){2}}
\put(4,2){\line(-1,1){2}} \put(2,0){\line(1,1){2}}
\put(0.75,0.55){\em \small k} \put(2.7,0.7){\em \small l}
\put(0.7,2.7){\em \small i} \put(2.75,2.55){\em \small j}
\end{picture} \end{center}} is an edge-colored subgraph of the Hasse 
diagram for $P$, then $i = l$ and $j = k$.   
A necessary and sufficient condition for an edge-colored distributive 
lattice $L$ to be isomorphic (as an edge-colored poset) to 
$J_{color}(P)$ for some vertex-colored poset $P$ is for $L$ to have 
the diamond coloring property.  
Then for $\selt \in L$ and $i \in I$, 
one can see that $\comp_{i}(\selt)$ is the Hasse 
diagram for a 
distributive lattice. In particular, $\comp_{i}(\selt)$ is a 
distributive sublattice of $L$ in the induced order, and a covering 
relation in $\comp_{i}(\selt)$ is also a covering relation in $L$. 
 
Let $n \geq 1$.  Let $\mathcal{D}$ be a Dynkin diagram with $n$ nodes 
which are indexed by the elements of a set $I$ such that $|I| = n$.  
The associated Cartan matrix is denoted $(\mathcal{D}_{i,j})_{i,j \in 
I}$.  
Throughout this paper $\mathfrak{g}$ will denote the complex 
semisimple Lie algebra of rank $n$
with Chevalley generators $\{x_i,y_i,h_i\}_{i \in I}$ satisfying
the Serre relations specified by the Cartan matrix for the Dynkin 
diagram at hand.  
Usually $I = \{1,\ldots,n\}$. 
In any Cartan matrix, $\mathcal{D}_{i,i} = 2$ for $i \in I$. 
\CartanMatrixEntriesTable\ presents the off-diagonal entries 
$\mathcal{D}_{i,j}$, $i \not= j$, for the rank two semisimple Dynkin 
diagrams $A_{1} \oplus A_{1}$, $A_{2}$, $C_{2}$, and 
$G_{2}$. 
Two Dynkin diagrams $\mathcal{D}$ and $\mathcal{D'}$ are {\em 
isomorphic} if under some one-to-one correspondence $\sigma : I 
\longrightarrow I'$ we have  
$\mathcal{D}_{i,j} = \mathcal{D}'_{\sigma(i),\sigma(j)}$ and 
$\mathcal{D}_{j,i} = \mathcal{D}'_{\sigma(j),\sigma(i)}$. 
Let $E$ denote the Euclidean space equipped with an inner product 
$\langle \cdot,\cdot \rangle$ which contains the root system $\Phi$ 
associated to $\mathcal{D}$.  The set of simple roots is denoted 
$\{\alpha_{i}\}_{i \in I}$. 
For a root $\alpha$, the {\em coroot} is 
$\alpha^{\vee} := \frac{2\alpha}{\langle \alpha,\alpha \rangle}$.  
The $(i,j)$-element $\mathcal{D}_{i,j}$ of the Cartan matrix is  
$\langle \alpha_{i},\alpha_{j}^{\vee} \rangle$.  
The {\em fundamental weights} $\{\omega_{1},\ldots,\omega_{n}\}$ 
form the basis for $E$ dual to the simple coroots 
$\{\alpha_{i}^{\vee}\}_{i=1}^{n}$: $\langle 
\omega_{j},\alpha_{i}^{\vee} \rangle = \delta_{i,j}$.  
The {\em lattice of weights} $\Lambda$ is the set 
of all integral linear combinations of the fundamental weights.  
We coordinatize 
$\Lambda$ to obtain a one-to-one correspondence with 
$\mathbb{Z}^{n}$ as follows: identify $\omega_{i}$  
with the axis vector $(0,\ldots,1,\ldots,0)$, 
where ``1'' is in the $i$th position. 
For $i \in I$, $\alpha_{i} = \sum_{j \in I}\mathcal{D}_{i,j}\omega_{j}$. So
the simple root 
$\alpha_{i}$  can be identified with the $i$th row vector of 
the Cartan matrix.  
The {\em Weyl group} 
$W$ is generated by the {\em simple reflections} $s_{i}: E \rightarrow E$ 
for all $i \in I$:
Here $s_{i}(v) = 
v - \langle v,\alpha_{i}^{\vee} \rangle\alpha_{i}$ for $v \in E$. 

\begin{figure}[ht]
\hspace*{0.65in} \CartanMatrixEntriesTable\ \hspace*{0.35in} 
{\small%
\begin{tabular}{|c||c|c|c|c|}
\hline Subgraph 
& \parbox[c][0.3in][t]{0.6in}{\DynkinZero} 
& \parbox[c][0.3in][t]{0.6in}{\DynkinOne}
& \parbox[c][0.3in][t]{0.6in}{\DynkinTwo}
& \parbox[c][0.3in][t]{0.6in}{\DynkinThree}\\
\hline 
$\mathcal{D}_{i,j}\ \ , \ \ \mathcal{D}_{j,i}$ 
& $0\ \ , \ \ 0$ 
& $-1\ \ , \ \ -1$ 
& $-1\ \ , \ \ -2$ 
& $-1\ \ , \ \ -3$\\
\hline 
\end{tabular}
}
\end{figure}

Vector spaces in this 
paper are complex and 
finite-dimensional. 
If $V$ is a $\mathfrak{g}$-module, then there 
is at least one basis 
$\mathcal{B} := \{v_{\selt}\}_{\selt \in R}$ (where $R$ is an 
indexing set with
$|R| = \dim{V}$) consisting of eigenvectors for the actions of the 
$h_{i}$'s: for any $\selt$ in $R$ and $i \in I$, there exists 
an integer $k_{i}(\selt)$ such that $h_{i}.v_{\selt} = 
k_{i}(\selt)v_{\selt}$.  The {\em weight} of the basis vector 
$v_{\selt}$ is the sum $wt(v_{\selt}) := 
\sum_{i \in I}k_{i}(\selt)\omega_{i}$.  We say $\mathcal{B}$ is a 
{\em weight basis} for $V$.  If $\mu$ is a weight in $\Lambda$, then 
we let $V_{\mu}$ be the subspace of $V$ spanned by all basis vectors 
$v_{\selt} \in \mathcal{B}$ 
such that $wt(v_{\selt}) = \mu$. The subspace $V_{\mu}$ is independent 
of the choice of weight basis $\mathcal{B}$.  
The finite-dimensional irreducible $\mathfrak{g}$-modules are indexed 
by their ``highest weights'' $\lambda$ as these highest weights 
$\lambda$ run through the 
{\em dominant} weights $\Lambda^{+}$ (the 
nonnegative linear combinations of the fundamental weights). 
The Lie algebra $\mathfrak{g}$ acts on the dual space $V^{*}$ by the rule 
$(z.f)(v) = -f(z.v)$ for all $v \in V$, $f \in V^{*}$, and $z \in 
\mathfrak{g}$. 

Let $R$ be a ranked poset whose Hasse diagram edges are colored with 
colors taken from $I$, $|I| = n$.  For $i \in I$, find the connected 
components of the subgraph with edges $\mathcal{E}_{i}(R)$.  For $i \in I$ and 
$\selt$ in $R$, set $m_{i}(\selt) := \rho_{i}(\selt) - \delta_{i}(\selt) = 
2\rho_{i}(\selt) - l_{i}(\selt)$. 
Let $wt_{R}(\selt)$ be the $n$-tuple $(\, m_{i}(\selt)\, )_{i \in 
I}$.  See \FundLatticeTableauxFigure.  
Given a matrix $M = (M_{p,q})_{p,q \in I}$, then for fixed 
$i \in I$ let 
$M^{(i)}$ be the $n$-tuple $(M_{i,j})_{j \in I}$, the $i$th row vector 
for $M$. 
We say $R$ {\em satisfies the structure condition for} $M$ if 
$wt_{R}(\selt) + M^{(i)} = wt_{R}(\telt)$ whenever $\selt 
\myarrow{i} \telt$ for some $i \in I$, that is, 
for all $j \in I$ we have 
$m_{j}(\selt) + M_{i,j} = m_{j}(\telt)$. 
Following 
\cite{DLP1}, we say $R$ {\em satisfies the} $\mathfrak{g}$-{\em structure 
condition} if $M$ is the Cartan matrix for the Dynkin diagram 
$\mathcal{D}$ associated to $\mathfrak{g}$.  
In this case view $wt_{R}: R \longrightarrow 
\Lambda$ as the function given by $wt_{R}(\selt) = 
\sum_{j \in I}m_{j}(\selt)\omega_{j}$.  Then $R$ satisfies the 
$\mathfrak{g}$-structure condition if and only if for each simple 
root $\alpha_{i}$ we have 
$wt_{R}(\selt) + \alpha_{i} = 
wt_{R}(\telt)$ whenever $\selt \myarrow{i} \telt$ in $R$.  (In 
\cite{DonSupp} the edges of $R$ were said to ``preserve 
weights''.)  
This condition requires the color structure of $R$ to be compatible 
with the structure of the set of weights for a representation of 
$\mathfrak{g}$. 
The largest edge-colored distributive 
lattice of \FundLatticeTableauxFigure\ satisfies the structure 
condition for the 
%matrix $M = \left(\begin{array}{cc} 2 & -1\\ -3 & 
%2\end{array}\right)$ 
$G_{2}$ Cartan matrix (\CartanTable) 
and therefore satisfies 
the $G_{2}$-structure condition. 

The following obvious lemma is used 
when the Dynkin diagram has symmetry 
or when other numberings of the Dynkin diagram are 
convenient. 

\noindent
{\bf \DynkinDiagramSymmetryLemma}\ \ {\sl Let $\mathcal{D}$ and 
$\mathcal{D}'$ be Dynkin diagrams with nodes indexed by 
$I$ and $I'$ 
such that $\mathcal{D}$ and $\mathcal{D}'$ 
are isomorphic under a one-to-one correspondence $\sigma: I 
\longrightarrow I'$.  Let $\mathfrak{g}$ and $\mathfrak{g}'$ be 
the respective semisimple Lie algebras.  Let $R$ be a ranked poset 
with edges colored by the set $I$, and consider the recoloring 
$R^{\sigma}$.  Then $R$ satisfies 
the $\mathfrak{g}$-structure condition if and only if $R^{\sigma}$ satisfies 
the $\mathfrak{g}'$-structure condition. 
} 

Let $w_{0}$ be the longest element of the Weyl group $W$ associated 
to $\mathfrak{g}$, as in  
Exercise 10.9 of \cite{Hum}.  
When $w_{0}$ acts on $\Lambda$, then for each $i$ it sends $\alpha_{i} \mapsto 
-\alpha_{\sigma_{0}(i)}$ and $\omega_{i} \mapsto 
-\omega_{\sigma_{0}(i)}$, where  
$\sigma_{0} : I \longrightarrow I$ is some permutation 
of the node labels of the 
Dynkin diagram $\mathcal{D}$. Here $\sigma_{0}$ must be a symmetry 
of the Dynkin diagram, and since $w_{0}^{2} = id$ in $W$ 
it is the case that $\sigma_{0}^{2}$ is 
the identity permutation.  
For any weight $\mu = \sum a_{i}\omega_{i}$ we have 
$-w_{0}\mu = \sum a_{i}\omega_{\sigma_{0}(i)}$.  
For connected Dynkin diagrams, $\sigma_{0}$ 
is trivial except in the cases $A_{n}$ $(n \geq 2)$, 
$D_{2k+1}$ $(k \geq 2)$, and 
$E_{6}$; in these cases it is the only nontrivial Dynkin diagram 
automorphism.  
Given an edge-colored poset $R$ with edges colored by the set $I$ of 
indices for the Dynkin diagram $\mathcal{D}$, 
we let $R^{\triangle}$ be the edge-colored poset 
$(R^{*})^{\sigma_{0}}$ and call $R^{\triangle}$ the $\sigma_{0}$-{\em 
recolored dual} of $R$.  
Observe that $(R^{\triangle})^{\triangle} = R$. 
We allow  
``$\triangle$'' to be applied to any vertex-colored  
poset $Q$ whose vertex colors correspond to nodes 
of a Dynkin diagram. 

The group ring 
$\mathbb{Z}[\Lambda]$ has vector space basis $\{e_{\mu}\, |\, \mu \in \Lambda\}$  
and multiplication rule $e_{\mu + \nu} = e_{\mu}e_{\nu}$.  The Weyl 
group ${W}$ acts on $\mathbb{Z}[\Lambda]$ by the rule 
$\sigma.e_{\mu} := e_{\sigma\mu}$.  
The {\em character ring} 
$\mathbb{Z}[\Lambda]^{W}$ for $\mathfrak{g}$ is the ring of 
$W$-invariant elements of $\mathbb{Z}[\Lambda]$;  
elements of $\mathbb{Z}[\Lambda]^{W}$ are {\em characters for} 
$\mathfrak{g}$.  If $V$ is a 
representation of $\mathfrak{g}$, then the {\em Weyl character for} $V$ is 
$\chi(V) := \sum_{\mu \in \Lambda}\, (\dim V_{\mu})e_{\mu} \in 
\mathbb{Z}[\Lambda]^{W}$. If $V$ is 
irreducible with highest weight 
$\lambda$, let 
$\displaystyle 
\chi_{_{\lambda}} := \chi(V)$. We call $\chi_{_{\lambda}}$ an {\em 
irreducible character}.  
Let $\displaystyle A_{\mu} := \sum_{\sigma \in 
{W}}\det(\sigma)e_{\sigma\mu}$. 
Let $\varrho := \omega_{1} + 
\cdots + \omega_{n}$. 
It is well-known 
that  
$\displaystyle A_{\varrho} = e_{\varrho}\Pi(1-e_{-\alpha})$,  
product taken over the positive roots $\alpha$.  
Weyl's character formula says that 
$\chi_{_{\lambda}}$ is the unique element 
of $\mathbb{Z}[\Lambda]^{W}$ for which 
$A_{\varrho}\chi_{_{\lambda}} = A_{\varrho + \lambda}$.   

Let $V$ be a representation of $\mathfrak{g}$. 
A {\em splitting system for} $V$ (or for $\chi(V)$) is a pair 
$(\mathcal{T},weight)$, where $\mathcal{T}$ is a set and 
$weight:\mathcal{T} \longrightarrow \Lambda$  is a {\em 
weight function} such that 
$\displaystyle \chi(V) := \sum_{\telt \in \mathcal{T}}e_{weight(\telt)}$.  
If $R$ is a ranked 
poset with edges colored by the set $\{1,\ldots,n\}$, if $R$ 
satisfies the structure condition for $\mathfrak{g}$, 
and if $(R,wt_{R})$ is a splitting system for $V$, 
then we say $R$ is a {\em splitting poset for} $V$ (or for $\chi(V)$).  
This concept appears unnamed on p.\ 266 of \cite{DonSupp} and as ``labelling 
poset'' in Corollary 5.3 of \cite{ADLP}. 
An edge-colored ranked poset $R$ for which $(R,wt_{R})$ is a 
splitting system for an irreducible representation can fail to 
satisfy the structure condition for $\mathfrak{g}$. 
We use $z_{i}$ to denote $e_{\omega_{i}}$.   
If $R$ is a splitting poset for $V$, then 
$\displaystyle \chi(V) = \sum_{\telt \in 
R}(z_{1},\ldots,z_{n})^{wt_{R}(\telt)}$, 
where $(z_{1},\ldots,z_{n})^{wt_{R}(\telt)} := 
z_{1}^{m_{1}(\telt)}\cdots{z_{n}^{m_{n}(\telt)}}$.  
Here $\chi_{_{\lambda}}$ is a Laurent polynomial in the 
indeterminates $z_{i}$ with nonnegative integer coefficients. We denote 
this polynomial by 
$char_{\mathfrak{g}}(\lambda; z_{1},\ldots,z_{n})$.  

\noindent
{\bf \DualCharacterLemma}\ \ {\sl Let $V$ be a representation for a 
semisimple Lie algebra $\mathfrak{g}$. Let $\mathfrak{g}'$ be a 
semisimple Lie algebra isomorphic to $\mathfrak{g}$ obtained from an 
isomorphism $\sigma$ of Dynkin diagrams as in the statement of 
\DynkinDiagramSymmetryLemma.  Suppose $R$ is a splitting 
poset for $V$.  Then the edge-colored poset $R^{*}$ is a splitting 
poset for the dual representation $V^{*}$ of $\mathfrak{g}$, 
$R^{\sigma}$ is a splitting 
poset for the $\mathfrak{g}'$-module $V$, and $R^{\triangle}$ is a 
splitting poset for the $\mathfrak{g}$-module $V$.} 

{\em Proof.}  The only assertion that does not 
immediately follow from 
the definitions and \DynkinDiagramSymmetryLemma\ is 
that $R^{\triangle}$ is a splitting poset {\em for the} 
$\mathfrak{g}$-{\em module} $V$. 
Write $V \cong V_{1} \oplus \cdots \oplus V_{k}$, a 
decomposition of $V$ into irreducible $\mathfrak{g}$-modules $V_{i}$ 
such that $V_{i}$ has highest weight $\mu_{i}$.  
The dual $\mathfrak{g}$-module $V^{*}$ has $R^{*}$ as a supporting 
graph; $V^{*}$ decomposes as  
$V_{1}^{*} \oplus \cdots \oplus V_{k}^{*}$, where each $V_{i}^{*}$ is 
irreducible with highest weight $-w_{0}(\mu_{i})$ (cf.\ Exercise 21.6 
of \cite{Hum}).  
Recolor $R^{*}$ by applying the 
permutation $\sigma_{0}$ to obtain $R^{\triangle}$.  Now view 
$V^{*}$ as a new $\mathfrak{g}$-module $U$ induced by the action 
$x_{i}.v := x_{\sigma_{0}(i)}.v$ and $y_{i}.v := y_{\sigma_{0}(i)}.v$ for 
each $i \in I$ and $v \in V^{*}$.  It is apparent that 
$R^{\triangle}$ is a splitting poset for the $\mathfrak{g}$-module $U$. 
Let $U_{i}$ be the (irreducible) 
$\mathfrak{g}$-submodule of $U$ corresponding to 
$V_{i}^{*}$.  One can see that the highest weight of 
$U_{i}$ is now $-w_{0}(-w_{0}(\mu_{i}))$, 
which is just $\mu_{i}$.  Hence $U$ is isomorphic to $V$.\hfill\QED  

\noindent 
{\bf \ConnectedLabellingPosetExists}\ \ 
{\sl Let $V$ be an irreducible $\mathfrak{g}$-module.  Then there is a 
connected splitting poset for $V$.} 

{\em Proof.} By Lemmas 3.1.A and 3.2.A of \cite{DonSupp}, any 
supporting graph 
for $V$ will do.  By Lemma 3.1.F of \cite{DonSupp}, since $V$ is 
irreducible then any  
supporting graph for $V$ is 
necessarily connected.\hfill\QED

This paragraph and \RGFProp\ borrow from Sections 5 and 6 of 
\cite{PrEur}.  
If we set \[x := 2\sum_{i=1}^{n}\left[\sum_{j=1}^{n}\frac{2\langle 
\omega_{i},\omega_{j} \rangle}{\langle \alpha_{j},\alpha_{j} 
\rangle}\right]x_{i},\ \ \  
y := \sum y_{i},\ \ \ \mbox{ and }\ \  h := 
2\sum_{i=1}^{n}\left[\sum_{j=1}^{n}\frac{2\langle 
\omega_{i},\omega_{j} \rangle}{\langle \alpha_{j},\alpha_{j} 
\rangle}\right]h_{i},\] then $\mathfrak{s} := \mathrm{span}\{x,y,h\}$ is 
a three-dimensional subalgebra of 
$\mathfrak{g}$ isomorphic to $\mathfrak{sl}(2,\mathbb{C})$. It is called 
a ``principal three-dimensional subalgebra''.    
Set $\varrho^{\vee} := \sum_{i=1}^{n}\frac{2\omega_{i}}{\langle 
\alpha_{i},\alpha_{i} \rangle}$. 
Observe that $\langle 
\alpha_{i},\varrho^{\vee} \rangle = 1$ for $1 \leq i \leq n$.  
Let $V$ be a $\mathfrak{g}$-module. 
Let $R$ be a splitting poset for $V$. Then there exists a weight basis 
for $V$ which can be indexed by the elements of $R$, say 
$\{v_{\telt}\}_{\telt \in R}$, so that the weight of the basis 
vector $v_{\telt}$ is $wt_{R}(\telt)$.  
One can check that $h.v_{\telt} = 2\langle wt_{R}(\telt),\varrho^{\vee} 
\rangle v_{\telt}$, so the set $\{2\langle 
wt_{R}(\xelt),\varrho^{\vee} \rangle\}_{\xelt\in{R}}$ consists of the 
integral weights for $V$ regarded as an $\mathfrak{s}$-module.  
Choose an element 
$\mathbf{max}$ in $R$ such that $2\langle wt_{R}(\mathbf{max}),\varrho^{\vee} 
\rangle$ is largest in the set $\{2\langle 
wt_{R}(\xelt),\varrho^{\vee} \rangle\}_{\xelt\in{R}}$, and choose 
$\mathbf{min}$ such that $2\langle wt_{R}(\mathbf{min}),\varrho^{\vee} 
\rangle$ is smallest. Symmetry of the integral weights for $V$ under 
the action of 
$\mathfrak{s} \cong \mathfrak{sl}(2,\mathbb{C})$ 
implies that 
$2\langle wt_{R}(\mathbf{max}),\varrho^{\vee} 
\rangle = -2\langle wt_{R}(\mathbf{min}),\varrho^{\vee} 
\rangle$.   Set $l := 2\langle 
wt_{R}(\mathbf{max}),\varrho^{\vee} \rangle$.  
Since $R$ satisfies the $\mathfrak{g}$-structure condition, 
it follows that if $\selt \myarrow{i} \telt$ is an edge in $R$, 
then $wt_{R}(\selt) + 
\alpha_{i} = wt_{R}(\telt)$; therefore $\langle 
wt_{R}(\selt),\varrho^{\vee} \rangle + 1 = \langle 
wt_{R}(\telt),\varrho^{\vee} \rangle$. 
Suppose for the moment that $R$ is connected.  
Then the weights $\{2\langle 
wt_{R}(\xelt),\varrho^{\vee} \rangle\}_{\xelt\in{R}}$ all have the 
same parity. 
Consider the function $\rho: R 
\longrightarrow \mathbb{Z}$ given by $\rho(\telt) := 
\frac{l}{2}+\langle wt_{R}(\telt),\varrho^{\vee} \rangle$.  Based on 
what we have seen so far, the range of $\rho$ is the set of integers 
$\{0,\ldots,l\}$, and hence $\rho$ is the rank function for $R$.  
Next consider the case that $V$ is irreducible with highest weight 
$\lambda$.  Then $R$ need not be connected.  However, since $V$ has a 
connected splitting poset by \ConnectedLabellingPosetExists, 
then the weights $\{2\langle 
wt_{R}(\xelt),\varrho^{\vee} \rangle\}_{\xelt\in{R}}$ all have the 
same parity. Thus the function $\rho: R 
\longrightarrow \mathbb{Z}$ given by $\rho(\telt) := 
\frac{l}{2}+\langle wt_{R}(\telt),\varrho^{\vee} \rangle$ will be a 
rank function for $R$ with range $\{0,\ldots,l\}$. Call $\rho$ the 
{\em natural rank function} for $R$.  Since $V$ is irreducible, we 
can see that  
$\mathbf{max}$ is the unique element of $R$ with weight 
$wt_{R}(\mathbf{max}) = \lambda$. Hence $l = 2\langle 
\lambda,\varrho^{\vee} \rangle$.  
Next we define the {\em rank generating function} for $R$ to be 
$\mathit{RGF}_{\mathfrak{g}}(\lambda, q) := 
\sum_{i=0}^{l}|\rho^{-1}(i)|q^{i} = \sum_{\telt \in R} 
q^{\rho(\telt)}$.  This is the usual rank generating function for the 
ranked poset $R$. 
We do not refer to $R$ in the notation 
$\mathit{RGF}_{\mathfrak{g}}(\lambda, q)$ because: 
We have 
$\Big|\rho^{-1}(i)\Big| = \Big|\{\telt \in R\, |\, \frac{l}{2} + 
\langle wt_{R}(\telt),\varrho^{\vee} \rangle = i\}\Big| = 
\sum_{\mu}\dim(V_{\mu})$, where the latter sum is over all 
weights $\mu$ such that $\frac{l}{2}+\langle \mu,\varrho^{\vee} 
\rangle = i$.  
Thus if $R'$ is another 
naturally-ranked splitting poset for $V$, then corresponding ranks of 
$R$ and $R'$ have the same size.  
To obtain the rank generating function identity in the following 
result we use the ``principal specialization'' of Weyl's character 
formula from Section 6 of \cite{PrEur}. 

\noindent 
{\bf \RGFProp}\ \ {\sl Let $V$ be an irreducible $\mathfrak{g}$-module 
with highest weight $\lambda$, and let $R$ be a splitting poset for 
$V$ with the natural rank function identified in the preceding paragraph. 
(If $R$ is connected, then the natural rank function is the unique 
rank function.) Then $R$ is rank symmetric and rank unimodal, and}
\[\mathit{RGF}_{\mathfrak{g}}(\lambda, q) = \frac{\mbox{$\displaystyle 
\Pi_{\alpha \in 
\Phi^{+}}$}(1-q^{\langle \lambda+\varrho,\alpha \rangle})}{\Pi_{\alpha \in 
\Phi^{+}}(1-q^{\langle \varrho,\alpha \rangle})}\]

{\em Proof.} Choose a connected splitting poset $R'$ for $V$; the 
natural rank function for $R'$ is the unique rank function.  Then 
by Proposition 3.5 of \cite{DonSupp}, it follows that $R'$ is rank 
symmetric and rank unimodal.   From the observation of the 
next-to-last sentence of the paragraph preceding the proposition, 
we conclude that  
the naturally ranked poset $R$ is rank symmetric and rank unimodal.  
The principal specialization obtained 
from \cite{PrEur} pp.\ 337-338 is for simple Lie algebras, but the 
same arguments are valid for semisimple Lie algebras.  Apply this to 
obtain the rank generating function identity of the proposition 
statement.\hfill\QED

%==================================================================
%\newpage
\vspace{2ex} 

\noindent
{\Large \bf \GridNum.\ \ Grid posets and \dichromatic grid posets}

Here we introduce general grid posets and two-color grid posets with 
purely combinatorial definitions.  From Section \SemiNum\ onward we will consider 
only the particular two-color grid posets called  ``$\mathfrak{g}$-semistandard'' 
posets, whose structures are indexed by rank two Dynkin diagrams.  Some 
(uncolored) grid posets are displayed in \ThreeElementGridPosetFigure;  
the poset  $P$  in 
\PosetAndLatticeFig\ is a two-colored grid poset.  In the general setting of this 
section, \WeightsLemma\ and its related definitions provide for the 
decomposition of two-color grid posets into manageable pieces. 
Given $m \geq 1$, set $[m] := \{1,2,\ldots,m\}$. 

\setcounter{myfn}{1}

Given a finite poset $(P,\leq_{_{P}})$,  
a {\em chain function for} $P$ is a function  
$\mathbf{chain}: P 
\longrightarrow [m]$ for some positive integer $m$ 
such that (1) $\mathbf{chain}^{-1}(i)$ is a (possibly empty) 
chain in $P$ for $1 \leq 
i \leq m$, and (2) given any cover  
$u \rightarrow v$ in $P$, it is the case that either 
$\mathbf{chain}(u) = \mathbf{chain}(v)$ or $\mathbf{chain}(u) = 
\mathbf{chain}(v) + 1$.  
A {\em grid poset} 
is a finite poset $(P,\leq_{_{P}})$ 
together with a chain function  
$\mathbf{chain}: P \longrightarrow [m]$ 
for some $m \geq 1$.  Depending on context, 
the notation $P$ can refer to 
the grid poset $(P,\leq_{_{P}}, 
\mathbf{chain}: P \longrightarrow [m])$ 
or the underlying poset $(P,\leq_{_{P}})$. 
The conditions on $\mathbf{chain}$ imply that an element in a grid 
poset covers no more than two elements and is covered by no more than 
two elements.\myfootnote{Motivation for terminology: 
For $m, n \geq 1$, let $\mathcal{G}$ be 
the directed graph with  
$\mathcal{V}(\mathcal{G}) = \{(p,q) \in \mathbb{Z} \times \mathbb{Z} | 
1 \leq p \leq n, 1 \leq q \leq m\}$ and with 
$\mathcal{E}(\mathcal{G}) = \{(p,q) \rightarrow (r,s) \mbox{ if } 
(r,s) - (p,q) = (1,0) \mbox{ or } (0,1)\}$.  Refer to 
$\mathcal{G}$ as a ``directed grid graph''.  Here 
$\mathcal{G}$ is the Hasse diagram for a poset obtained 
by rotating the plane 
counterclockwise through an angle of $45^{\circ}$, so that the  
vertex $(1,1)$ of $\mathcal{G}$ 
is the minimal element.  A grid poset $(P,\leq_{_{P}}, 
\mathbf{chain}: P \longrightarrow [m])$ can be 
obtained as a subgraph of a directed grid graph for an appropriately 
large $n$ by removing some vertices and some ``NW'' edges.}   
Observe that if $i$ is the smallest (respectively largest) integer such that 
$\mathbf{chain}^{-1}(i)$ is nonempty and   
if $u$ is the maximal (respectively minimal) element of 
$\mathbf{chain}^{-1}(i)$, then $u$ is a maximal 
(respectively minimal) element of the poset $P$.  
%======
% Needed only for Mars stuff
%If $P$ is nonempty, the {\em rightmost maximal 
%element} is the maximal element $z$ with the property that 
%$\mathbf{chain}(z) \geq \mathbf{chain}(u)$ for all maximal elements 
%$u$ in $P$. A {\em face vertex} of the grid poset 
%$P$ is an element $u$ of $P$  
%such that whenever $u \rightarrow v$, then 
%$\mathbf{chain}(u) = \mathbf{chain}(v)$. 
%======
A grid poset $P$ is 
{\em connected} if and only if the Hasse diagram for the poset $P$ is 
connected.  For $1 \leq i \leq m$ we set 
$\mathcal{C}_{i} := \mathbf{chain}^{-1}(i)$. When we depict grid 
posets, the chains $\mathcal{C}_{i}$ will be directed from SW to NE.  
See \ThreeElementGridPosetFigure. 

\vspace*{-0.15in}
\begin{center}
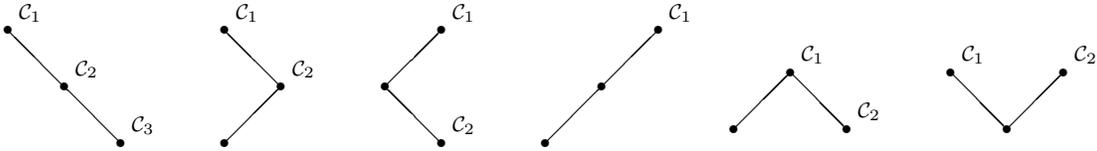
\ThreeElementGridPosetFigure: The six non-isomorphic connected 
grid posets with three elements. \\ 
%=========================================
\setlength{\unitlength}{0.75cm}
\begin{picture}(3.5,3.5)
\put(0,0.25){
\begin{picture}(3,3)
%Vertices
\put(2,0){\circle*{0.15}} 
\put(1,1){\circle*{0.15}} 
\put(0,2){\circle*{0.15}} 
%LInes
\put(2,0){\line(-1,1){2}} 
%Labels
\put(2.2,0.2){\scriptsize $\mathcal{C}_{3}$}
\put(1.2,1.2){\scriptsize $\mathcal{C}_{2}$}
\put(0.2,2.2){\scriptsize $\mathcal{C}_{1}$}
\end{picture}
}
\end{picture}
%=========================================
\setlength{\unitlength}{0.75cm}
\begin{picture}(2.5,3.5)
\put(0,0.25){
\begin{picture}(2.5,3)
%Vertices
\put(0,0){\circle*{0.15}} 
\put(1,1){\circle*{0.15}} 
\put(0,2){\circle*{0.15}} 
%Lines
\put(0,0){\line(1,1){1}}
\put(1,1){\line(-1,1){1}} 
%Labels
\put(1.2,1.2){\scriptsize $\mathcal{C}_{2}$}
\put(0.2,2.2){\scriptsize $\mathcal{C}_{1}$}
\end{picture}
}
\end{picture}
%=========================================
\setlength{\unitlength}{0.75cm}
\begin{picture}(2.5,3.5)
\put(0,0.25){
\begin{picture}(2.5,3)
%Vertices
\put(1,0){\circle*{0.15}} 
\put(0,1){\circle*{0.15}} 
\put(1,2){\circle*{0.15}} 
%Lines
\put(1,0){\line(-1,1){1}}
\put(0,1){\line(1,1){1}} 
%Labels
\put(1.2,0.2){\scriptsize $\mathcal{C}_{2}$}
\put(1.2,2.2){\scriptsize $\mathcal{C}_{1}$}
\end{picture}
}
\end{picture}
%=========================================
\setlength{\unitlength}{0.75cm}
\begin{picture}(3,3.5)
\put(0,0.25){
\begin{picture}(3,3)
%Vertices
\put(0,0){\circle*{0.15}} 
\put(1,1){\circle*{0.15}} 
\put(2,2){\circle*{0.15}} 
%Lines
\put(0,0){\line(1,1){2}}
%Labels
\put(2.2,2.2){\scriptsize $\mathcal{C}_{1}$}
\end{picture}
}
\end{picture}
%=========================================
\setlength{\unitlength}{0.75cm}
\begin{picture}(3.5,3.5)
\put(0,0.5){
\begin{picture}(3,3)
%Vertices
\put(0,0){\circle*{0.15}} 
\put(1,1){\circle*{0.15}} 
\put(2,0){\circle*{0.15}} 
%Lines
\put(0,0){\line(1,1){1}}
\put(1,1){\line(1,-1){1}}
%Labels
\put(2.2,0.2){\scriptsize $\mathcal{C}_{2}$}
\put(1.2,1.2){\scriptsize $\mathcal{C}_{1}$}
\end{picture}
}
\end{picture}
%=========================================
\setlength{\unitlength}{0.75cm}
\begin{picture}(2.5,3.5)
\put(0,0.5){
\begin{picture}(3,3)
%Vertices
\put(0,1){\circle*{0.15}} 
\put(1,0){\circle*{0.15}} 
\put(2,1){\circle*{0.15}} 
%Lines
\put(0,1){\line(1,-1){1}}
\put(1,0){\line(1,1){1}}
%Labels
\put(2.2,1.2){\scriptsize $\mathcal{C}_{2}$}
\put(0.2,1.2){\scriptsize $\mathcal{C}_{1}$}
\end{picture}
}
\end{picture}
\end{center}
%\end{figure}

\vspace*{-0.15in}
Let $(P,\leq_{_{P}}, 
\mathbf{chain}: P \longrightarrow [m])$ be a grid 
poset.  
The {\em dual grid poset} $P^{*}$ is the dual poset 
$P^{*}$ together with the chain function 
$\mathbf{chain}^{*}: P^{*} \longrightarrow [m]$ given by 
$\mathbf{chain}^{*}(u^{*}) = m+1-\mathbf{chain}(u)$ for all $u \in P$. 
For $i = 1,2$, let $P_{i}$ 
be a grid poset with chain function $\mathbf{chain}_{i}: P_{i} 
\longrightarrow [m_{i}]$ for some $m_{i} \geq 1$.  
A one-to-one correspondence $\phi: P_{1} \longrightarrow P_{2}$ is an 
{\em isomorphism of grid posets} if we have $u \rightarrow v$ in 
$P_{1}$ with $\mathbf{chain}_{1}(u) = \mathbf{chain}_{1}(v)$ (respectively 
$\mathbf{chain}_{1}(u) = \mathbf{chain}_{1}(v) + 1$) if and only if  
$\phi(u) \rightarrow \phi(v)$ in 
$P_{2}$ with $\mathbf{chain}_{2}(\phi(u)) = \mathbf{chain}_{2}(\phi(v))$ 
(respectively 
$\mathbf{chain}_{2}(\phi(u)) = \mathbf{chain}_{2}(\phi(v)) + 1$).  
\ThreeElementGridPosetFigure\ depicts each of the isomorphism classes 
of connected grid posets with three elements apiece. 
Given a nonempty grid poset 
$(P,\leq_{_{P}},\mathbf{chain}: P \longrightarrow [m])$, 
there exists some $m' \geq 1$ and a 
surjective chain function $\mathbf{chain}':P \longrightarrow [m']$ 
such that the grid poset $P$ is isomorphic to $(P, {\leq_{_{P}},} 
\mathbf{chain'}: P \longrightarrow [m'])$. If $P$ is 
connected, then this surjective chain function $\mathbf{chain}'$ is 
unique. We say $Q$ is a {\em grid subposet} of a given grid poset 
$P$ if (1) 
$Q$ is a subposet of $P$ in the induced order, and (2) whenever 
$u \rightarrow v$ is a covering relation in $Q$ then it is also a 
covering relation in $P$.  In this case, we regard $Q$ with the chain 
function $\mathbf{chain}|_{Q}$ to be a grid poset on its own.  

For a grid poset $(P,\leq_{_{P}}, 
\mathbf{chain}: P \longrightarrow [m])$, 
let $\mathcal{T}_{P}$ be the totally ordered set whose elements are 
the elements of $P$ and whose ordering is given by the following rule: 
for distinct $u$ and $v$ in $P$ write 
$u <_{_{\mathcal{T}_{P}}} v$ 
if and only if (1) $\mathbf{chain}(u) < \mathbf{chain}(v)$ or 
(2) $\mathbf{chain}(u) = \mathbf{chain}(v)$ with  
$v <_{_{P}} u$. 
Let $l := |P|$. Number the vertices of $P$ $v_{1}, 
v_{2},\ldots, v_{l}$ so that $v_{p} 
<_{_{\mathcal{T}_{P}}} v_{q}$ 
whenever $1 \leq p < q \leq l$. Let $L := J(P)$ 
be the distributive lattice of order ideals of $P$.  We 
simultaneously think of order ideals of $P$ as subsets of $P$ 
and as elements of $L$.  

A {\em two-color function} for a grid poset $(P,\leq_{_{P}}, 
\mathbf{chain}: P \longrightarrow [m])$ is a function 
$\mathbf{color}: P \longrightarrow \Delta$ such that (1) $|\Delta| = 
2$, (2) $\mathbf{color}(u) = \color(v)$ 
if $\mathbf{chain}(u) = \mathbf{chain}(v)$, and (3) 
if $u$ and $v$ are in the same connected component 
of $P$ with $\mathbf{chain}(u) = \mathbf{chain}(v)+1$, 
then $\color(u) \not= \color(v)$. A {\em \dichromatic grid poset} is a 
grid poset 
$(P, \leq_{_{P}}, \mathbf{chain}: P \longrightarrow [m])$ 
together with a two-color function 
$\color: P \longrightarrow \Delta$. In some contexts we will use the 
notation $P$ to refer to the \dichromatic grid poset 
$\digriddelta$. 
Two-color grid posets are 
vertex-colored posets. 
To a  
\dichromatic grid poset $P$ we associate 
the edge-colored distributive lattice $L := 
J_{color}(P)$, as in Section \DefsNum.  
The number of nonempty chains $\mathcal{C}_{i}$ in $P$ of color 
$\gamma \in \Delta$ gives an upper bound for the number of ancestors 
(respectively, descendants) an element in $L$ can have along edges of 
color $\gamma$. One can also see that any color $\gamma$ 
component of $L$ is poset-isomorphic to a product of chains. 
The {\em dual \dichromatic grid poset} $P^{*}$ is the dual 
grid poset $P^{*}$ together with the two-color function 
$\mathbf{color}^{*}: P^{*} \longrightarrow \Delta$ given by 
$\mathbf{color}^{*}(u^{*}) = \mathbf{color}(u)$ for all $u \in P$. 
If $Q$ is 
a grid subposet of the \dichromatic grid poset $P$, then $Q$ is a 
\dichromatic grid poset with chain function $\mathbf{chain}|_{Q}$ and 
two-color 
function $\mathbf{color}|_{Q}$. In this case we call $Q$ a {\em 
\dichromatic grid subposet of} $P$.  
Two \dichromatic grid posets 
$(P_{i},\leq_{_{P_{i}}}, 
\mathbf{chain}_{i}: P_{i} \longrightarrow [m_{i}], 
\mathbf{color}_{i}: P_{i} \longrightarrow \Delta)$ for $i = 1,2$ 
are {\em isomorphic} if there is an isomorphism $\phi: P_{1}  
\longrightarrow P_{2}$ of grid posets such that 
$\mathbf{color}_{2}(\phi(u)) = \mathbf{color}_{1}(u)$ for all $u$ in 
$P_{1}$.  
%=====================================================================
%A ``weak'' form of \dichromatic grid poset isomorphism that doesn't 
%require that the target sets be the same. 
%$\mathbf{color}_{2}(\phi(u)) = 
%\mathbf{color}_{2}(\phi(v))$ if and only if $\mathbf{color}_{1}(u) = 
%\mathbf{color}_{1}(v)$.  
%=====================================================================
We will often  
take $\Delta := \{\alpha,\beta\}$.  
When we {\em switch} (or {\em reverse}) the vertex 
colors of $P$ we replace the color function 
$\mathbf{color}: P \longrightarrow \{\alpha, \beta\}$ with the 
color function $\mathbf{color}': P \longrightarrow \{\alpha, \beta\}$ given 
by: $\mathbf{color}'(v) = \alpha$ if $\mathbf{color}(v) = \beta$, and 
$\mathbf{color}'(v) = \beta$ if $\mathbf{color}(v) = \alpha$.  
Similarly, one can {\em switch} (or {\em reverse}) the edge colors of 
$L$. 
In \GridPosetsFigureList\ 
we depict eight \dichromatic grid 
posets; the numbering of the vertices for each 
poset $P$ follows the 
total ordering $\mathcal{T}_{P}$. The vertex-colored poset $P$ of 
\PosetAndLatticeFig\ is a two-color grid poset.  
The lattice $L$ in that figure is $J_{color}(P)$. 

In this paper the following definition is needed only for a comment 
in Section \SemiNum\ and for preview statements of 
\ClassificationTheorems.  (It is also needed in \cite{ADLP}.)  
We say a \dichromatic grid poset $P$  
has the {\em max property} if $P$ is isomorphic to a \dichromatic grid 
poset $(Q,\leq_{_{Q}}, 
\mathbf{chain}: Q \longrightarrow [m], 
\mathbf{color}: Q \longrightarrow \Delta)$ with a surjective chain 
function such that (1) if $u$ is any maximal element in 
the poset $Q$, then $\mathbf{chain}(u) \leq 2$, and (2) if $v \not= u$ is 
another maximal element in $Q$, then $\mathbf{color}(u) \not= 
\mathbf{color}(v)$. 
Note that the dual \dichromatic grid poset $P^{*}$ might fail 
to have the max property.  
The \dichromatic grid posets of \MaxPropertyFigList\ have the max 
property. 

Let  $P$ be a grid poset with chain function $\mathbf{chain}: P 
\longrightarrow [m]$.  Suppose $P_{1}$ is 
a nonempty order ideal such that $P_{1} \not= P$.  
Regard $P_{1}$ and $P_{2} := P \setminus P_{1}$ to be subposets of 
the poset $P$ in the induced order.  Suppose that whenever $u$ is a maximal 
(respectively minimal) element of $P_{1}$ and $v$ is a maximal 
(respectively minimal) element 
of $P_{2}$, then $\mathbf{chain}(u) \leq \mathbf{chain}(v)$.  Then we 
say that $P$ 
{\em decomposes into} $P_{1} \triangleleft P_{2}$, and we write 
$P = P_{1} \triangleleft P_{2}$. 
If no such order ideal $P_{1}$ exists, 
then we say 
the grid poset $P$ is {\em indecomposable}. See 
\NotLieFundPosets.  
Note that if $P = P_{1} \triangleleft P_{2}$ and 
$u < v$ in $P$ with $u \in P_{2}$, then $v 
\in P_{2}$. 
Moreover, if $u \rightarrow v$ in $P$ with $u \in P_{1}$ and $v \in 
P_{2}$, then $\mathbf{chain}(u) = \mathbf{chain}(v)$. 
Also, if $u \rightarrow v$ is a covering 
relation in the poset $P_{i}$ for $i \in \{1,2\}$, note that $u \rightarrow 
v$ is also a covering relation in $P$.  Hence  
each $P_{i}$ is a grid subposet of $P$.  
If $P$ is a grid poset that decomposes into $P_{1} \triangleleft Q$, 
and if $Q$ decomposes into $P_{2} \triangleleft P_{3}$, then  
$P = P_{1} \triangleleft (P_{2} \triangleleft P_{3})$.  But now 
observe that $P = (P_{1} \triangleleft P_{2}) \triangleleft P_{3}$.  
So we may write $P = P_{1} \triangleleft P_{2} \triangleleft P_{3}$.  
In general, if 
$P = P_{1} \triangleleft P_{2} \triangleleft \cdots \triangleleft 
P_{k}$, then each $P_{i}$ with chain function 
$\mathbf{chain}|_{P_{i}}$ is a grid subposet of $P$. Also,   
an order ideal $\selt$ of $P$ may be expressed as 
the disjoint union $(\selt \cap P_{1}) \cup (\selt \cap P_{2}) 
\cup \cdots \cup (\selt \cap P_{k})$, where each $\selt \cap P_{i}$ is 
an order ideal in $P_{i}$. 
If in addition $P = P_{1} \triangleleft P_{2} \triangleleft \cdots \triangleleft 
P_{k}$ is a \dichromatic grid poset with two-color 
function $\color$, then each $P_{i}$ with chain function 
$\mathbf{chain}|_{P_{i}}$ and two-color function 
$\color|_{P_{i}}$ is a \dichromatic grid subposet of $P$. Here 
$P_{1} \triangleleft P_{2} \triangleleft 
\cdots \triangleleft P_{k}$ is a 
decomposition of $P$ into \dichromatic grid posets. 

Consider a \dichromatic grid poset $\digrid$ with edge-colored distributive 
lattice $L = J_{color}(P)$. For each $\selt$ in $L$ 
we can view the quantity $wt_{L}(\selt)$ 
as the pair $(2\rho_{\alpha}(\selt) - 
l_{\alpha}(\selt)\, ,\, 2\rho_{\beta}(\selt) - 
l_{\beta}(\selt))$ in $\mathbb{Z} \times \mathbb{Z}$.  
The mapping $wt_{L}: L \longrightarrow \mathbb{Z} \times \mathbb{Z}$ 
is the {\em lattice weight 
function} for $L$.  If 
$P = P_{1} \triangleleft P_{2} \triangleleft \cdots \triangleleft 
P_{k}$, then for each $i$ we let 
$L_{i} := J_{color}(P_{i})$ be the edge-colored lattice for the \dichromatic 
grid subposet $P_{i}$ of $P$.  Then 
$wt_{L_{i}}$ denotes the lattice weight function for $L_{i}$.   
We let $\rho^{(i)}_{\alpha}$ and 
$l^{(i)}_{\alpha}$ (respectively $\rho^{(i)}_{\beta}$ and 
$l^{(i)}_{\beta}$) denote the rank and length functions for color 
$\alpha$ (respectively color $\beta$) for $L_{i}$.  

\noindent 
{\bf \WeightsLemma}\ \ {\sl Let $\digrid$ be a \dichromatic grid 
poset, and suppose $P$ decomposes into 
$P = P_{1} \triangleleft P_{2} \triangleleft \cdots \triangleleft 
P_{k}$. Keep the notation of the preceding paragraph.  (1) 
Let $\gamma \in \Delta = \{\alpha,\beta\}$, and 
let $\selt$ be an element of $L = J_{color}(P)$.  
Then } \[\rho_{\gamma}(\selt) = \sum_{i=1}^{k}  
\rho^{(i)}_{\gamma}(\selt \cap P_{i}), \hspace*{0.15in} 
l_{\gamma}(\selt) = \sum_{i=1}^{k} 
l^{(i)}_{\gamma}(\selt \cap P_{i}), \hspace*{0.15in} \mbox{and}  
\hspace*{0.15in} 
wt_{L}(\selt) = \sum_{i=1}^{k} wt_{L_{i}}(\selt \cap 
P_{i}).\] {\sl (2) Consequently, if 
there is a $2 \times 2$ matrix $M = (M_{\iota,\kappa})_{(\iota,\kappa) 
\in \Delta \times \Delta}$ such that 
each edge-colored distributive lattice 
$L_{i} = J_{color}(P_{i})$ satisfies the structure condition for 
$M$, then $L$ satisfies the structure condition for 
$M$ as well.} 

{\em Proof.} 
First we show how (2) follows from (1).  
Given an edge $\selt \myarrow{\gamma} \telt$ in $L$, 
then it is the case that for some $j$ with $1 \leq j \leq k$ we 
have $\selt \cap P_{j} \myarrow{\gamma} \telt \cap P_{j}$ in $L_{j}$ 
while for 
$1 \leq i \leq k$ with $i \not= j$ we have 
$\selt \cap P_{i} = \telt \cap P_{i}$. 
Since $L_{j}$ satisfies the structure condition, we
see that 
$wt_{L_{j}}(\selt \cap P_{j}) + M^{(\gamma)} = 
wt_{L_{j}}(\telt \cap P_{j})$. For $i \not= j$ we have 
$wt_{L_{i}}(\selt \cap P_{i}) = wt_{L_{i}}(\telt \cap P_{i})$.   By  
(1) it follows that $wt_{L}(\selt) + M^{(\gamma)} = wt_{L}(\telt)$.

The results in (1) for general $k$ follow by induction 
once we prove the results for $k = 2$.  So let $k = 2$, 
$\selt \in L$, and $\gamma \in \{\alpha, \beta\}$. It 
suffices to show that  
$\rho_{\gamma}(\selt) = \rho^{(1)}_{\gamma}(\selt \cap P_{1}) + 
\rho^{(2)}_{\gamma}(\selt \cap P_{2})$ and 
$l_{\gamma}(\selt) = 
l^{(1)}_{\gamma}(\selt \cap P_{1}) + l^{(2)}_{\gamma}(\selt \cap 
P_{2})$. 
Let $\relt_{0}, \relt_{1}, \ldots$ be the sequence with 
$\relt_{0} := \selt$ and 
$\relt_{j+1} := \relt_{j} \setminus \{v_{i_{j+1}}\}$, 
where $j \geq 0$ and 
$v_{i_{j+1}}$ is the smallest vertex in $\mathcal{T}_{P}$ of color 
$\gamma$ that can be removed from $\relt_{j}$ so 
that $\relt_{j+1}$ is 
an order ideal of $P$.  Let $\relt_{q}$ be the terminal element of 
the sequence. Observe that $i_{1} < i_{2} < 
\cdots < i_{q}$. We have $\relt_{q} \myarrow{\gamma} 
\relt_{q-1} \myarrow{\gamma} \cdots \myarrow{\gamma} \relt_{1} 
\myarrow{\gamma} \relt_{0} = \selt$. Similarly define a sequence 
$\uelt_{0}, \uelt_{1}, \ldots$ where $\uelt_{0} := \selt$ and 
$\uelt_{s+1} := \uelt_{s} \cup \{v_{r_{s+1}}\}$ where $s \geq 0$ and 
$v_{r_{s+1}}$ is the largest element in $\mathcal{T}_{P}$ of color 
$\gamma$ not in 
$\uelt_{s}$ that can be added to $\uelt_{s}$ so that $\uelt_{s+1}$ is 
an order ideal of $P$.  Let $\uelt_{p}$ be the terminal element of 
the sequence.  Observe that $r_{1} > r_{2} > \cdots > r_{p}$.  
We have $\selt = \uelt_{0} 
\myarrow{\gamma} \cdots \myarrow{\gamma} 
\uelt_{p-1} \myarrow{\gamma} \uelt_{p}$.  Since 
$\comp_{\gamma}(\selt)$ is the Hasse diagram for a distributive 
lattice, \ 
and since $\relt_{q}$ and $\uelt_{p}$ are respectively 
%============================
% Begin pagebreak for figures
%============================

%\begin{figure}[h,t]
\newpage
\begin{center}
\GridPosets: Depicted below are four \dichromatic grid posets each 
possessing the max 
property.\\  
{\small (Each is the $\mathfrak{g}$-semistandard poset 
$P_{\mathfrak{g}}^{\beta\alpha}(2,2)$ of \S \SemiNum\  
for the 
indicated rank two semisimple Lie algebra $\mathfrak{g}$.)} 

\setlength{\unitlength}{1cm}
\begin{picture}(6,7)
\put(0,7){\fbox{$\mathfrak{g} = A_{1} \oplus A_{1}$}}
\put(0.5,2.75){
\begin{picture}(3,3.5)
%Vertex labels
\put(3,2){\circle*{0.15}} 
\put(2.4,1.9){\footnotesize $v_{4}$}
\put(3.2,1.9){\footnotesize $\alpha$}
\put(4,3){\circle*{0.15}} 
\put(3.4,2.9){\footnotesize $v_{3}$}
\put(4.2,2.9){\footnotesize $\alpha$}
\put(0,2){\circle*{0.15}} 
\put(-0.6,1.9){\footnotesize $v_{2}$}
\put(0.2,1.9){\footnotesize $\beta$}
\put(1,3){\circle*{0.15}} 
\put(0.4,2.9){\footnotesize $v_{1}$}
\put(1.2,2.9){\footnotesize $\beta$}
%Chains
\put(3,2){\line(1,1){1}}
\put(0,2){\line(1,1){1}} 
%Chain labels
\put(1.5,3.5){$\mathcal{C}_{1}$}
\put(4.5,3.5){$\mathcal{C}_{2}$}
\put(1.8,2.375){$\bigoplus$}
\end{picture}
}
\put(0,2){\fbox{$\mathfrak{g} = A_{2}$}}
\put(0,-2.5){
\begin{picture}(3,3.5)
%Vertex labels
\put(4,2){\circle*{0.15}} 
\put(3.4,1.9){\footnotesize $v_{8}$}
\put(4.2,1.9){\footnotesize $\beta$}
\put(5,3){\circle*{0.15}} 
\put(4.4,2.9){\footnotesize $v_{7}$}
\put(5.2,2.9){\footnotesize $\beta$}
\put(1,1){\circle*{0.15}} 
\put(0.4,0.9){\footnotesize $v_{6}$}
\put(1.2,0.9){\footnotesize $\alpha$}
\put(2,2){\circle*{0.15}} 
\put(1.4,1.9){\footnotesize $v_{5}$}
\put(2.2,1.9){\footnotesize $\alpha$}
\put(3,3){\circle*{0.15}} 
\put(2.4,2.9){\footnotesize $v_{4}$}
\put(3.2,2.9){\footnotesize $\alpha$}
\put(4,4){\circle*{0.15}} 
\put(3.4,3.9){\footnotesize $v_{3}$}
\put(4.2,3.9){\footnotesize $\alpha$}
\put(0,2){\circle*{0.15}} 
\put(-0.6,1.9){\footnotesize $v_{2}$}
\put(0.2,1.9){\footnotesize $\beta$}
\put(1,3){\circle*{0.15}} 
\put(0.4,2.9){\footnotesize $v_{1}$}
\put(1.2,2.9){\footnotesize $\beta$}
%Chains
\put(1,1){\line(1,1){3}}
\put(4,2){\line(1,1){1}}
\put(0,2){\line(1,1){1}} 
%Other edges 
\put(1,1){\line(-1,1){1}} 
\put(2,2){\line(-1,1){1}} 
\put(3,3){\line(1,-1){1}}
\put(4,4){\line(1,-1){1}}
%Chain labels
\put(1.5,3.5){$\mathcal{C}_{1}$}
\put(4.5,4.5){$\mathcal{C}_{2}$}
\put(5.5,3.5){$\mathcal{C}_{3}$}
\end{picture}
}
\end{picture}
\hspace*{1.5cm}
\setlength{\unitlength}{1cm}
\begin{picture}(8,8)
\put(0,6){\fbox{$\mathfrak{g} = C_{2}$}}
\put(0,0.25){
\begin{picture}(3,3.5)
%Vertex labels
\put(6,3){\circle*{0.15}} 
\put(5.4,2.9){\footnotesize $v_{14}$}
\put(6.2,2.9){\footnotesize $\alpha$}
\put(7,4){\circle*{0.15}} 
\put(6.4,3.9){\footnotesize $v_{13}$}
\put(7.2,3.9){\footnotesize $\alpha$}
\put(1,0){\circle*{0.15}} 
\put(0.4,-0.1){\footnotesize $v_{12}$}
\put(1.2,-0.1){\footnotesize $\beta$}
\put(3,2){\circle*{0.15}} 
\put(2.4,1.9){\footnotesize $v_{11}$}
\put(3.2,1.9){\footnotesize $\beta$}
\put(5,4){\circle*{0.15}} 
\put(4.4,3.9){\footnotesize $v_{10}$}
\put(5.2,3.9){\footnotesize $\beta$}
\put(6,5){\circle*{0.15}} 
\put(5.4,4.9){\footnotesize $v_{9}$}
\put(6.2,4.9){\footnotesize $\beta$}
\put(0,1){\circle*{0.15}} 
\put(-0.6,0.9){\footnotesize $v_{8}$}
\put(0.2,0.9){\footnotesize $\alpha$}
\put(1,2){\circle*{0.15}} 
\put(0.4,1.9){\footnotesize $v_{7}$}
\put(1.2,1.9){\footnotesize $\alpha$}
\put(2,3){\circle*{0.15}} 
\put(1.4,2.9){\footnotesize $v_{6}$}
\put(2.2,2.9){\footnotesize $\alpha$}
\put(3,4){\circle*{0.15}} 
\put(2.4,3.9){\footnotesize $v_{5}$}
\put(3.2,3.9){\footnotesize $\alpha$}
\put(4,5){\circle*{0.15}} 
\put(3.4,4.9){\footnotesize $v_{4}$}
\put(4.2,4.9){\footnotesize $\alpha$}
\put(5,6){\circle*{0.15}} 
\put(4.4,5.9){\footnotesize $v_{3}$}
\put(5.2,5.9){\footnotesize $\alpha$}
\put(0,3){\circle*{0.15}} 
\put(-0.6,2.9){\footnotesize $v_{2}$}
\put(0.2,2.9){\footnotesize $\beta$}
\put(2,5){\circle*{0.15}} 
\put(1.4,4.9){\footnotesize $v_{1}$}
\put(2.2,4.9){\footnotesize $\beta$}
%Chains
\put(0,1){\line(1,1){5}}
\put(1,0){\line(1,1){5}}
\put(0,3){\line(1,1){2}} 
\put(6,3){\line(1,1){1}} 
%Other edges 
\put(1,0){\line(-1,1){1}} 
\put(1,2){\line(-1,1){1}} 
\put(3,2){\line(-1,1){1}} 
\put(3,4){\line(-1,1){1}} 
\put(4,5){\line(1,-1){2}}
\put(5,6){\line(1,-1){2}}
%Chain labels
\put(2.5,5.5){$\mathcal{C}_{1}$}
\put(5.5,6.5){$\mathcal{C}_{2}$}
\put(6.5,5.5){$\mathcal{C}_{3}$}
\put(7.5,4.5){$\mathcal{C}_{4}$}
\end{picture}
}
\end{picture}

\setlength{\unitlength}{1cm}
\begin{picture}(11,13)
\put(-0.5,9){\fbox{$\mathfrak{g} = G_{2}$}}
\put(2,0){\circle*{0.15}} 
\put(1.4,-0.1){\footnotesize $v_{30}$}
\put(2.2,-0.15){\footnotesize $\beta$}
\put(1,1){\circle*{0.15}} 
\put(0.4,0.9){\footnotesize $v_{26}$}
\put(1.2,0.9){\footnotesize $\alpha$}
\put(2,2){\circle*{0.15}} 
\put(1.4,1.9){\footnotesize $v_{25}$}
\put(2.2,1.9){\footnotesize $\alpha$}
\put(1,3){\circle*{0.15}} 
\put(0.4,2.9){\footnotesize $v_{16}$}
\put(1.2,2.9){\footnotesize $\beta$}
\put(3,3){\circle*{0.15}} 
\put(2.4,2.9){\footnotesize $v_{24}$}
\put(3.2,2.9){\footnotesize $\alpha$}
\put(5,3){\circle*{0.15}} 
\put(4.4,2.9){\footnotesize $v_{29}$}
\put(5.2,2.85){\footnotesize $\beta$}
\put(0,4){\circle*{0.15}} 
\put(-0.6,3.9){\footnotesize $v_{10}$}
\put(0.2,3.9){\footnotesize $\alpha$}
\put(2,4){\circle*{0.15}} 
\put(1.4,3.9){\footnotesize $v_{15}$}
\put(2.2,3.9){\footnotesize $\beta$}
\put(4,4){\circle*{0.15}} 
\put(3.4,3.9){\footnotesize $v_{23}$}
\put(4.2,3.9){\footnotesize $\alpha$}
\put(1,5){\circle*{0.15}} 
\put(0.4,4.9){\footnotesize $v_{9}$}
\put(1.2,4.9){\footnotesize $\alpha$}
\put(5,5){\circle*{0.15}} 
\put(4.4,4.9){\footnotesize $v_{22}$}
\put(5.2,4.9){\footnotesize $\alpha$}
\put(9,5){\circle*{0.15}} 
\put(8.4,4.9){\footnotesize $v_{32}$}
\put(9.2,4.9){\footnotesize $\alpha$}
\put(2,6){\circle*{0.15}} 
\put(1.4,5.9){\footnotesize $v_{8}$}
\put(2.2,5.9){\footnotesize $\alpha$}
\put(4,6){\circle*{0.15}} 
\put(3.4,5.9){\footnotesize $v_{14}$}
\put(4.2,5.9){\footnotesize $\beta$}
\put(6,6){\circle*{0.15}} 
\put(5.4,5.9){\footnotesize $v_{21}$}
\put(6.2,5.9){\footnotesize $\alpha$}
\put(8,6){\circle*{0.15}} 
\put(7.4,5.9){\footnotesize $v_{28}$}
\put(8.2,5.9){\footnotesize $\beta$}
\put(1,7){\circle*{0.15}} 
\put(0.4,6.9){\footnotesize $v_{2}$}
\put(1.2,6.9){\footnotesize $\beta$}
\put(3,7){\circle*{0.15}} 
\put(2.4,6.9){\footnotesize $v_{7}$}
\put(3.2,6.9){\footnotesize $\alpha$}
\put(5,7){\circle*{0.15}} 
\put(4.4,6.9){\footnotesize $v_{13}$}
\put(5.2,6.9){\footnotesize $\beta$}
\put(7,7){\circle*{0.15}} 
\put(6.4,6.9){\footnotesize $v_{20}$}
\put(7.2,6.9){\footnotesize $\alpha$}
\put(11,7){\circle*{0.15}} 
\put(10.4,6.9){\footnotesize $v_{31}$}
\put(11.2,6.9){\footnotesize $\alpha$}
\put(4,8){\circle*{0.15}} 
\put(3.4,7.9){\footnotesize $v_{6}$}
\put(4.2,7.9){\footnotesize $\alpha$}
\put(8,8){\circle*{0.15}} 
\put(7.4,7.9){\footnotesize $v_{19}$}
\put(8.2,7.9){\footnotesize $\alpha$}
\put(10,8){\circle*{0.15}} 
\put(9.4,7.9){\footnotesize $v_{27}$}
\put(10.2,7.9){\footnotesize $\beta$}
\put(5,9){\circle*{0.15}} 
\put(4.4,8.9){\footnotesize $v_{5}$}
\put(5.2,8.9){\footnotesize $\alpha$}
\put(7,9){\circle*{0.15}} 
\put(6.4,8.9){\footnotesize $v_{12}$}
\put(7.2,8.9){\footnotesize $\beta$}
\put(9,9){\circle*{0.15}} 
\put(8.4,8.9){\footnotesize $v_{18}$}
\put(9.2,8.9){\footnotesize $\alpha$}
\put(4,10){\circle*{0.15}} 
\put(3.4,9.9){\footnotesize $v_{1}$}
\put(4.2,9.9){\footnotesize $\beta$}
\put(6,10){\circle*{0.15}} 
\put(5.4,9.9){\footnotesize $v_{4}$}
\put(6.2,9.9){\footnotesize $\alpha$}
\put(10,10){\circle*{0.15}} 
\put(9.4,9.9){\footnotesize $v_{17}$}
\put(10.2,9.9){\footnotesize $\alpha$}
\put(9,11){\circle*{0.15}} 
\put(8.4,10.9){\footnotesize $v_{11}$}
\put(9.2,10.9){\footnotesize $\beta$}
\put(8,12){\circle*{0.15}} 
\put(7.4,11.9){\footnotesize $v_{3}$}
\put(8.2,11.9){\footnotesize $\alpha$}
%Chains
\put(9,5){\line(1,1){2}} 
\put(2,0){\line(1,1){8}} 
\put(1,1){\line(1,1){9}} 
\put(1,3){\line(1,1){8}} 
\put(0,4){\line(1,1){8}} 
\put(1,7){\line(1,1){3}}
%Other edges 
\put(2,0){\line(-1,1){1}} 
\put(2,2){\line(-1,1){2}} 
\put(3,3){\line(-1,1){2}} 
\put(5,3){\line(-1,1){1}} 
\put(2,6){\line(-1,1){1}} 
\put(5,5){\line(-1,1){2}} 
\put(6,6){\line(-1,1){2}} 
\put(9,5){\line(-1,1){2}} 
\put(5,9){\line(-1,1){1}} 
\put(8,8){\line(-1,1){2}} 
\put(11,7){\line(-1,1){2}} 
\put(10,10){\line(-1,1){2}} 
%Chain labels
\put(11.5,7.5){$\mathcal{C}_{6}$}
\put(10.5,8.5){$\mathcal{C}_{5}$}
\put(10.5,10.5){$\mathcal{C}_{4}$}
\put(9.5,11.5){$\mathcal{C}_{3}$}
\put(8.5,12.5){$\mathcal{C}_{2}$}
\put(4.5,10.5){$\mathcal{C}_{1}$}
%Vertex labels
\end{picture} 
\end{center}
%\end{figure}

%\begin{figure}[h,t]
\newpage
\begin{center}
\GridPosetsII: Depicted below are four \dichromatic grid posets each 
possessing the max 
property.\\  
{\small (Each is the $\mathfrak{g}$-semistandard poset 
$P_{\mathfrak{g}}^{\alpha\beta}(2,2)$ of \S \SemiNum\  
for the 
indicated rank two semisimple Lie algebra $\mathfrak{g}$.)} 

\setlength{\unitlength}{1cm}
\begin{picture}(6,7)
\put(0,7){\fbox{$\mathfrak{g} = A_{1} \oplus A_{1}$}}
\put(0.5,2.75){
\begin{picture}(3,3.5)
%Vertex labels
\put(3,2){\circle*{0.15}} 
\put(2.4,1.9){\footnotesize $v_{4}$}
\put(3.2,1.9){\footnotesize $\beta$}
\put(4,3){\circle*{0.15}} 
\put(3.4,2.9){\footnotesize $v_{3}$}
\put(4.2,2.9){\footnotesize $\beta$}
\put(0,2){\circle*{0.15}} 
\put(-0.6,1.9){\footnotesize $v_{2}$}
\put(0.2,1.9){\footnotesize $\alpha$}
\put(1,3){\circle*{0.15}} 
\put(0.4,2.9){\footnotesize $v_{1}$}
\put(1.2,2.9){\footnotesize $\alpha$}
%Chains
\put(3,2){\line(1,1){1}}
\put(0,2){\line(1,1){1}} 
%Chain labels
\put(1.5,3.5){$\mathcal{C}_{1}$}
\put(4.5,3.5){$\mathcal{C}_{2}$}
\put(1.8,2.375){$\bigoplus$}
\end{picture}
}
\put(0,2){\fbox{$\mathfrak{g} = A_{2}$}}
\put(0,-2.5){
\begin{picture}(3,3.5)
%Vertex labels
\put(4,2){\circle*{0.15}} 
\put(3.4,1.9){\footnotesize $v_{8}$}
\put(4.2,1.9){\footnotesize $\alpha$}
\put(5,3){\circle*{0.15}} 
\put(4.4,2.9){\footnotesize $v_{7}$}
\put(5.2,2.9){\footnotesize $\alpha$}
\put(1,1){\circle*{0.15}} 
\put(0.4,0.9){\footnotesize $v_{6}$}
\put(1.2,0.9){\footnotesize $\beta$}
\put(2,2){\circle*{0.15}} 
\put(1.4,1.9){\footnotesize $v_{5}$}
\put(2.2,1.9){\footnotesize $\beta$}
\put(3,3){\circle*{0.15}} 
\put(2.4,2.9){\footnotesize $v_{4}$}
\put(3.2,2.9){\footnotesize $\beta$}
\put(4,4){\circle*{0.15}} 
\put(3.4,3.9){\footnotesize $v_{3}$}
\put(4.2,3.9){\footnotesize $\beta$}
\put(0,2){\circle*{0.15}} 
\put(-0.6,1.9){\footnotesize $v_{2}$}
\put(0.2,1.9){\footnotesize $\alpha$}
\put(1,3){\circle*{0.15}} 
\put(0.4,2.9){\footnotesize $v_{1}$}
\put(1.2,2.9){\footnotesize $\alpha$}
%Chains
\put(1,1){\line(1,1){3}}
\put(4,2){\line(1,1){1}}
\put(0,2){\line(1,1){1}} 
%Other edges 
\put(1,1){\line(-1,1){1}} 
\put(2,2){\line(-1,1){1}} 
\put(3,3){\line(1,-1){1}}
\put(4,4){\line(1,-1){1}}
%Chain labels
\put(1.5,3.5){$\mathcal{C}_{1}$}
\put(4.5,4.5){$\mathcal{C}_{2}$}
\put(5.5,3.5){$\mathcal{C}_{3}$}
\end{picture}
}
\end{picture}
\hspace*{1.5cm}
\setlength{\unitlength}{1cm}
\begin{picture}(8,8)
\put(0,6){\fbox{$\mathfrak{g} = C_{2}$}}
\put(0,0.25){
\begin{picture}(3,3.5)
%Vertex labels
\put(1,3){\VertexForPosets{1}{\alpha}}
\put(0,2){\VertexForPosets{2}{\alpha}}
\put(6,6){\VertexForPosets{3}{\beta}}
\put(4,4){\VertexForPosets{4}{\beta}}
\put(2,2){\VertexForPosets{5}{\beta}}
\put(1,1){\VertexForPosets{6}{\beta}}
\put(7,5){\VertexForPosets{7}{\alpha}}
\put(6,4){\VertexForPosets{8}{\alpha}}
\put(5,3){\VertexForPosets{9}{\alpha}}
\put(4,2){\VertexForPosets{10}{\alpha}}
\put(3,1){\VertexForPosets{11}{\alpha}}
\put(2,0){\VertexForPosets{12}{\alpha}}
\put(7,3){\VertexForPosets{13}{\beta}}
\put(5,1){\VertexForPosets{14}{\beta}}
%Chains
\put(0,2){\line(1,1){1}}
\put(1,1){\line(1,1){5}}
\put(2,0){\line(1,1){5}} 
\put(5,1){\line(1,1){2}} 
%Other edges 
\put(2,0){\line(-1,1){2}} 
\put(3,1){\line(-1,1){2}} 
\put(5,1){\line(-1,1){1}} 
\put(5,3){\line(-1,1){1}} 
\put(7,3){\line(-1,1){1}}
\put(7,5){\line(-1,1){1}}
%Chain labels
\put(1.5,3.5){$\mathcal{C}_{1}$}
\put(6.5,6.5){$\mathcal{C}_{2}$}
\put(7.5,5.5){$\mathcal{C}_{3}$}
\put(7.5,3.5){$\mathcal{C}_{4}$}
\end{picture}
}
\end{picture}

\setlength{\unitlength}{1cm}
\begin{picture}(11,13)
\put(-0.5,9){\fbox{$\mathfrak{g} = G_{2}$}}
%Vertices
\put(2,7){\VertexForPosets{1}{\alpha}}
\put(0,5){\VertexForPosets{2}{\alpha}}
\put(9,12){\VertexForPosets{3}{\beta}}
\put(6,9){\VertexForPosets{4}{\beta}}
\put(3,6){\VertexForPosets{5}{\beta}}
\put(1,4){\VertexForPosets{6}{\beta}}
\put(10,11){\VertexForPosets{7}{\alpha}}
\put(9,10){\VertexForPosets{8}{\alpha}}
\put(8,9){\VertexForPosets{9}{\alpha}}
\put(7,8){\VertexForPosets{10}{\alpha}}
\put(6,7){\VertexForPosets{11}{\alpha}}
\put(5,6){\VertexForPosets{12}{\alpha}}
\put(4,5){\VertexForPosets{13}{\alpha}}
\put(3,4){\VertexForPosets{14}{\alpha}}
\put(2,3){\VertexForPosets{15}{\alpha}}
\put(1,2){\VertexForPosets{16}{\alpha}}
\put(10,9){\VertexForPosets{17}{\beta}}
\put(9,8){\VertexForPosets{18}{\beta}}
\put(7,6){\VertexForPosets{19}{\beta}}
\put(6,5){\VertexForPosets{20}{\beta}}
\put(4,3){\VertexForPosets{21}{\beta}}
\put(2,1){\VertexForPosets{22}{\beta}}
\put(11,8){\VertexForPosets{23}{\alpha}}
\put(10,7){\VertexForPosets{24}{\alpha}}
\put(9,6){\VertexForPosets{25}{\alpha}}
\put(8,5){\VertexForPosets{26}{\alpha}}
\put(7,4){\VertexForPosets{27}{\alpha}}
\put(6,3){\VertexForPosets{28}{\alpha}}
\put(5,2){\VertexForPosets{29}{\alpha}}
\put(3,0){\VertexForPosets{30}{\alpha}}
\put(10,5){\VertexForPosets{31}{\beta}}
\put(7,2){\VertexForPosets{32}{\beta}}
%Chains
\put(7,2){\line(1,1){3}} 
\put(3,0){\line(1,1){8}} 
\put(2,1){\line(1,1){8}} 
\put(1,2){\line(1,1){9}} 
\put(1,4){\line(1,1){8}} 
\put(0,5){\line(1,1){2}}
%Other edges 
\put(3,0){\line(-1,1){2}} 
\put(2,3){\line(-1,1){2}} 
\put(5,2){\line(-1,1){2}} 
\put(4,5){\line(-1,1){2}} 
\put(7,2){\line(-1,1){1}} 
\put(7,4){\line(-1,1){2}} 
\put(8,5){\line(-1,1){2}} 
\put(7,8){\line(-1,1){1}} 
\put(10,5){\line(-1,1){1}} 
\put(10,7){\line(-1,1){2}} 
\put(11,8){\line(-1,1){2}} 
\put(10,11){\line(-1,1){1}} 
%Chain labels
\put(10.5,5.5){$\mathcal{C}_{6}$}
\put(11.5,8.5){$\mathcal{C}_{5}$}
\put(10.5,9.5){$\mathcal{C}_{4}$}
\put(10.5,11.5){$\mathcal{C}_{3}$}
\put(9.5,12.5){$\mathcal{C}_{2}$}
\put(2.5,7.5){$\mathcal{C}_{1}$}
%Vertex labels
\end{picture} 
\end{center}
%\end{figure}

%==========================
% End pagebreak for figures
%==========================
\noindent
a minimal and a maximal element in $\comp_{\gamma}(\selt)$, then it 
follows that $\relt_{q}$ and $\uelt_{p}$ are respectively the unique minimal 
and the unique maximal element of $\comp_{\gamma}(\selt)$.  Then 
$\rho_{\gamma}(\selt) = q$ and $l_{\gamma}(\selt) = p+q$. 

Reorganize the sequence $(v_{i_{1}},\ldots,v_{i_{q}})$ as follows: 
write $(v_{k_{1}},\ldots,v_{k_{q'}},v_{k_{q'+1}},\ldots,v_{k_{q}})$, 
where the vertices 
$v_{k_{1}},\ldots,v_{k_{q'}}$ are all in $P_{2}$ with $k_{1} < 
\cdots < k_{q'}$, and the vertices $v_{k_{q'+1}},\ldots,v_{k_{q}}$ 
are all in $P_{1}$ with $k_{q'+1} < \cdots < k_{q}$. Set 
$\relt'_{0} := \selt$ and for $j \geq 0$ set $\relt'_{j+1} := 
\relt'_{j} \setminus \{v_{k_{j+1}}\}$.  We claim that each $\relt'_{j+1}$ 
is an order ideal of $P$, and if $0 \leq j < q'$ (respectively $q' 
\leq j < q$) then $v_{k_{j+1}}$ is the smallest element in 
$\mathcal{T}_{P}$ of color $\gamma$ 
that is also in $P_{2}$ (respectively $P_{1}$) that 
can be removed from $\relt'_{j}$ so that $\relt'_{j+1}$ is an order 
ideal of $P$.  (If so, we have a path 
$\relt_{q} = \relt'_{q} \myarrow{\gamma} 
\relt'_{q-1} \myarrow{\gamma} \cdots \myarrow{\gamma} \relt'_{1} 
\myarrow{\gamma} \relt'_{0} = \relt_{0} = \selt$ in $L$.)  
Proceed by induction on $j$.  
The statement 
follows if we can show that $v_{k_{j+1}}$ is maximal in $\relt'_{j}$. 
First suppose $0 \leq j < q'$.  
If $v_{k_{j+1}}$ is not maximal in $\relt'_{j}$, then $v_{k_{j+1}} < v$ 
for some other maximal element $v$ in $\relt'_{j}$.  It must be the 
case that $v$ is one of 
$(v_{k_{j+2}},\ldots,v_{k_{q'}},v_{k_{q'+1}},\ldots,v_{k_{q}})$; 
otherwise one could not descend from $\relt'_{j}$ to $\relt_{q}$ 
in $L = J_{color}(P)$ along edges corresponding to vertices from 
$(v_{k_{j+2}},\ldots,v_{k_{q'}},v_{k_{q'+1}},\ldots,v_{k_{q}})$. 
It cannot be the case that $v$ is 
one of $(v_{k_{j+2}},\ldots,v_{k_{q'}})$; otherwise $k_{1}< \cdots 
< k_{j+1} < k_{j+2} < 
\cdots < k_{q'}$ implies that $v$ is larger than $v_{k_{j+1}}$ in the 
total order $\mathcal{T}_{P}$, violating the fact that $v_{k_{j+1}} < 
v$ in $P$.  And it cannot be the case that $v$ is one of 
$(v_{k_{q'+1}},\ldots,v_{k_{q}})$ since these are elements 
of $P_{1}$ and $v_{k_{j+1}}$ is in $P_{2}$.  So for $0 \leq j < q'$, the vertex 
$v_{k_{j+1}}$ is maximal in $\relt'_{j}$. Second, suppose that 
$q' \leq j < 
q$.  If $v_{k_{j+1}}$ is not maximal in $\relt'_{j}$, then by 
reasoning similar to the preceding case we have $v_{k_{j+1}} < v$ for 
some element $v$ 
from $(v_{k_{j+2}},\ldots,v_{k_{q}})$. But this violates the fact 
that  $v_{k_{j+1}}$ precedes $v$ in the total order $\mathcal{T}_{P}$ 
since $k_{q'} < \cdots < k_{j+1} < k_{j+2} < \cdots < k_{q}$.  
So for $q' \leq j < q$, the vertex 
$v_{k_{j+1}}$ is maximal in $\relt'_{j}$.  This concludes our 
induction on $j$. 

Let $\relt^{(1)}$ be the unique minimal element in the 
$\gamma$-component $\comp_{\gamma}^{(1)}(\selt \cap P_{1})$ 
of $\selt \cap P_{1}$ in the edge-colored distributive lattice 
$L_{1} = J_{color}(P_{1})$.  We 
claim that $\relt^{(1)} = \xelt$, where $\xelt := (\selt \cap P_{1}) \setminus 
\{v_{k_{q'+1}},\ldots,v_{k_{q}}\}$.  Now $\xelt$ is an order ideal 
of $P_{1}$ since $\xelt = \relt'_{q} \cap P_{1} = \relt_{q} \cap P_{1}$.  
Also, $\xelt \in 
\comp_{\gamma}^{(1)}(\selt \cap P_{1})$ since $\selt \cap P_{1} = 
\relt'_{q'} \cap P_{1}$ and the path $\xelt \myarrow{\gamma} 
(\relt'_{q-1} \cap P_{1}) \myarrow{\gamma} \cdots \myarrow{\gamma} 
(\relt'_{q'+1} \cap P_{1}) \myarrow{\gamma} (\relt'_{q'} \cap P_{1})$ stays in 
$\comp_{\gamma}^{(1)}(\selt \cap P_{1})$.  
If $\relt^{(1)} \not= \xelt$, then $\relt^{(1)} < \xelt$.  In this 
case let $u \in \xelt$ be any color $\gamma$ vertex  
such that $\xelt \setminus \{u\}$ is 
an order ideal of $P_{1}$.  
Let $\mathcal{C}_{i}$ be the chain in $P$ that contains $u$.  
Note that $u$ is not 
maximal in $\relt_{q} \subseteq P$, and hence $u \rightarrow u'$ is a covering 
relation in $P$ for some $u'$ in $\relt_{q}$. 
We refer to the following as observation ({\tt *}): If $w$ is any element of 
$P$ such that $u \rightarrow w$ and $w \in \relt_{q}$, then $w \in 
P_{2}$.  (Otherwise $w \in P_{1}$, so that $w \in \xelt$, and 
then $\xelt \setminus \{u\}$ cannot be an 
order ideal of $P_{1}$.)  In particular, $u' \in P_{2}$.  We claim that 
$u' \not\in \{v_{k_{1}},\ldots,v_{k_{q'}}\}$.  Indeed, if 
$u' \in \{v_{k_{1}},\ldots,v_{k_{q'}}\}$, then since $u \not\in 
\{v_{k_{q'+1}},\ldots,v_{k_{q}}\}$, it must be the case that  
$u \rightarrow u''$ for some $u'' \in \relt_{q}$ in 
$\mathcal{C}_{i-1}$.  By observation ({\tt *}), the element 
$u''$ is in $P_{2}$. But a covering relation in $P$ between 
elements of $P_{1}$ 
and elements of $P_{2}$ can only occur along the 
chains $\mathcal{C}_{1},\ldots,\mathcal{C}_{m}$.  Therefore $u'' \in 
\mathcal{C}_{i}$, which contradicts the fact that $u'' \in 
\mathcal{C}_{i-1}$.  So it must be the 
case that $u' \not\in \{v_{k_{1}},\ldots,v_{k_{q'}}\}$. It follows 
that $u' < u''$ for some $u'' \in \relt_{q}$ in 
$\mathcal{C}_{i-1}$.  Let $v$ be a maximal element in $P$ such that 
$u'' \leq v$. Note that $v \in P_{2}$ since $u'' \in P_{2}$.  Moreover, $v \in 
\mathcal{C}_{j}$ with $j \leq i-1$.  Next suppose $u \rightarrow z$ 
for some $z \in P_{1}$.  Since $z \in P_{1}$, then $z \not= u'$. 
Therefore $z \in \mathcal{C}_{i-1}$.  Therefore $z \leq u''$.  But 
since $u''$ is in the order ideal $\relt_{q}$, it follows that $z \in 
\relt_{q}$.  But by observation ({\tt *}), it now follows that $z \in P_{2}$.  
This contradicts our hypothesis that $z \in P_{1}$.  In particular, 
$u$ must be a maximal element in $P_{1}$.  So $u \in \mathcal{C}_{i}$ 
is a maximal element in $P_{1}$ and $v \in \mathcal{C}_{j}$ 
is a maximal element in $P_{2}$, and $j < i$.  This violates 
the fact that $P$ decomposes into $P_{1} \triangleleft 
P_{2}$.  So $\relt^{(1)} = \xelt$, and hence 
$\rho^{(1)}_{\gamma}(\selt \cap P_{1}) = q - q'$. 

Let $\relt^{(2)}$ be the unique minimal element in the 
$\gamma$-component $\comp_{\gamma}^{(2)}(\selt \cap P_{2})$ 
of $\selt \cap P_{2}$ in the edge-colored distributive lattice 
$L_{2} = J_{color}(P_{2})$.  We 
claim that $\relt^{(2)} = \yelt$, where $\yelt := (\selt \cap P_{2}) \setminus 
\{v_{k_{1}},\ldots,v_{k_{q'}}\}$.  Now $\yelt$ is an order ideal 
of $P_{2}$ since $\yelt = \relt'_{q'} \cap P_{2} = \relt_{q} \cap P_{2}$.  
Also, $\yelt \in 
\comp_{\gamma}^{(2)}(\selt \cap P_{2})$ since the path $\yelt \myarrow{\gamma} 
(\relt'_{q'-1} \cap P_{2}) \myarrow{\gamma} \cdots \myarrow{\gamma} 
(\relt'_{1} \cap P_{2}) \myarrow{\gamma} (\relt'_{0} \cap P_{2})$ stays in 
$\comp_{\gamma}^{(2)}(\selt \cap P_{2})$.  
If $\relt^{(2)} \not= \yelt$, then $\relt^{(2)} < \yelt$.  In this 
case let $u \in \yelt$ be any color $\gamma$ vertex such 
that $\yelt \setminus \{u\}$ is 
an order ideal of $P_{2}$.  In particular, $u$ is a maximal element 
in $\yelt$.  Let $w$ be any element of $\relt_{q}$ with $u \not= w$.  If 
$w \in P_{2}$, then $w \in \yelt$, so $u \not< w$.  If $w \in P_{1}$, 
then by properties of the decomposition of $P$ into $P_{1} 
\triangleleft P_{2}$, it cannot be the case that $u < w$.  Therefore $u$ is a 
maximal element of $\relt_{q}$ of color $\gamma$.  But this 
contradicts the fact that $\relt_{q}$ is the minimal element in 
$\comp_{\gamma}(\selt)$.  So it is not the case that $\relt^{(2)} < 
\yelt$.  Therefore $\relt^{(2)} = \yelt$, and so 
$\rho^{(2)}_{\gamma}(\selt \cap P_{2}) = q'$.  Combine this with 
$\rho^{(1)}_{\gamma}(\selt \cap P_{1}) = q - q'$ to see that 
$\rho_{\gamma}(\selt) = q = (q-q') + q' = 
\rho^{(1)}_{\gamma}(\selt \cap P_{1}) + 
\rho^{(2)}_{\gamma}(\selt \cap P_{2})$. 

The dual $P^{*}$ may be 
viewed as a \dichromatic grid poset that decomposes into 
$P_{2}^{*} \triangleleft P_{1}^{*}$.  Order ideals of $P^{*}$ are 
complements of order ideals of $P$.  
Then arguments analogous to those above apply to the complements of 
elements of the sequence $\selt = \uelt_{0}, \uelt_{1},\ldots, 
\uelt_{p}$.  So we obtain: $\rho^{*}_{\gamma}(P \setminus \selt) = 
\rho^{* (2)}_{\gamma}((P \setminus \selt) \cap P_{2}) +
\rho^{* (1)}_{\gamma}((P \setminus \selt) \cap P_{1})$.  Note that 
$(P \setminus \selt) 
\cap P_{i} = P_{i} \setminus (\selt \cap P_{i})$ for $i \in \{1,2\}$.  Now  
$l_{\gamma}(\selt) = \rho_{\gamma}(\selt) + 
\rho_{\gamma}^{*}(P \setminus \selt)$, 
$l^{(1)}_{\gamma}(\selt \cap P_{1}) = \rho^{(1)}_{\gamma}(\selt \cap P_{1}) + 
\rho^{* (1)}_{\gamma}(P_{1} \setminus (\selt \cap P_{1}))$, and 
$l^{(2)}_{\gamma}(\selt \cap P_{2}) = \rho^{(2)}_{\gamma}(\selt \cap P_{2}) + 
\rho^{* (2)}_{\gamma}(P_{2} \setminus (\selt \cap P_{2}))$.   
Therefore $l_{\gamma}(\selt) = l^{(1)}_{\gamma}(\selt \cap P_{1}) + 
l^{(2)}_{\gamma}(\selt \cap P_{2})$.\hfill\QED

It can be shown that if either of the conditions on the maximal and 
minimal elements on  $P_{1}$  and  $P_{2}$  required for the statement  
``$P = P_{1} \triangleleft P_{2}$''  fail,  
then so does at least one of the decomposition equations 
in \WeightsLemma\ for $\rho_{\gamma}(\selt)$ and $l_{\gamma}(\selt)$.

%==================================================================
%\newpage
\vspace{2ex} 

\noindent
{\bf \Large \SemiNum.\ \ $\mathfrak{g}$-semistandard posets, lattices, and tableaux}

We define special two-color grid posets  $P$, the ``$\mathfrak{g}$-semistandard 
posets''.   
Then we define corresponding lattices  $L = J_{color}(P)$, the 
``$\mathfrak{g}$-semistandard''  
lattices.  In the second half of the section, ``$\mathfrak{g}$-semistandard'' tableau 
descriptions of the elements of these lattices are developed.

For the remainder of this paper, $\mathfrak{g}$ denotes 
a rank two semisimple Lie
algebra: $\mathfrak{g} \in \{A_{1} \oplus A_{1}, A_{2}, C_{2}, G_{2}\}$. 
We identify $\alpha$ with a short simple 
root for $\mathfrak{g}$ 
and $\beta$ as the other simple root.   
The vertex colors for the posets and the edge colors for the lattices 
which we now introduce correspond to the simple roots of 
$\mathfrak{g}$.  
So here the index set $I$  
of Section \DefsNum\ becomes  $I = \{\alpha,\beta\}$. 
Let $\omega_{\alpha} = \omega_{1} = (1,0)$ 
and 
$\omega_{\beta} = \omega_{2} = (0,1)$ respectively denote the 
corresponding fundamental weights. 
Then any weight $\mu$ in $\Lambda$ of the form $\mu = p\omega_{\alpha} + 
q\omega_{\beta}$ (where $p$ and $q$ are integers) is now identified 
with the pair $(p,q)$ in $\mathbb{Z} \times \mathbb{Z}$. \hfill 
In particular, $\alpha$ and $\beta$ are respectively identified with the first 
and second 
row vectors 
from the Cartan matrix $M$ for $\mathfrak{g}$.  These 
matrices, displayed in  \CartanTable, 
specify the  $\mathfrak{g}$-structure condition of Section \DefsNum\ 
for edge-colored ranked 
posets. 

The $\mathfrak{g}$-{\em 
fundamental posets} $P_{\mathfrak{g}}(1,0)$ and 
$P_{\mathfrak{g}}(0,1)$ are defined to be the  
\dichromatic grid posets of \FundPosets. 
The corresponding $\mathfrak{g}$-{\em fundamental 
lattices} are defined to be the edge-colored 
lattices $L_{\mathfrak{g}}(1,0) := 
J_{color}(P_{\mathfrak{g}}(1,0))$ 
and 
$L_{\mathfrak{g}}(0,1) := 
J_{color}(P_{\mathfrak{g}}(0,1))$. 
\hspace*{0.1in}  
See \FundLatticeIdealsFigure.  
For the remainder of
%=================================
% Text before figures
%=================================

\newpage

\hspace*{0.25in}\CartanTable\ \hspace*{0.25in}
\begin{tabular}{|c|c|c|c|}
\hline 
$A_{1} \oplus A_{1}$ & $A_{2}$ & $C_{2}$ & $G_{2}$\\
\hline 
\hline
\rule[-7mm]{0mm}{16mm} $\displaystyle \left(\begin{array}{cc} 2 & 0 \\ 
0 & 2\end{array}\right)$ 
&  
\rule[-7mm]{0mm}{16mm} $\displaystyle \left(\begin{array}{cc} 2 & -1 \\ 
-1 & 2\end{array}\right)$ 
& 
\rule[-7mm]{0mm}{16mm} $\displaystyle \left(\begin{array}{cc} 2 & -1 \\ 
-2 & 2\end{array}\right)$ 
& 
\rule[-7mm]{0mm}{16mm} $\displaystyle \left(\begin{array}{cc} 2 & -1 \\ 
-3 & 2\end{array}\right)$\\
\hline 
\end{tabular}

%=================
% Figure
%=================
%\begin{figure}[htb]
\begin{center}
\FundPosets: $\mathfrak{g}$-fundamental posets. 

\

\vspace*{-0.15in} 

\begin{tabular}{|c||c|c|}
\hline
Algebra $\mathfrak{g}$ & $P_{\mathfrak{g}}(1,0)$ & 
$P_{\mathfrak{g}}(0,1)$\\
\hline
\hline
\parbox[b]{0.5in}{\ 

$A_{1} \oplus A_{1}$

\ 

\vspace*{-0.1in}
} & \AOneAOneAlpha & \AOneAOneBeta\\
\hline
\parbox[b]{0.225in}{$A_{2}$
\vspace*{0.4in}} 
 & \ATwoAlpha & \ATwoBeta\\
\hline
\parbox[b]{0.225in}{$C_{2}$
\vspace*{1.1in}} 
 & \BTwoAAlpha & \BTwoBBeta\\
\hline
\parbox[b]{0.225in}{$G_{2}$
\vspace*{2.5in}} 
 & \GTwoAAlpha & \GTwoBBeta\\
\hline
\end{tabular}
\end{center}
%\end{figure}

%=================================
% Text after figures
%=================================
\noindent
this section, everything presented for the {\em simple} 
cases ($A_{2}$, $C_{2}$, and $G_{2}$) has an easy $A_{1} \oplus 
A_{1}$ analog.  The details for $A_{1} \oplus A_{1}$ are omitted 
to save space, beginning with \FundLatticeIdealsFigure. 

Let $\lambda = (a,b)$, with $a, b \geq 0$.  
The $\mathfrak{g}$-{\em semistandard poset} 
$P_{\mathfrak{g}}^{\beta\alpha}(\lambda)$ associated to $\lambda$ is 
defined to be the \dichromatic grid poset $P$ which has the 
decomposition 
$P_{1} \triangleleft P_{2} \triangleleft \cdots \triangleleft 
P_{a+b}$, where $P_{i}$ is vertex-color isomorphic to  
$P_{\mathfrak{g}}(0,1)$ for $1 \leq i \leq b$ and to 
$P_{\mathfrak{g}}(1,0)$ for $1+b \leq i \leq a+b$. 
\hfill It can be seen that $P$ is 
%=================================
% Text before figure
%=================================

\newpage
%\begin{figure}[ht]
\begin{center}
\setlength{\unitlength}{1.1cm}
\begin{picture}(14,13.3)
\put(0.845,12.9){\FundLatticeIdealsFigure: Elements of $\mathfrak{g}$-fundamental  
lattices as order ideals of $\mathfrak{g}$-fundamental posets.}
\put(0.765,12.4){\parbox{6.15in}{\small \begin{center}
(Each order ideal 
is identified by the 
indices of its maximal vertices.
) 
\end{center}
}}
\put(0,0){\line(0,1){12}}
\put(6,0){\line(0,1){12}}
\put(14,0){\line(0,1){12}}
\put(0,7){\line(1,0){6}}
\put(0,0){\line(1,0){14}}
\put(0,12){\line(1,0){14}}
\put(0.5,11.25){\fbox{\Large $A_{2}$}}
\put(0,6.75){\AtwoAlphaIdeals}
\put(3,6.75){\AtwoBetaIdeals}
\put(0.5,6.25){\fbox{\Large $C_{2}$}}
\put(0,0.75){\BtwoAlphaIdeals}
\put(3,0.25){\BtwoBetaIdeals}
\put(6.5,11.25){\fbox{\Large $G_{2}$}}
\put(6,1.75){\GtwoAlphaIdeals}
\put(8,-0.25){\GtwoBetaIdeals}
\end{picture}
\end{center} 
%\end{figure}

%=================================
% Text after figure
%=================================
\noindent 
unique up to isomorphism.  
For each semisimple Lie algebra $\mathfrak{g}$, the poset 
$P_{\mathfrak{g}}^{\beta\alpha}(2,2)$ is depicted in \GridPosets.  
The $\mathfrak{g}$-{\em semistandard poset} 
$P_{\mathfrak{g}}^{\alpha\beta}(\lambda)$ associated to $\lambda$ is 
analogously defined, except with $P_{i}$ vertex-color isomorphic to  
$P_{\mathfrak{g}}(1,0)$ for $1 \leq i \leq a$ and to 
$P_{\mathfrak{g}}(0,1)$ for
$a+1 \leq i \leq a+b$.  
See \GridPosetsII\ for the corresponding 
$P_{\mathfrak{g}}^{\alpha\beta}(2,2)$. 
Note that 
$P_{\mathfrak{g}}^{\beta\alpha}(1,0) = 
P_{\mathfrak{g}}^{\alpha\beta}(1,0) = 
P_{\mathfrak{g}}(1,0)$, and   
$P_{\mathfrak{g}}^{\beta\alpha}(0,1) = 
P_{\mathfrak{g}}^{\alpha\beta}(0,1) = 
P_{\mathfrak{g}}(0,1)$. If $a = b = 0$, then 
$P_{\mathfrak{g}}^{\beta\alpha}(\lambda)$ and 
$P_{\mathfrak{g}}^{\alpha\beta}(\lambda)$ are the empty set. 
The $\mathfrak{g}$-{\em semistandard lattices} associated to 
$\lambda$ are the edge-colored lattices 
$\Lba := 
J_{color}(P_{\mathfrak{g}}^{\beta\alpha}(\lambda))$ 
and $L_{\mathfrak{g}}^{\alpha\beta}(\lambda) := 
J_{color}(P_{\mathfrak{g}}^{\alpha\beta}(\lambda))$. 
Note that 
$L_{\mathfrak{g}}^{\beta\alpha}(1,0) = 
L_{\mathfrak{g}}^{\alpha\beta}(1,0) = 
L_{\mathfrak{g}}(1,0)$, and   
$L_{\mathfrak{g}}^{\beta\alpha}(0,1) = 
L_{\mathfrak{g}}^{\alpha\beta}(0,1) = 
L_{\mathfrak{g}}(0,1)$. We will not consider 
``mixed'' concatenations, where some copies of 
$P_{\mathfrak{g}}(0,1)$ are interlaced amongst copies of $P_{\mathfrak{g}}(1,0)$. 
Any such concatenation will not have the max property, 
which {\sl is} possessed by all of the $\mathfrak{g}$-semistandard 
posets.

Each $\mathfrak{g}$-semistandard lattice is an edge-colored poset.  
From now on we write $wt(\selt)$ for $wt_{L}(\selt)$ when $L$ is 
$\mathfrak{g}$-semistandard.  Let $\selt \in L$. Let $\gamma \in 
\{\alpha,\beta\}$.  By definition, the 
$\gamma$-entry of the 2-tuple $wt(\selt)$ is the rank of $\selt$ 
within the $\gamma$-colored connected component of $\selt$ diminished 
by the depth of $\selt$ in that component. 

\noindent
{\bf \StructureResultForFundamentals}\ \ {\sl 
Let $\selt \myarrow{\gamma} \telt$ be an 
edge of color $\gamma \in \{\alpha,\beta\}$ in a 
$\mathfrak{g}$-fundamental lattice $L$.  
Then $wt(\selt) + \gamma = wt(\telt)$. Hence each 
$\mathfrak{g}$-fundamental lattice satisfies the 
$\mathfrak{g}$-structure condition.} 

{\em Proof.}  Note that  $\selt$  and  $\telt$  
are in the same $\gamma$-component.  
Since  $\telt$  covers  $\selt$ in this component,  the $\gamma$-entry of  
$wt(\telt)$   
is  2  more than the  $\gamma$-entry of  $wt(\selt)$.  
But adding the simple root  $\gamma$  
to  $wt(\selt)$  adds  2  to the  $\gamma$-entry of  $wt(\selt)$,  
since  $M_{\gamma,\gamma} = 2$  always.  
Let  $\gamma'$  in  $I$  be such that  $\gamma' \not= \gamma$.  
Using \FundLatticeTableauxFigure, one can quickly 
check by hand that the  $\gamma'$-entry of $wt(\selt)$  changes by  
$M_{\gamma,\gamma'}$  or by  $M_{\gamma',\gamma}$  
(as appropriate) for each edge 
within each  $\gamma$-component of a $\mathfrak{g}$-fundamental 
lattice.\hfill\QED

\noindent
{\bf \StructureResult}\ \ {\sl Let 
$\lambda = (a,b)$, with $a, b \geq 0$. Let $L$ be 
one of the 
$\mathfrak{g}$-semistandard lattices $L_{\mathfrak{g}}^{\beta\alpha}(\lambda)$ 
or $L_{\mathfrak{g}}^{\alpha\beta}(\lambda)$.  
Let $\selt \myarrow{\gamma} \telt$ be an 
edge of color $\gamma \in \{\alpha,\beta\}$ in $L$.  Then 
$wt(\selt) + \gamma = wt(\telt)$, and hence $L$ satisfies the 
$\mathfrak{g}$-structure condition.}

{\em Proof.} In light of  
{\StructureResultForFundamentals}, apply part (2) of 
{\WeightsLemma}.\hfill\QED

\noindent 
{\bf \NewRemark}\ \ If $\mathfrak{g} = A_{1} \oplus 
A_{1}$,   
then we have $P_{\mathfrak{g}}^{\beta\alpha}(\lambda) \cong 
P_{\mathfrak{g}}^{\alpha\beta}(\lambda)$ as vertex-colored posets:  
their Hasse diagrams are vertex-color isomorphic to 
$P_{\mathfrak{g}}^{\beta\alpha}(a,0) 
\oplus P_{\mathfrak{g}}^{\beta\alpha}(0,b) \cong \mathbf{a} \oplus 
\mathbf{b}$.     
Hence $L_{\mathfrak{g}}^{\beta\alpha}(\lambda)$ and 
$L_{\mathfrak{g}}^{\alpha\beta}(\lambda)$ are edge-color isomorphic to 
$L_{\mathfrak{g}}^{\beta\alpha}(a,0) \times 
L_{\mathfrak{g}}^{\beta\alpha}(0,b) \cong (\mathbf{a+1}) \times 
(\mathbf{b+1})$. For $\mathfrak{g} = C_{2}$ or $\mathfrak{g} = G_{2}$, observe 
that $P_{\mathfrak{g}}^{\alpha\beta}(\lambda)$ 
is vertex-color isomorphic to 
$(P_{\mathfrak{g}}^{\beta\alpha}(\lambda))^{*}$, 
and thus $L_{\mathfrak{g}}^{\alpha\beta}(\lambda)$ and 
$(L_{\mathfrak{g}}^{\beta\alpha}(\lambda))^{*}$ are isomorphic as 
edge-colored posets.  For $\mathfrak{g} = A_{2}$, 
$P_{\mathfrak{g}}^{\alpha\beta}(\lambda)$ and 
$(P_{\mathfrak{g}}^{\beta\alpha}(\lambda))^{*}$ are isomorphic as posets, 
but their vertex colors 
are reversed; disregarding edge colors, it follows that 
$L_{\mathfrak{g}}^{\alpha\beta}(\lambda)$ and 
$(L_{\mathfrak{g}}^{\beta\alpha}(\lambda))^{*}$ are isomorphic as posets. 
In all cases, $L_{\mathfrak{g}}^{\alpha\beta}(\lambda) \cong 
(L_{\mathfrak{g}}^{\beta\alpha}(\lambda))^{\triangle}$.\hfill\QED

The easy proof of the following statement will be omitted: 

\noindent{\bf \DistinctLemma}\ \ {\sl Let $\lambda = (a,b)$, with $a, 
b \geq 0$. 
If $\mathfrak{g}$ is simple,  
then $L_{\mathfrak{g}}^{\beta\alpha}(\lambda) 
\cong L_{\mathfrak{g}}^{\alpha\beta}(\lambda)$ as edge-colored posets if 
and only if $a = 0$ or $b = 0$.}

Now we develop tableau labels for the elements of half of the  
$\mathfrak{g}$-semistandard lattices, the  $\Lba$.  Comments relating these tableaux 
to tableaux developed by some of us and other authors appear in Section 5. 
We associate to the fundamental weight $\omega_{\alpha} = (1,0)$  
the shape $\mathbf{shape}(1,0) = \!$\lengthonecolumn{};  
we associate to $\omega_{\beta} = (0,1)$  
the shape 
$\mathbf{shape}(0,1) = \!$\lengthtwocolumn{}{}. 
For $a, b \geq 0$, we associate to 
$\lambda = (a,b)$ the 
shape (Ferrers diagram) with $b$ columns of length two and $a$ 
columns of length one.  
A {\em tableau of shape} $\lambda$ 
is a filling of the boxes of $\mathbf{shape}(\lambda)$ with 
entries from some totally ordered set. 
For a tableau $T$ of 
shape $\lambda$, we 
write $T = (T^{(1)},\ldots,T^{(a+b)})$, where $T^{(i)}$ is the $i$th 
column of $T$ from the left.  
We let $T^{(i)}_{j}$ denote the $j$th entry of the 
column $T^{(i)}$, counting from the top.   
The tableau $T$ is {\em semistandard} if  the entries 
weakly increase across rows and 
strictly increase down columns.
To each element $\telt$ 
of a $\mathfrak{g}$-fundamental lattice from \FundLatticeIdealsFigure\ 
we associate 
the one-column semistandard tableau $\mathbf{tableau}(\telt)$ 
of \FundLatticeTableauxFigure. 
For an order ideal $\telt$ of $P_{\mathfrak{g}}^{\beta\alpha}(\lambda)$, 
let $\mathbf{tableau}(\telt)$ be the tableau $T = 
(T^{(1)},\ldots,T^{(a+b)})$ with $T^{(i)} = 
\mathbf{tableau}(\telt \cap P_{i})$. 
A tableau $T$ of 
shape $\lambda$ obtained in this way is a 
$\mathfrak{g}$-{\em semistandard tableau of shape} 
$\lambda$.  We let 
$\TabSet$ denote the set of all 
$\mathfrak{g}$-semistandard tableaux of shape 
$\lambda$.  The function 
$\mathbf{tableau} :\, L_{\mathfrak{g}}^{\beta\alpha}(\lambda) \longrightarrow 
\TabSet$ is a 
one-to-one correspondence. See \BtwoMolevFigureList\ for a $C_{2}$ 
example. 

\noindent
{\bf \TableauProp}\ \ {\sl Let $a, b \geq 0$, 
and let $\lambda = (a,b)$.  Then:} 

$\mathcal{S}_{A_{2}}(\lambda) = \Big\{$semistandard tableau $T$ of 
shape $\lambda$ with entries from $\{1,2,3\}\, \Big\}$

%\newpage 
%\begin{figure}[ht]
\begin{center}
\setlength{\unitlength}{1.1cm}
\begin{picture}(14,12.8)
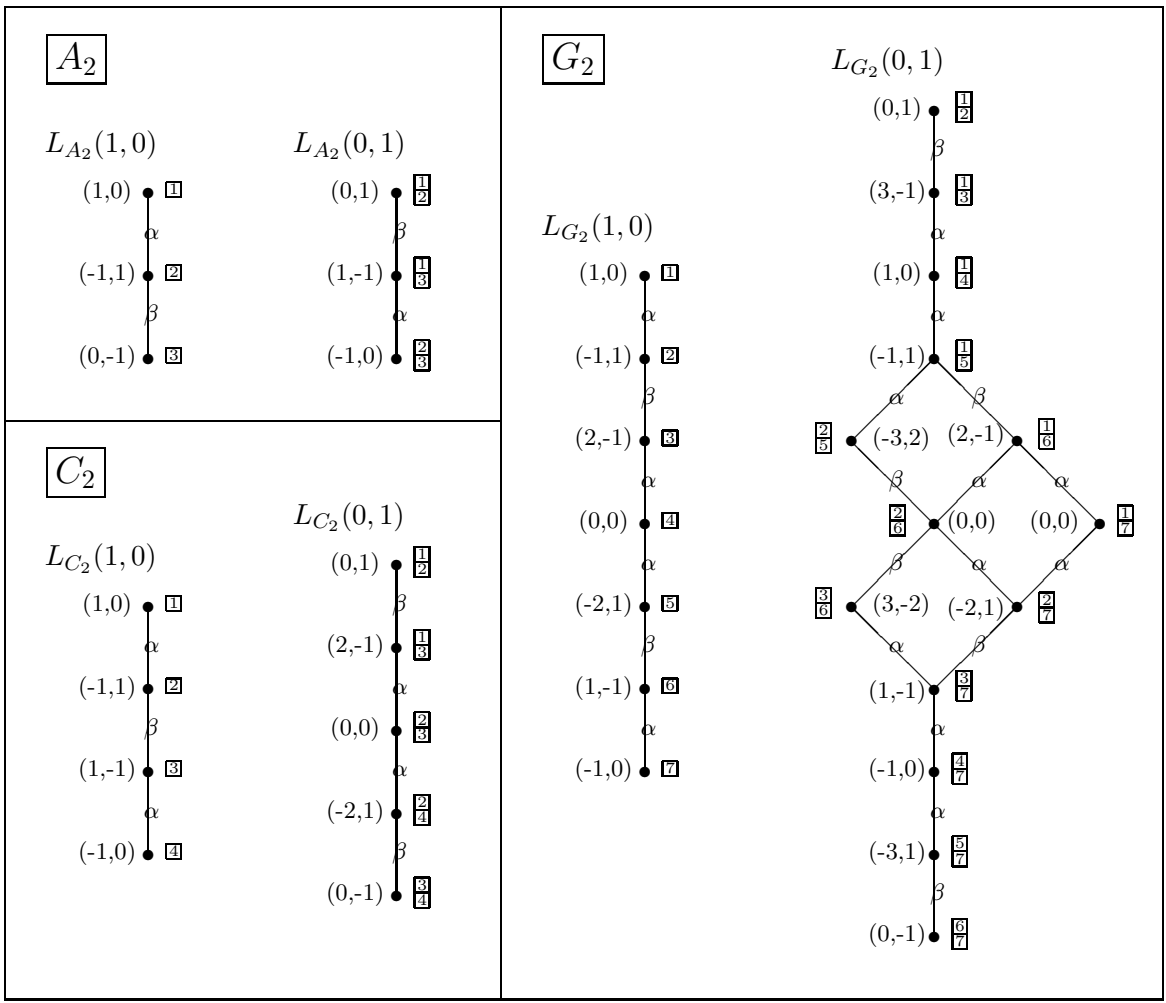
\put(2.25,12.4){\FundLatticeTableauxFigure: Weights and tableaux for 
$\mathfrak{g}$-fundamental  
lattices.}
\put(0,0){\line(0,1){12}}
\put(6,0){\line(0,1){12}}
\put(14,0){\line(0,1){12}}
\put(0,7){\line(1,0){6}}
\put(0,0){\line(1,0){14}}
\put(0,12){\line(1,0){14}}
\put(0.5,11.25){\fbox{\Large $A_{2}$}}
\put(0,6.75){\AtwoAlphaTableaux}
\put(3,6.75){\AtwoBetaTableaux}
\put(0.5,6.25){\fbox{\Large $C_{2}$}}
\put(0,0.75){\BtwoAlphaTableaux}
\put(3,0.25){\BtwoBetaTableaux}
\put(6.5,11.25){\fbox{\Large $G_{2}$}}
\put(6,1.75){\GtwoAlphaTableaux}
\put(8,-0.25){\GtwoBetaTableaux}
\end{picture}
\end{center} 
%\end{figure}

$\mathcal{S}_{C_{2}}(\lambda) = \Big\{$semistandard 
tableau $T$ of 
shape $\lambda$ with entries from $\{1,2,3,4\}\,\, \Big| $\\ 
\hspace*{2.25in}\lengthtwocolumn{1}{4}  
is not a column of $T$, and$\!\!$ \lengthtwocolumn{2}{3}  
appears at most once 
in $T \Big\}$

$\mathcal{S}_{G_{2}}(\lambda) = \Big\{$semistandard 
tableau $T$ of 
shape $\lambda$ with entries from $\{1,2,3,4,5,6,7\}\,\, \Big| $\\ 
\vspace*{0.1in}   
\hspace*{1in}the column$\!\!$ \lengthonecolumn{4} appears at most once in 
$T$;  
\lengthtwocolumn{2}{3},   
\lengthtwocolumn{2}{4},   
\lengthtwocolumn{3}{4},   
\lengthtwocolumn{3}{5},   
\lengthtwocolumn{4}{5},   
\lengthtwocolumn{4}{6}, and    
\lengthtwocolumn{5}{6}\\ 
\hspace*{1.05in}are not columns of $T$; plus the restrictions of 
\RestrictionChart\ $\Big\}$

\begin{figure}[h,t]
\begin{center}
\RestrictionChart: Some restrictions for any given 
$G_{2}$-semistandard tableau $T$.

\

\vspace*{-0.10in} 

\begin{tabular}{|c|c|}
\hline
\rule[-3mm]{0mm}{8mm}Column $T^{(i)}$ of $T$ & 
Then the succeeding column $T^{(i+1)}$ 
of $T$ cannot be\ldots \\
\hline
\hline
\rule[-2mm]{0mm}{8mm}\lengthonecolumn{4} & \lengthonecolumn{4}\\
\hline
\rule[-3mm]{0mm}{10mm}\lengthtwocolumn{1}{4} & 
\lengthonecolumn{1}, 
\lengthtwocolumn{1}{4}, \lengthtwocolumn{1}{5}, 
\lengthtwocolumn{1}{6}, \lengthtwocolumn{1}{7}\\
\hline
\rule[-3mm]{0mm}{10mm}\lengthtwocolumn{1}{5} & 
\lengthonecolumn{1}, 
\lengthtwocolumn{1}{5}, 
\lengthtwocolumn{1}{6}, \lengthtwocolumn{1}{7}\\
\hline
\rule[-3mm]{0mm}{10mm}\lengthtwocolumn{1}{6} & 
\lengthonecolumn{1}, \lengthonecolumn{2}, 
\lengthtwocolumn{1}{6}, \lengthtwocolumn{1}{7},
\lengthtwocolumn{2}{6}, \lengthtwocolumn{2}{7}\\
\hline
\rule[-3mm]{0mm}{10mm}\lengthtwocolumn{2}{6} & 
\lengthonecolumn{2}, 
\lengthtwocolumn{2}{6}, \lengthtwocolumn{2}{7}\\
\hline
\rule[-3mm]{0mm}{10mm}\lengthtwocolumn{1}{7} & 
\lengthonecolumn{1}, \lengthonecolumn{2}, 
\lengthonecolumn{3}, \lengthonecolumn{4}, 
\lengthtwocolumn{1}{7}, \lengthtwocolumn{2}{7},
\lengthtwocolumn{3}{7}, \lengthtwocolumn{4}{7}\\
\hline
\rule[-3mm]{0mm}{10mm}\lengthtwocolumn{2}{7} & 
\lengthonecolumn{2}, 
\lengthonecolumn{3}, \lengthonecolumn{4}, 
\lengthtwocolumn{2}{7},
\lengthtwocolumn{3}{7}, \lengthtwocolumn{4}{7}\\
\hline
\rule[-3mm]{0mm}{10mm}\lengthtwocolumn{3}{7} & 
\lengthonecolumn{3}, \lengthonecolumn{4}, 
\lengthtwocolumn{3}{7}, \lengthtwocolumn{4}{7}\\
\hline
\rule[-3mm]{0mm}{10mm}\lengthtwocolumn{4}{7} & 
\lengthonecolumn{4}, 
\lengthtwocolumn{4}{7}\\
\hline
\end{tabular}
\end{center}

\vspace*{-0.25in}
\end{figure}

{\em Proof.} The association of one-column $\mathfrak{g}$-semistandard tableaux with 
order ideals of $\mathfrak{g}$-fundamental posets is given in 
\IdealTableauFigureList.  Consider 
the $\mathfrak{g} = C_{2}$ case.  We want to show that the set 
$\mathcal{S}_{C_{2}}(\lambda)$ 
is the same as the 
stated set, which we denote $\mathcal{S}$. 
Let $T \in 
\mathcal{S}_{C_{2}}(\lambda)$, so $T = \mathbf{tableau}(\telt)$ for 
some order ideal $\telt$ of  
$\PBtwoba$.  Write $\PBtwoba = P_{1} \triangleleft 
\cdots \triangleleft P_{a+b}$, as depicted in 
\CaseBFig. Following \FundPosets, we label the vertices of $P_{j}$ as 
$w_{1,j}$, $w_{2,j}$, $w_{3,j}$, and $w_{4,j}$ with $w_{1,j} > 
w_{2,j} > w_{3,j} > w_{4,j}$ whenever $1 \leq j \leq b$, and   
we label the vertices of $P_{j}$ as 
$z_{1,j}$, $z_{2,j}$, and $z_{3,j}$ with $z_{1,j} > 
z_{2,j} > z_{3,j}$ whenever $1+b \leq j \leq a+b$.  
By definition, $T = (T^{(1)},\ldots,T^{(a+b)})$ 
with $T^{(i)} = \mathbf{tableau}(\telt \cap P_{i})$.   
The entries for $T^{(i)}$ are from the set $\{1,2,3,4\}$, and no 
$T^{(i)}$ is the column  
\lengthtwocolumn{1}{4}.  To see how the 
semistandard and other 
restrictions occur, suppose (for example) 
that $T^{(i)}$ is the column \lengthtwocolumn{2}{3} 
for some $1 \leq i \leq b$.  
Note that $\telt \cap P_{i} = \{w_{3,i}, w_{4,i}\}$.  It follows that 
$w_{1,j}$, $w_{2,j}$, and $w_{3,j}$ are not in $\telt$ for 
for $i < j \leq b$, and moreover $z_{1,j}$ is not in $\telt$ for 
$1+b \leq j \leq a+b$.  In particular, it follows that $T^{(i+1)}$ 
cannot be \lengthtwocolumn{1}{2}, \lengthtwocolumn{1}{3}, 
\lengthtwocolumn{2}{3}, or \lengthonecolumn{1}, assuming the column 
$T^{(i+1)}$ exists.  That is, the pair of 
columns $T^{(i)}$ and $T^{(i+1)}$ meets the requirements for 
inclusion in the set $\mathcal{S}$.  
The other eight cases for 
$T^{(i)}$ can be handled in a similar fashion.  
We conclude that $T \in \mathcal{S}$.  
So $\mathcal{S}_{C_{2}}(\lambda) \subseteq 
\mathcal{S}$.  

\begin{figure}[htb]
\begin{center}
\CaseBFig: $\PBtwoba = P_{1} \triangleleft 
\cdots \triangleleft P_{a+b}$

\setlength{\unitlength}{1cm}
\begin{picture}(10,10.25)
%Vertex labels
\put(1,0){\circle*{0.1}} \put(1.2,-0.1){\footnotesize $\beta$}
\put(0,1){\circle*{0.1}} \put(0.2,0.9){\footnotesize $\alpha$}
\put(1,2){\circle*{0.1}} \put(1.2,1.9){\footnotesize $\alpha$}
\put(0,3){\circle*{0.1}} \put(0.2,2.9){\footnotesize $\beta$}
\put(3,2){\circle*{0.1}} \put(3.2,1.9){\footnotesize $\beta$}
\put(2,3){\circle*{0.1}} \put(2.2,2.9){\footnotesize $\alpha$}
\put(3,4){\circle*{0.1}} \put(3.2,3.9){\footnotesize $\alpha$}
\put(2,5){\circle*{0.1}} \put(2.2,4.9){\footnotesize $\beta$}
\put(5,4){\circle*{0.1}} \put(5.2,3.9){\footnotesize $\beta$}
\put(4,5){\circle*{0.1}} \put(4.2,4.9){\footnotesize $\alpha$}
\put(5,6){\circle*{0.1}} \put(5.2,5.9){\footnotesize $\alpha$}
\put(4,7){\circle*{0.1}} \put(4.2,6.9){\footnotesize $\beta$}
\put(8,5){\circle*{0.1}} \put(8.2,4.9){\footnotesize $\alpha$}
\put(7,6){\circle*{0.1}} \put(7.2,5.9){\footnotesize $\beta$}
\put(6,7){\circle*{0.1}} \put(6.2,6.9){\footnotesize $\alpha$}
\put(9,6){\circle*{0.1}} \put(9.2,5.9){\footnotesize $\alpha$}
\put(8,7){\circle*{0.1}} \put(8.2,6.9){\footnotesize $\beta$}
\put(7,8){\circle*{0.1}} \put(7.2,7.9){\footnotesize $\alpha$}
\put(10,7){\circle*{0.1}} \put(10.2,6.9){\footnotesize $\alpha$}
\put(9,8){\circle*{0.1}} \put(9.2,7.9){\footnotesize $\beta$}
\put(8,9){\circle*{0.1}} \put(8.2,8.9){\footnotesize $\alpha$}
%Edges
\put(3,2){\line(1,1){5}} \put(2,3){\line(1,1){5}}
\put(2,5){\line(1,1){2}} \put(8,5){\line(1,1){1}}
\put(3,2){\line(-1,1){1}} \put(3,4){\line(-1,1){1}}
\put(5,4){\line(-1,1){1}} \put(5,6){\line(-1,1){1}}
\put(8,5){\line(-1,1){2}} \put(9,6){\line(-1,1){2}}
\put(1,0){\line(-1,1){1}} \put(0,1){\line(1,1){1}}
\put(1,2){\line(-1,1){1}} \put(10,7){\line(-1,1){2}}
\multiput(1.1,0.1375)(0.1,0.1){18}{.} \multiput(1.1,2.1375)(0.1,0.1){8}{.}
\multiput(0.1,3.1375)(0.1,0.1){18}{.} \multiput(9.1,6.1375)(0.1,0.1){8}{.}
\multiput(8.1,7.1375)(0.1,0.1){8}{.} \multiput(7.1,8.1375)(0.1,0.1){8}{.}
\put(-0.4,3.4){$P_{1}$} \put(1.5,5.4){$P_{b-1}$}
\put(3.6,7.4){$P_{b}$} \put(5.5,7.4){$P_{1+b}$}
\put(6.5,8.4){$P_{2+b}$} \put(7.5,9.4){$P_{a+b}$}
\end{picture}
\end{center}

\vspace*{-0.25in}
\end{figure}

In the other direction, suppose $T = (T^{(1)},\ldots,T^{(a+b)})$ is in 
$\mathcal{S}$.  For each $i$, let $Q_{i}$ be the order ideal of $P_{i}$ 
corresponding to the one-column tableau 
$T^{(i)}$, and let $\telt := \cup_{i} Q_{i}$.  By 
examining cases as in the previous paragraph, one can check that the 
restrictions on $T$ as an element of $\mathcal{S}$ guarantee that 
$\telt$ will be an order ideal of $\PBtwoba$ with $Q_{i} = \telt \cap 
P_{i}$ for each $i$.  Hence $T \in \mathcal{S}_{C_{2}}(\lambda)$.  
It follows that $\mathcal{S} \subseteq 
\mathcal{S}_{C_{2}}(\lambda)$, which completes the proof for the 
$C_{2}$ case.  The $A_{2}$ and $G_{2}$  
cases can be handled by similar arguments.\hfill\QED

\noindent
{\bf \OrderRemark}\ \ \ In passing we note that 
the partial ordering and the covering relations in 
$L_{\mathfrak{g}}^{\beta\alpha}(\lambda)$ are easy to describe with 
the ``coordinates'' of $\mathfrak{g}$-semistandard 
tableaux.  For 
$\selt$ and $\telt$ in $L_{\mathfrak{g}}^{\beta\alpha}(\lambda)$, let 
$S := \mathbf{tableau}(\selt)$ and $T := \mathbf{tableau}(\telt)$.  
Then $\selt \leq \telt$ if and only if $S^{(i)}_{j} \geq T^{(i)}_{j}$ for 
all $i, j$.  (This is the ``reverse componentwise'' order on 
tableaux.)  Moreover, $\selt \rightarrow \telt$ is a covering relation 
in the poset $L_{\mathfrak{g}}^{\beta\alpha}(\lambda)$ if and only if 
for some $i$ and $j$ we have $S^{(i)}_{j} = T^{(i)}_{j} + 1$ while 
$S^{(p)}_{q} = T^{(p)}_{q}$ for all $(p,q) \not= (i,j)$.  
For 
$\mathfrak{g} = A_{2}$, 
the edge gets color $\alpha$ if $T^{(i)}_{j}$ 
is 1 and color $\beta$ if $T^{(i)}_{j}$ is 2; 
for $\mathfrak{g} = C_{2}$, 
the edge gets color $\alpha$ if $T^{(i)}_{j}$ 
is 1 or 3 and color $\beta$ if $T^{(i)}_{j}$ is 2; 
for $\mathfrak{g} = G_{2}$, 
the edge gets color $\alpha$ if $T^{(i)}_{j}$ 
is 1 or 3 or 4 or 6 and color $\beta$ if $T^{(i)}_{j}$ is 2 or 5.  
\hfill\QED 

For a tableau $T$ of shape $\lambda = (a,b)$ and a positive 
integer $k$, we define $n_{k}(T)$ to be the number of times the entry 
$k$ appears in the tableau $T$. Observe that $\displaystyle n_{k}(T) = 
\sum_{i=1}^{a+b} n_{k}(T^{(i)})$. 
Define a function $\Twt: \TabSet 
\longrightarrow \mathbb{Z} \times \mathbb{Z}$ by the rules: 

{\small
\[\Twt(T) := \left\{\begin{array}{ll}
\Big(n_{1}(T) - n_{2}(T), n_{2}(T) - 
n_{3}(T)\Big) & \mbox{ if } \mathfrak{g} = A_{2}\\
\Big(n_{1}(T) - n_{2}(T) + n_{3}(T) - 
n_{4}(T), n_{2}(T) - 
n_{3}(T)\Big) & \mbox{ if } \mathfrak{g} = C_{2}\\
\Big(n_{1}(T) - n_{2}(T) + 2n_{3}(T) - 
2n_{5}(T) + n_{6}(T) - n_{7}(T),  & \\ \hspace*{1in}  
n_{2}(T) - 
n_{3}(T) + n_{5}(T) - n_{6}(T)\Big) & 
\mbox{ if } \mathfrak{g} = G_{2} 
\end{array}
\right.\]
}

The function  $wt(\selt)$  defined on  $\Lba$  in terms of the color components 
of  $\Lba$  can be expressed in terms of the tableau entry counts when 
the elements  $\selt$  of  $\Lba$  are viewed as tableaux  $\telt$ in 
$\TabSet$:

\noindent
{\bf \TableauWeightProp}\ \ {\sl Let $\lambda = (a,b)$, with $a, b \geq 0$. 
For $\telt \in L_{\mathfrak{g}}^{\beta\alpha}(\lambda)$, 
consider $T := \mathbf{tableau}(\telt) \in \TabSet$. Then $wt(\telt) = \Twt(T)$.}

{\em Proof.} With the help of \FigsForTProp, 
one can easily confirm 
the result by hand whenever $\lambda$ is a fundamental weight.  
Then more generally apply  \WeightsLemma\ to $wt(\telt)$, 
noting that $\Twt(T) = 
\displaystyle \sum_{i=1}^{a+b}\Twt(T^{(i)})$.\hfill\QED

This concludes our self-contained development of 
$\mathfrak{g}$-semistandard posets,  
$\mathfrak{g}$-semistandard lattices, and  
$\mathfrak{g}$-semistandard tableaux in Sections \GridNum\ and 
\SemiNum.

%==================================================================
%\newpage
\vspace{2ex} 

\noindent
{\bf \Large \CharNum.\ \ Weyl characters; Littelmann's tableaux; 
main results}

Our main result, \CharacterProposition, expresses the 
Weyl characters for the irreducible representations of the 
rank two semisimple Lie algebras as 
generating functions for $\mathfrak{g}$-semistandard lattices.  
We begin 
by recording some 
explicit data on roots, weights, Weyl groups, and irreducible characters 
for the rank two semisimple Lie algebras. 
Then we describe certain 
tableaux obtained by Littelmann in \cite{Lit1}, and in 
\WeightPresBijection\ we match these with our 
$\mathfrak{g}$-semistandard tableaux. \CharacterCorollary\ gives 
the product expressions 
for the rank generating functions. 
 
In rank two we denote the elements $e_{\omega_{\alpha}}$ 
and $e_{\omega_{\beta}}$ of the group ring $\mathbb{Z}[\Lambda]$ 
by $x$ and $y$. 
%At times write $(x,y)^{(r,t)}$ for $x^{r}y^{t}$. 
The reflections 
$s_{\alpha}$ and $s_{\beta}$ in ${W}$ act on the fundamental weights 
as follows: 
$s_{\alpha}\omega_{\alpha} = \omega_{\alpha} - \alpha$, 
$s_{\alpha}\omega_{\beta} = \omega_{\beta}$, 
$s_{\beta}\omega_{\alpha} = \omega_{\alpha}$, and 
$s_{\beta}\omega_{\beta} = \omega_{\beta} - \beta$.  
\RankTwoData\ has data for the simple roots, positive roots, and Weyl group 
for each of the rank two simple Lie algebras. 
Recall from Section \DefsNum\ that the denominator $A_{\varrho}$ of the Weyl 
character formula 
can be expressed as a product over the positive roots. Also 
recall that the numerator $A_{\varrho + \lambda}$ is 
an alternating sum over the elements of the Weyl group.  
Using the data of \RankTwoData\ one obtains for 
$A_{2}$ 
\begin{eqnarray*}
A_{\varrho} & = & xy(1-x^{-2}y)(1-xy^{-2})(1-x^{-1}y^{-1})\\
 & = & xy - x^{-1}y^{2}-x^{2}y^{-1}+x^{-2}y+xy^{-2}-x^{-1}y^{-1}\\
A_{\varrho+\lambda} & = & 
x^{a+1}y^{b+1}-x^{-(a+1)}y^{a+b+2}-x^{a+b+2}y^{-(b+1)}\\
 & & +x^{-(a+b+2)}y^{a+1}+x^{b+1}y^{-(a+b+2)}-x^{-(b+1)}y^{-(a+1)}
\end{eqnarray*}
For $C_{2}$ we get:
\begin{eqnarray*}
A_{\varrho} & = & xy(1-x^{-2}y)(1-x^{2}y^{-2})(1-y^{-1})(1-x^{-2})\\
 & = & xy-x^{-1}y^{2}-x^{3}y^{-1}+x^{-3}y^{2}+x^{3}y^{-2}-x^{-3}y^{1}-xy^{-2}+x^{-1}y^{-1}\\
A_{\varrho+\lambda} & = & 
x^{a+1}y^{b+1}-x^{-(a+1)}y^{a+b+2}-x^{a+2b+3}y^{-(b+1)}+x^{-(a+2b+3)}y^{a+b+2}\\
& & +x^{a+2b+3}y^{-(a+b+2)}-x^{-(a+2b+3)}y^{b+1}-x^{a+1}y^{-(a+b+2)}+x^{-(a+1)}y^{-(b+1)}
\end{eqnarray*}
And for $G_{2}$ we have: 
\begin{eqnarray*}
A_{\varrho} & = & 
xy(1-x^{-2}y)(1-x^{3}y^{-2})(1-xy^{-1})(1-x^{-1})(1-x^{-3}y)(1-y^{-1})\\
 & = & 
xy-x^{-1}y^{2}-x^{4}y^{-1}+x^{-4}y^{3}+x^{5}y^{-2}-x^{-5}y^{3}-x^{5}y^{-3}+x^{-5}y^{2}+x^{4}y^{-3}\\
& & -x^{-4}y-xy^{-2}+x^{-1}y^{-1}\\
A_{\varrho+\lambda} & = & 
x^{a+1}y^{b+1}-x^{-(a+1)}y^{a+b+2}-x^{a+3b+4}y^{-(b+1)}+x^{-(a+3b+4)}y^{a+2b+3}+x^{2a+3b+5}y^{-(a+b+2)}\\
& & 
-x^{-(2a+3b+5)}y^{a+2b+3}-x^{2a+3b+5}y^{-(a+2b+3)}+x^{-(2a+3b+5)}y^{a+b+2}+x^{a+3b+4}y^{-(a+2b+3)}\\
& & 
-x^{-(a+3b+4)}y^{b+1}-x^{a+1}y^{-(a+b+2)}+x^{-(a+1)}y^{-(b+1)}
\end{eqnarray*}

\begin{figure}[t]
\begin{center}
{\RankTwoData: Roots and Weyl groups for the rank two simple Lie 
algebras.}

\

\vspace*{-0.10in}

\begin{tabular}{|c|c|c|c|}
\hline
Algebra & Simple roots & Positive roots & \parbox{3in}{
\begin{center}
Weyl group $W$ 

\vspace*{-0.05in}

{\scriptsize (By generators and relations; as reduced words)}
\end{center}}\\
\hline
\hline
$A_{2}$ 
& \parbox{1.15in}{\begin{center}
$\alpha = 2\omega_{\alpha}-\omega_{\beta}$ 

$\beta = -\omega_{\alpha}+2\omega_{\beta}$
\end{center}} 
& $\alpha, \beta, \alpha+\beta$ 
& \rule[-7mm]{0mm}{16mm} 
\parbox{3in}{\begin{center}
$\langle s_{\alpha}, s_{\beta} | s_{\alpha}^{2} = 
s_{\beta}^{2} = id, (s_{\alpha}s_{\beta})^{3} = id \rangle$ 

$\{id, s_{\alpha}, s_{\beta}, s_{\alpha}s_{\beta}, s_{\beta}s_{\alpha}, 
s_{\alpha}s_{\beta}s_{\alpha} = s_{\beta}s_{\alpha}s_{\beta}\}$
\end{center}}\\
\hline 
$C_{2}$ 
& \parbox{1.15in}{\begin{center}
$\alpha = 2\omega_{\alpha}-\omega_{\beta}$ 

$\beta = -2\omega_{\alpha}+2\omega_{\beta}$
\end{center}} 
& \parbox{0.9in}{\begin{center}
$\alpha, \beta, \alpha+\beta,$ 

$2\alpha+\beta$ 
\end{center}} 
& \rule[-7mm]{0mm}{16mm} 
\parbox{3in}{\begin{center}
$\langle s_{\alpha}, s_{\beta} | s_{\alpha}^{2} = 
s_{\beta}^{2} = id, (s_{\alpha}s_{\beta})^{4} = id \rangle$ 

$\{id, s_{\alpha}, s_{\beta}, s_{\alpha}s_{\beta}, s_{\beta}s_{\alpha}, 
s_{\alpha}s_{\beta}s_{\alpha}, s_{\beta}s_{\alpha}s_{\beta},$

$s_{\alpha}s_{\beta}s_{\alpha}s_{\beta} = 
s_{\beta}s_{\alpha}s_{\beta}s_{\alpha}\}$
\end{center}}\\
\hline 
$G_{2}$ 
& \parbox{1.15in}{\begin{center}
$\alpha = 2\omega_{\alpha}-\omega_{\beta}$ 

$\beta = -3\omega_{\alpha}+2\omega_{\beta}$
\end{center}} 
& \parbox{0.9in}{\begin{center}
$\alpha, \beta, \alpha+\beta,$ 

$2\alpha+\beta,$ 

$3\alpha+\beta,$ 

$3\alpha+2\beta$ 
\end{center}} 
& \rule[-7mm]{0mm}{16mm} 
\parbox{3in}{\begin{center}
$\langle s_{\alpha}, s_{\beta} | s_{\alpha}^{2} = 
s_{\beta}^{2} = id, (s_{\alpha}s_{\beta})^{6} = id \rangle$ 

$\{id, s_{\alpha}, s_{\beta}, s_{\alpha}s_{\beta}, s_{\beta}s_{\alpha}, 
s_{\alpha}s_{\beta}s_{\alpha}, s_{\beta}s_{\alpha}s_{\beta},$

$s_{\alpha}s_{\beta}s_{\alpha}s_{\beta},  
s_{\beta}s_{\alpha}s_{\beta}s_{\alpha}, 
s_{\alpha}s_{\beta}s_{\alpha}s_{\beta}s_{\alpha},$

$s_{\beta}s_{\alpha}s_{\beta}s_{\alpha}s_{\beta}, 
s_{\alpha}s_{\beta}s_{\alpha}s_{\beta}s_{\alpha}s_{\beta} = 
s_{\beta}s_{\alpha}s_{\beta}s_{\alpha}s_{\beta}s_{\alpha}\}$
\end{center}}\\
\hline 
\end{tabular}
\end{center}

\vspace*{-0.25in}
\end{figure} 

We now seek a correspondence between our $\mathfrak{g}$-semistandard tableaux 
and certain tableaux of Littelmann \cite{Lit1}.   
Littelmann's tableaux are ``translations'' of the standard monomial 
theory tableaux of 
Lakshmibai and Seshadri.   
The roles of his columns and rows are reversed with respect 
to this paper.  We pre-process Littelmann's tableaux in two steps.  
First, we reflect them across the main diagonal $i = j$.  Then we 
group $k$ of his columns at a time into 
a ``block'' of $k$
columns, where $k=1$ for $A_{2}$, $k=2$ for $C_{2}$, and $k=6$ for $G_{2}$. 
We define $\mathbf{shape}(k \times \lambda) := \mathbf{shape}(\mu)$, 
where $\mu = ka\omega_{\alpha} + kb\omega_{\beta} = (ka,kb)$.  A 
{\em $k$-tableau of shape} $\lambda$ is a filling of  
$\mathbf{shape}(k \times \lambda)$ with entries from some totally 
ordered set. The {\em semistandard} condition on $k$-tableaux 
is the same as the semistandard condition of Section \SemiNum.  
For a $k$-tableau $T$ of shape $\lambda$, we 
write $T = (T^{(1)},\ldots,T^{(a+b)})$, where $T^{(i)}$ is the $i$th 
block of $k$ columns of $T$ counting from the left.  
Only certain
fillings of these $k$-column blocks will be ``admissible''.  
Here are our processed versions of Littelmann's tableaux:

\noindent
{\bf \LitTableauProp}\ \ {Let $\lambda = (a,b)$, with $a, b \geq 0$.  Then:} 

$\mathcal{LT}_{A_{2}}(\lambda) := \Big\{$semistandard  
$1$-tableau $T$ of 
shape $\lambda$ with entries from $\{1,2,3\}\,\, \Big| $\\
\hspace*{1.5in}admissible 1-column blocks of $T$ come from \LitAdmissibleColumns\ $\Big\}$

$\mathcal{LT}_{C_{2}}(\lambda) := \Big\{$semistandard  
$2$-tableau $T$ of 
shape $\lambda$ with entries from $\{1,2,3,4\}\,\, \Big| $\\ 
\hspace*{1.5in}admissible 2-column blocks of $T$ come from \LitAdmissibleColumns\ $\Big\}$

$\mathcal{LT}_{G_{2}}(\lambda) := \Big\{$semistandard  
$6$-tableau $T$ of 
shape $\lambda$ with entries from $\{1,2,3,4,5,6\}\,\, \Big| $\\ 
\hspace*{1.5in}admissible 6-column blocks of $T$ come from \LitAdmissibleColumns\ $\Big\}$

\noindent
{Moreover, the weight $wt_{Lit}(T)$ of a Littelmann tableau $T$ is 
given by} 
{\small
\[wt_{Lit}(T) := \left\{\begin{array}{cl}
\Big(n_{1}(T) - n_{2}(T)\Big)\omega_{\alpha} + \Big(n_{2}(T) - 
n_{3}(T)\Big)\omega_{\beta} & \mbox{ if } \mathfrak{g} = A_{2}\\
\frac{1}{2}\Big[\Big(n_{1}(T) - n_{2}(T) + n_{3}(T) - 
n_{4}(T)\Big)\omega_{\alpha} + \Big(n_{2}(T) - 
n_{3}(T)\Big)\omega_{\beta}\Big] & \mbox{ if } \mathfrak{g} = C_{2}\\
\frac{1}{6}\Big[\Big(n_{1}(T) - n_{2}(T) + 2n_{3}(T) - 
2n_{4}(T) + n_{5}(T) - n_{6}(T)\Big)\omega_{\alpha}  & \\ \hspace*{1.42in}+ 
\Big(n_{2}(T) - 
n_{3}(T) + n_{4}(T) - n_{5}(T)\Big)\omega_{\beta}\Big] & 
\mbox{ if } \mathfrak{g} = G_{2}
\end{array}
\right.\]
}

To obtain these tableaux  
for $A_{2}$, see Section 2 of \cite{Lit1}; for $C_{2}$, 
see the Appendix of \cite{Lit1}; and for $G_{2}$ see Section 3 of 
that paper. Littelmann expresses his weight function in 
terms of a basis $\{\varepsilon_{1},\varepsilon_{2}\}$ for $\Lambda$, where 
$\varepsilon_{1} = \omega_{\alpha}$ and $\varepsilon_{2} = 
\omega_{\beta} - \omega_{\alpha}$. 
A consequence of standard monomial theory is: 

\noindent
{\bf Theorem (Littelmann, Lakshmibai, Seshadri)}\ \ 
{\sl Let $\lambda = (a,b)$, with $a, b \geq 0$.  
Let $\mathfrak{g}$ be a rank two simple Lie algebra. Then 
$(\mathcal{LT}_{\mathfrak{g}}(\lambda), wt_{Lit})$ is a splitting 
system for the irreducible character $\chi_{_{\lambda}}$.}
  
Next we describe a 
weight-preserving bijection $\phi : \TabSet \longrightarrow 
\mathcal{LT}_{\mathfrak{g}}(\lambda)$.  For fundamental 
weights, the 
correspondence between the $\mathfrak{g}$-semistandard tableaux of Section \SemiNum\ 
and Littelmann $k$-tableaux of this section is given in 
\LitAdmissibleColumns.  Given a rank two simple Lie algebra 
$\mathfrak{g}$, a dominant weight $\lambda = 
a\omega_{\alpha}+b\omega_{\beta} = (a,b)$, and a 
tableau $T$ in $\mathcal{S}_{\mathfrak{g}}(\lambda)$, we 
let $U = \phi(T)$ be the Littelmann $k$-tableau of shape $\lambda$ whose $i$th 
$k$-column block $U^{(i)}$ corresponds to the $i$th column $T^{(i)}$ 
of $T$. Keeping in mind the restrictions on which columns can follow 
$T^{(i)}$ to form a $\mathfrak{g}$-semistandard tableau $T$ in 
$\mathcal{S}_{\mathfrak{g}}(\lambda)$, one can check that $U^{(i)}$ 
followed by $U^{(i+1)}$ 
obeys the semistandard requirement for Littelmann $k$-tableaux.  Hence 
$U$ is in $\mathcal{LT}_{\mathfrak{g}}(\lambda)$.  Similarly, given 
$U$ in $\mathcal{LT}_{\mathfrak{g}}(\lambda)$, let $T = \psi(U)$ be 
the tableau of shape $\lambda$ whose $i$th column $T^{(i)}$ 
corresponds to the $i$th $k$-column block $U^{(i)}$ of $U$. 
Keeping in mind the semistandard condition on the Littelmann 
$k$-tableaux in $\mathcal{LT}_{\mathfrak{g}}(\lambda)$, one can check that 
$T^{(i)}$ followed by 
$T^{(i+1)}$ obeys the restrictions for $\mathfrak{g}$-semistandard tableaux in 
$\mathcal{S}_{\mathfrak{g}}(\lambda)$, and hence $T$ is in 
$\mathcal{S}_{\mathfrak{g}}(\lambda)$.  Clearly the mappings $\phi$ 
and $\psi$ are inverses. 

\noindent
{\bf \WeightPresBijection}\ \ {\sl Keep the notation of the 
previous paragraph.  The mapping $\phi : \TabSet \longrightarrow 
\mathcal{LT}_{\mathfrak{g}}(\lambda)$ described above is a 
weight-preserving bijection:  
for any $T \in \TabSet$, $wt_{Lit}(\phi(T)) = \Twt(T)$.}

{\em Proof.}  
We must check that $\phi$ is weight-preserving. 
If $\lambda$ is a fundamental weight, simply inspect 
\LitAdmissibleColumns.  If $\lambda$ is a dominant weight and $T$ is 
in $\TabSet$, then $\Twt(T) = \sum \Twt(T^{(i)}) = \sum 
wt_{Lit}(\phi(T^{(i)}))$.  The characterization of 
$wt_{Lit}$ in \LitTableauProp\ implies that 
$wt_{Lit}(\phi(T)) = \sum wt_{Lit}(\phi(T^{(i)}))$.\hfill\QED

\noindent
{\bf \CharacterProposition}\ \ {\sl 
Let $\mathfrak{g}$ be a semisimple Lie algebra 
of rank two. 
Let $\lambda = (a,b)$, with $a, b \geq 0$.  
Let $L$ be one of the 
$\mathfrak{g}$-semistandard lattices $\Lba$ or $\Lab$. 
Then $L$ is a splitting poset for an 
irreducible representation of $\mathfrak{g}$ with highest weight 
$\lambda$.  In particular,}  
\[char_{\mathfrak{g}}(\lambda; x,y) = \sum_{\selt \in L}(x,y)^{wt(\selt)}\] 

{\em Proof.} \StructureResult\ states that $L$ satisfies the 
$\mathfrak{g}$-structure condition.  Suppose $\mathfrak{g}$ is simple.  Since 
$(\mathcal{LT}_{\mathfrak{g}}(\lambda), wt_{Lit})$ is a splitting 
system for $\chi_{_{\lambda}}$, it follows from \WeightPresBijection\ 
that $(\TabSet, \Twt)$ is as well.  
From \TableauWeightProp\ it now follows that 
$(\Lba,wt)$ is a splitting system for $\chi_{_{\lambda}}$. Since 
%\noindent 
$\Lab \cong (\Lba)^{\triangle}$, then by  
\DualCharacterLemma\ the result holds for $\Lab$ as well. 
The case $A_{1} \oplus A_{1}$ can be handled by constructing the 
corresponding representation. 
\hfill\QED 

The main results of \cite{McClard} and \cite{Alv} were closely 
related to \CharacterProposition\ for the cases of $G_{2}$ and $C_{2}$ 
respectively.  For these rank two simple Lie algebras $\mathfrak{g}$, 
the lattices $\Lba$ were obtained by taking natural partial orders on 
the corresponding $\mathfrak{g}$-semistandard tableaux of Section 
\SemiNum, and case analysis arguments were used to show that the 
mapping $\phi$ preserves weights and that the 
$\mathfrak{g}$-structure condition is satisfied.  However, 
$\mathfrak{g}$-semistandard posets did not arise in their approach.  If 
one is willing to depend entirely upon \cite{Lit1}, then in this 
manner one can obtain \SectFivePropsList\ 
%=============================
% Begin pagebreak for figure
%=============================
\newpage 
%======================
%\begin{figure}[ht]
\begin{center}
{\LitAdmissibleColumns: Admissible $k$-column blocks for 
Littelmann tableau, their weights,\\ and their corresponding 
$\mathfrak{g}$-semistandard columns.}

\

\vspace*{-0.10in}

\begin{tabular}{|c|c|c|}
\hline
\parbox{0.75in}{\small \begin{center}$A_{2}$ Admissible $1$-block $T$\end{center}} 
& \parbox{0.6in}{\small \begin{center}Weight $wt_{Lit}(T)$\end{center}} 
& \parbox{1.05in}{\small \begin{center}Corresponding $A_{2}$-semistandard 
tableau\end{center}}\\
\hline
\rule[-3mm]{0mm}{8mm} \lengthonecolumnA{1} 
& $\omega_{\alpha}$ 
& \lengthonecolumnA{1}\\
\hline
\rule[-3mm]{0mm}{8mm} \lengthonecolumnA{2} 
& $-\omega_{\alpha}+\omega_{\beta}$ 
& \lengthonecolumnA{2}\\
\hline
\rule[-3mm]{0mm}{8mm} \lengthonecolumnA{3} 
& $-\omega_{\beta}$ 
& \lengthonecolumnA{3}\\
\hline
\hline
\rule[-5mm]{0mm}{12mm} \lengthtwocolumnA{1}{2} 
& $\omega_{\beta}$ 
& \lengthtwocolumnA{1}{2}\\
\hline
\rule[-5mm]{0mm}{12mm} \lengthtwocolumnA{1}{3} 
& $\omega_{\alpha}-\omega_{\beta}$ 
& \lengthtwocolumnA{1}{3}\\
\hline
\rule[-5mm]{0mm}{12mm} \lengthtwocolumnA{2}{3} 
& $-\omega_{\alpha}$ 
& \lengthtwocolumnA{2}{3}\\
\hline
\end{tabular}
\begin{tabular}{|c|c|c|}
\hline
%======================
\parbox{0.75in}{\small \begin{center}$C_{2}$ Admissible $2$-block $T$\end{center}} 
& \parbox{0.6in}{\small \begin{center}Weight $wt_{Lit}(T)$\end{center}} 
& \parbox{1.05in}{\small \begin{center}Corresponding $C_{2}$-semistandard 
tableau\end{center}}\\
\hline
\rule[-3mm]{0mm}{8mm} \lengthonecolumnB{1} 
& $\omega_{\alpha}$ 
& \lengthonecolumnA{1}\\
\hline
\rule[-3mm]{0mm}{8mm} \lengthonecolumnB{2} 
& $-\omega_{\alpha}+\omega_{\beta}$ 
& \lengthonecolumnA{2}\\
\hline
\rule[-3mm]{0mm}{8mm} \lengthonecolumnB{3} 
& $\omega_{\alpha}-\omega_{\beta}$ 
& \lengthonecolumnA{3}\\
\hline
\rule[-3mm]{0mm}{8mm} \lengthonecolumnB{4} 
& $-\omega_{\alpha}$ 
& \lengthonecolumnA{4}\\
\hline
\hline
\rule[-5mm]{0mm}{12mm} \lengthtwocolumnB{1}{2}{1}{2} 
& $\omega_{\beta}$ 
& \lengthtwocolumnA{1}{2}\\
\hline
\rule[-5mm]{0mm}{12mm} \lengthtwocolumnB{1}{3}{1}{3} 
& $2\omega_{\alpha}-\omega_{\beta}$ 
& \lengthtwocolumnA{1}{3}\\
\hline
\rule[-5mm]{0mm}{12mm} \lengthtwocolumnB{1}{3}{2}{4} 
& $0\omega_{\alpha} + 0\omega_{\beta}$ 
& \lengthtwocolumnA{2}{3}\\
\hline
\rule[-5mm]{0mm}{12mm} \lengthtwocolumnB{2}{4}{2}{4} 
& $-2\omega_{\alpha}+\omega_{\beta}$ 
& \lengthtwocolumnA{2}{4}\\
\hline
\rule[-5mm]{0mm}{12mm} \lengthtwocolumnB{3}{4}{3}{4} 
& $-\omega_{\beta}$ 
& \lengthtwocolumnA{3}{4}\\
\hline
\end{tabular}

\vspace*{0.25in}
\begin{tabular}{|c|c|c|}
\hline
\parbox{0.75in}{\small \begin{center}$G_{2}$ Admissible $6$-block $T$\end{center}} 
& \parbox{0.6in}{\small \begin{center}Weight $wt_{Lit}(T)$\end{center}} 
& \parbox{1.05in}{\small \begin{center}Corresponding $G_{2}$-semistandard 
tableau\end{center}}\\
\hline
\rule[-3mm]{0mm}{8mm} \lengthonecolumnG{1}{1} 
& $\omega_{\alpha}$ 
& \lengthonecolumnA{1}\\
\hline
\rule[-3mm]{0mm}{8mm} \lengthonecolumnG{2}{2} 
& $-\omega_{\alpha}+\omega_{\beta}$ 
& \lengthonecolumnA{2}\\
\hline
\rule[-3mm]{0mm}{8mm} \lengthonecolumnG{3}{3} 
& $2\omega_{\alpha}-\omega_{\beta}$ 
& \lengthonecolumnA{3}\\
\hline
\rule[-3mm]{0mm}{8mm} \lengthonecolumnG{3}{4} 
& $0\omega_{\alpha}+0\omega_{\beta}$ 
& \lengthonecolumnA{4}\\
\hline
\rule[-3mm]{0mm}{8mm} \lengthonecolumnG{4}{4} 
& $-2\omega_{\alpha}+\omega_{\beta}$ 
& \lengthonecolumnA{5}\\
\hline
\rule[-3mm]{0mm}{8mm} \lengthonecolumnG{5}{5} 
& $\omega_{\alpha}-\omega_{\beta}$ 
& \lengthonecolumnA{6}\\
\hline
\rule[-3mm]{0mm}{8mm} \lengthonecolumnG{6}{6} 
& $-\omega_{\alpha}$ 
& \lengthonecolumnA{7}\\
\hline
\end{tabular}
\end{center}
%\end{figure} 

%=============================
% End pagebreak for figure
%=============================
\noindent 
from \WeightPresBijection\ and 
Littelmann's analog to \CharacterProposition\ without using 
\WeightsLemma.  But this approach would take at least as much 
(related) work and would not be as uniformly stated. 

%=============================
% Begin pagebreak for figure
%=============================
\newpage
%===================
%\begin{figure}[ht]
\begin{center}
{\LitAdmissibleColumns\ (continued): Admissible $k$-column blocks for 
Littelmann tableaux, their weights,\\ and their corresponding 
$\mathfrak{g}$-semistandard columns.}

\ 

\vspace*{-0.10in}

\begin{tabular}{|c|c|c|}
\hline
\parbox{0.75in}{\small \begin{center}$G_{2}$ Admissible $6$-block $T$\end{center}} 
& \parbox{0.6in}{\small \begin{center}Weight $wt_{Lit}(T)$\end{center}} 
& \parbox{1.05in}{\small \begin{center}Corresponding $G_{2}$-semistandard 
tableau\end{center}}\\
\hline
\rule[-5mm]{0mm}{12mm} \lengthtwocolumnG{1}{2}{1}{2}{1}{2}{1}{2} 
& $\omega_{\beta}$ 
& \lengthtwocolumnA{1}{2}\\
\hline
\rule[-5mm]{0mm}{12mm} \lengthtwocolumnG{1}{3}{1}{3}{1}{3}{1}{3} 
& $3\omega_{\alpha}-\omega_{\beta}$ 
& \lengthtwocolumnA{1}{3}\\
\hline
\rule[-5mm]{0mm}{12mm} \lengthtwocolumnG{1}{3}{1}{3}{1}{3}{2}{4} 
& $\omega_{\alpha}$ 
& \lengthtwocolumnA{1}{4}\\
\hline
\rule[-5mm]{0mm}{12mm} \lengthtwocolumnG{1}{3}{2}{4}{2}{4}{2}{4} 
& $-\omega_{\alpha}+\omega_{\beta}$
& \lengthtwocolumnA{1}{5}\\
\hline
\rule[-5mm]{0mm}{12mm} \lengthtwocolumnG{2}{4}{2}{4}{2}{4}{2}{4} 
& $-3\omega_{\alpha}+2\omega_{\beta}$ 
& \lengthtwocolumnA{2}{5}\\
\hline
\rule[-5mm]{0mm}{12mm} \lengthtwocolumnG{1}{3}{2}{4}{3}{5}{3}{5} 
& $2\omega_{\alpha}-\omega_{\beta}$ 
& \lengthtwocolumnA{1}{6}\\
\hline
\rule[-5mm]{0mm}{12mm} \lengthtwocolumnG{2}{4}{2}{4}{3}{5}{3}{5} 
& $0\omega_{\alpha}+0\omega_{\beta}$ 
& \lengthtwocolumnA{2}{6}\\
\hline
\rule[-5mm]{0mm}{12mm} \lengthtwocolumnG{1}{3}{2}{4}{3}{5}{4}{6} 
& $0\omega_{\alpha}+0\omega_{\beta}$
& \lengthtwocolumnA{1}{7}\\
\hline
\rule[-5mm]{0mm}{12mm} \lengthtwocolumnG{3}{5}{3}{5}{3}{5}{3}{5} 
& $3\omega_{\alpha}-2\omega_{\beta}$
& \lengthtwocolumnA{3}{6}\\
\hline
\rule[-5mm]{0mm}{12mm} \lengthtwocolumnG{2}{4}{2}{4}{3}{5}{4}{6} 
& $-2\omega_{\alpha}+\omega_{\beta}$
& \lengthtwocolumnA{2}{7}\\
\hline
\rule[-5mm]{0mm}{12mm} \lengthtwocolumnG{3}{5}{3}{5}{3}{5}{4}{6} 
& $\omega_{\alpha}-\omega_{\beta}$
& \lengthtwocolumnA{3}{7}\\
\hline
\rule[-5mm]{0mm}{12mm} \lengthtwocolumnG{3}{5}{4}{6}{4}{6}{4}{6} 
& $-\omega_{\alpha}$
& \lengthtwocolumnA{4}{7}\\
\hline
\rule[-5mm]{0mm}{12mm} \lengthtwocolumnG{4}{6}{4}{6}{4}{6}{4}{6} 
& $-3\omega_{\alpha}+\omega_{\beta}$
& \lengthtwocolumnA{5}{7}\\
\hline
\rule[-5mm]{0mm}{12mm} \lengthtwocolumnG{5}{6}{5}{6}{5}{6}{5}{6} 
& $-\omega_{\beta}$
& \lengthtwocolumnA{6}{7}\\
\hline
\end{tabular}
\end{center}
%\end{figure} 

%=============================
% End pagebreak for figure
%=============================

Our final result presents the $\mathfrak{g}$-semistandard lattices as 
answers to Stanley's Problem 3 \cite{StanUnim}: 

\noindent
{\bf \CharacterCorollary}\ \ {\sl  
Let 
$\mathfrak{g}$ 
be a simple Lie algebra 
of rank two. 
Let $\lambda = (a,b)$, with $a, b \geq 0$.  
Then the $\mathfrak{g}$-semistandard lattices $\Lba$ and 
$\Lab$ are rank symmetric and rank unimodal.  Moreover, the rank 
generating functions for these lattices are:}
\[\mathit{RGF}_{A_{2}}(\lambda,q) =  
\frac{(1-q^{a+1})(1-q^{b+1})(1-q^{a+b+2})}{(1-q)(1-q)(1-q^{2})}\]

\vspace*{-0.25in}
\[\mathit{RGF}_{C_{2}}(\lambda,q)  =  
\frac{(1-q^{a+1})(1-q^{b+1})(1-q^{a+b+2})(1-q^{a+2b+3})}
{(1-q)(1-q)(1-q^{2})(1-q^{3})}\]

\vspace*{-0.25in}
\[\mathit{RGF}_{G_{2}}(\lambda,q) =  
\frac{(1-q^{a+1})(1-q^{b+1})(1-q^{a+b+2})(1-q^{a+2b+3})
(1-q^{a+3b+4})(1-q^{2a+3b+5})}{(1-q)(1-q)(1-q^{2})(1-q^{3})
(1-q^{4})(1-q^{5})}\]
{\sl In each case $|\Lba| = |\Lab|$, and these counts may be found by 
letting $q \rightarrow 1$.}

{\em Proof.}  In light of \CharacterProposition, apply \RGFProp. 
We have specialized the right hand side quotient there using the data from 
\RankTwoData.
\hfill\QED

%==================================================================
%\newpage
\vspace{2ex} 

\noindent
{\bf \Large \ExamplesNum.\ \ Remarks}

Stanley's Exercises 4.25 and 3.27 on Gaussian and pleasant posets have 
attracted some attention \cite{Stanley}.  
A poset  $P$  with  $p$  elements is {\em Gaussian} if there exist positive integers  
$h_{1},\ldots,h_{p}$  such that for all  $m \geq 0$,  
the rank generating function of the lattice  $J(P \times \mathbf{m})$  
is  $\Pi_{i=1}^{p}(1-q^{m+h_{i}})/(1-q^{h_{i}})$.  
In \cite{PrEur}, the sixth author and Stanley gave a uniform proof of the 
Gaussian property for all Gaussian posets.   That proof used an analog 
of \CharacterProposition;  
it was based upon Seshadri's standard monomial basis theorem for the 
irreducible representations  $X_{n}(m\omega_{k})$,  where the representations  
$X_{n}(\omega_{k})$  are ``minuscule''.   
Now let  $P$  be our  $G_{2}$-fundamental poset  $P_{G_{2}}(0,1)$  of 
\FundPosets.   Please use Figure 3.2 to help visualize the  
$G_{2}$-semistandard poset  $P_{G_{2}}^{\beta_{\alpha}}(0,m)$  
for  $m \geq 0$.  
Note that $P_{G_{2}}^{\beta\alpha}(0,m)$ consists of  $P \times 
\mathbf{m}$  
together with some additional order relations.  
By \CharacterCorollary, the rank generating function for  
$L_{G_{2}}^{\beta\alpha}(0,m) = 
J_{color}(P_{G_{2}}^{\beta\alpha}(0,m))$  is  
\[\frac{(1 - q^{m+1})(1 - q^{m+2})(1 - q^{2m+3})(1 - q^{3m+4})(1 - 
q^{3m+5})}{(1 - 
q^{1})(1 - q^{2})(1 - q^{3})(1 - q^{4})(1 - q^{5})}.\]  
One could introduce a more general notion of ``quasi-Gaussian'' for a 
poset $P$ by requiring that the 
elements of  $P \times \mathbf{m}$  
remain distinct when some additional (if any) order relations are introduced, 
and by allowing a more general product form for the generating function identity.  
Then the fundamental posets  $P_{C_{2}}(0,1)$  and  $P_{G_{2}}(1,0)$  
of 
\FundPosets\ would also be quasi-Gaussian, but not Gaussian. 
In \cite{DW} more will be said about the order relations added to $P 
\times \mathbf{m}$ above and the juxtaposition rules for the fundamental posets 
shown in \GridPosetsFigureList.  For now, we note that these added order 
relations are similar to those added in the following example: 
The ``Catalan'' poset $P_{3}$ of \IntroFigOne\ can be obtained by adding 
order relations to 
the Gaussian poset $\mathbf{3} \times \mathbf{3}$; this corresponds 
to the restriction of $\mathfrak{sl}_{6}(\omega_{3})$ to 
$\mathfrak{sp}_{6}(\omega_{3})$.  

%\newpage
\begin{figure}[ht]
\begin{center}
\BtwoMolevFigureLattice: $P_{C_{2}}^{\beta\alpha}(1,1)$ and 
$L_{C_{2}}^{\beta\alpha}(1,1)$.\\
{\small (Vertices of $L_{C_{2}}^{\beta\alpha}(1,1)$ are indexed by 
$C_{2}$-semistandard tableaux.)}\\
%\vspace{.2cm} 
\setlength{\unitlength}{1.0cm}
\begin{picture}(5,9)
\put(-1,8){$P_{C_{2}}^{\beta\alpha}(1,1)$}
\put(-1.5,4.25){
\begin{picture}(3,3.5)
%Vertex labels
\put(4,1){\circle*{0.15}} 
\put(4.2,0.9){\footnotesize $\alpha$}
\put(1,0){\circle*{0.15}} 
\put(1.2,-0.1){\footnotesize $\beta$}
\put(3,2){\circle*{0.15}} 
\put(3.2,1.9){\footnotesize $\beta$}
\put(0,1){\circle*{0.15}} 
\put(0.2,0.9){\footnotesize $\alpha$}
\put(1,2){\circle*{0.15}} 
\put(1.2,1.9){\footnotesize $\alpha$}
\put(2,3){\circle*{0.15}} 
\put(2.2,2.9){\footnotesize $\alpha$}
\put(0,3){\circle*{0.15}} 
\put(0.2,2.9){\footnotesize $\beta$}
%Chains
\put(0,1){\line(1,1){2}}
\put(1,0){\line(1,1){2}} 
%Other edges 
\put(1,0){\line(-1,1){1}} 
\put(1,2){\line(-1,1){1}} 
\put(2,3){\line(1,-1){2}}
\end{picture}
}
\end{picture}
\hspace*{0.25in}
\setlength{\unitlength}{1.5cm}
\begin{picture}(4.5,7.5)
\put(-0.75,5){$L_{C_{2}}^{\beta\alpha}(1,1)$}
% edges
\put(1,0){\line(1,1){4}} \put(0,1){\line(1,1){4}}
\put(0,3){\line(1,1){3}} \put(1,6){\line(1,1){1}}
\put(1,0){\line(-1,1){1}} \put(2,1){\line(-1,1){2}}
\put(3,2){\line(-1,1){2}} \put(4,3){\line(-1,1){3}}
\put(5,4){\line(-1,1){3}}
% colors
\put(0.4,0.4){$\alpha$} \put(1.4,0.4){$\beta$} \put(0.4,1.4){$\beta$}
\put(1.4,1.4){$\alpha$} \put(2.4,1.4){$\alpha$} \put(0.4,2.4){$\beta$}
\put(1.4,2.4){$\alpha$} \put(2.4,2.4){$\alpha$} \put(3.4,2.4){$\alpha$}
\put(0.4,3.4){$\alpha$} \put(1.4,3.4){$\beta$} \put(2.4,3.4){$\alpha$}
\put(3.4,3.4){$\alpha$} \put(4.4,3.4){$\beta$} \put(1.4,4.4){$\alpha$}
\put(2.4,4.4){$\beta$} \put(3.4,4.4){$\beta$} \put(4.4,4.4){$\alpha$}
\put(1.4,5.4){$\alpha$} \put(2.4,5.4){$\beta$} \put(3.4,5.4){$\beta$}
\put(1.4,6.4){$\beta$} \put(2.4,6.4){$\alpha$}
% vertices
\put(1,0){\circle*{.125}} \put(0,1){\circle*{.125}}
\put(2,1){\circle*{.125}} \put(1,2){\circle*{.125}}
\put(3,2){\circle*{.125}} \put(0,3){\circle*{.125}}
\put(2,3){\circle*{.125}} \put(4,3){\circle*{.125}}
\put(1,4){\circle*{.125}} \put(3,4){\circle*{.125}}
\put(5,4){\circle*{.125}} \put(2,5){\circle*{.125}}
\put(4,5){\circle*{.125}} \put(1,6){\circle*{.125}}
\put(3,6){\circle*{.125}} \put(2,7){\circle*{.125}}
% vertex index
% vertex 1
\put(1.35,7){\abtab{1}{2}{1}} 
% vertex 2
\put(0.4,5.9){\abtab{1}{3}{1}}
% vertex 3
\put(3.3,5.9){\abtab{1}{2}{2}} 
% vertex 4
\put(1.3,4.8){\abtab{1}{3}{2}}
% vertex 5
\put(4.3,4.9){\abtab{1}{2}{3}} 
% vertex 6
\put(0.4,3.9){\abtab{2}{3}{2}}
% vertex 7
\put(3.35,3.9){\abtab{1}{3}{3}} 
% vertex 8
\put(5.2,3.9){\abtab{1}{2}{4}}
% vertex 9
\put(-0.6,2.9){\abtab{2}{4}{2}} 
% vertex 10
\put(2.35,2.9){\abtab{2}{3}{3}}
% vertex 11
\put(4.3,2.9){\abtab{1}{3}{4}} 
% vertex 12
\put(0.25,1.8){\abtab{2}{4}{3}}
% vertex 13
\put(3.3,1.9){\abtab{2}{3}{4}} 
% vertex 14
\put(-0.6,0.9){\abtab{3}{4}{3}}
% vertex 15
\put(2.3,0.9){\abtab{2}{4}{4}} 
% vertex 16
\put(1.2,-0.2){\abtab{3}{4}{4}}
\end{picture}
\end{center}

\vspace*{-0.15in}
\end{figure}

Here is the $C_{2}$ example promised in the middle of Section 
\SemiNum:  The 
$C_{2}$-semistandard poset $P_{C_{2}}^{\beta\alpha}(1,1)$ is displayed in 
\BtwoMolevFigureLattice.  Also displayed is the 
corresponding $C_{2}$-semistandard lattice $L_{C_{2}}^{\beta\alpha}(1,1)$, 
with vertices labelled by the $C_{2}$-semistandard tableaux of 
shape $(1,1)$.  The lattice $L_{C_{2}}^{\beta\alpha}(1,1)$ 
shown in \BtwoMolevFigureLattice\ looks similar in structure to the 
edge-colored lattice  $L$  displayed in \NotLieComboLattice.  In fact, this  $L = 
J_{color}(P)$  for the two-color grid poset  $P$  displayed in \NotLieFundPosets.  
Moreover,  $P = Q_{1} \triangleleft Q_{2}$  with  $Q_{1} \cong P_{1}$ and $Q_{2} 
\cong P_{2}$ for the 
indecomposable two-color grid posets  $P_{1}$ and $P_{2}$  displayed in \NotLieFundPosets.  
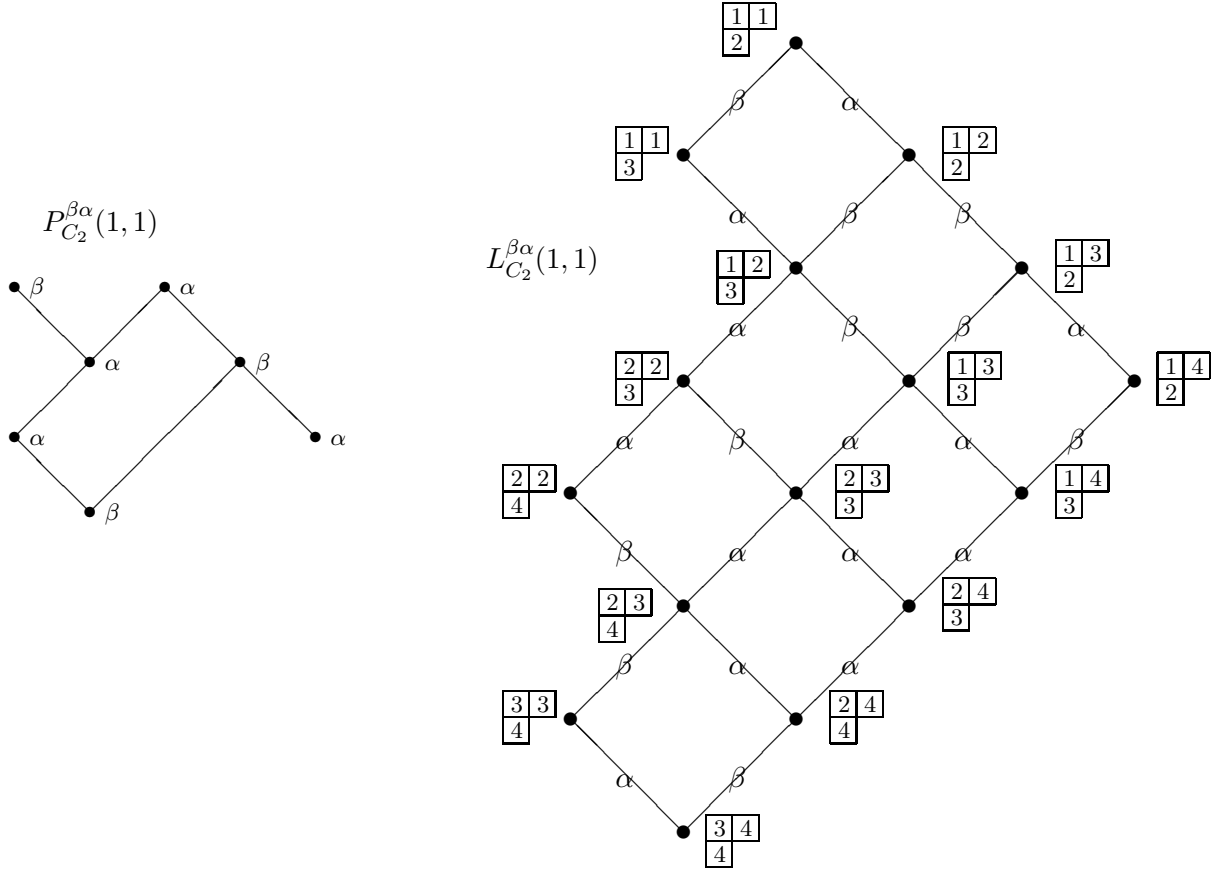
And  $P_{1}$ and $P_{2}$  look similar in structure to the fundamental  
$\mathfrak{g}$-semistandard posets presented in \FundPosets.  But it can be seen that  
$L$  does not satisfy the structure condition for any  $2 \times 2$  
matrix  $M$.  
Therefore  $L$  cannot be a splitting poset for a representation, and so 
there is no hope of applying \CharacterCorollary\ to $L$.  But  $L$  does have a 
``symmetric chain decomposition,'' and hence it is rank symmetric, rank unimodal, 
and ``strongly Sperner''.

It is possible to prove that the $\mathfrak{g}$-semistandard 
lattices, $\mathfrak{g} \in 
\{A_{1} \oplus A_{1}, A_{2}, C_{2}, G_{2}\}$, are the only lattices of the 
kind we have been 
considering which can have the $M$-structure property for any 
$2 \times 2$ integer 
matrix  $M$:

\noindent
{\bf \ClassificationTheoremOne\ \cite{DonTwoColor}}\ \ {\sl 
 Let  $P$  be a two-color grid poset which has the max 
property.  If  $L = J_{color}(P)$  has the $M$-structure property for some $2 \times 2$ integer 
matrix  $M$,  then  $L$  is a $\mathfrak{g}$-semistandard lattice,  
$\mathfrak{g} \in 
\{A_{1} \oplus A_{1}, A_{2}, C_{2}, G_{2}\}$.}

%\newpage
\begin{figure}[ht]
\begin{center}
\NotLieFundPosets: 
We can write $P = Q_{1} \triangleleft Q_{2}$ with $Q_{i} \cong 
P_{i}$ ($i = 1,2$).  Below, $L = J_{color}(P)$.\\
\vspace{.2cm} 
\setlength{\unitlength}{0.8cm}
\begin{picture}(10,6)
\put(-2.5,0.25){
\begin{picture}(3,3.5)
\put(1,4){$P_{1}$}
%Vertex labels
\put(1,0){\circle*{0.15}} 
\put(1.2,-0.1){\footnotesize $\beta$}
\put(0,1){\circle*{0.15}} 
\put(0.2,0.9){\footnotesize $\alpha$}
\put(1,2){\circle*{0.15}} 
\put(1.2,1.9){\footnotesize $\alpha$}
\put(2,3){\circle*{0.15}} 
\put(2.2,2.9){\footnotesize $\alpha$}
\put(3,4){\circle*{0.15}} 
\put(3.2,3.9){\footnotesize $\alpha$}
\put(2,5){\circle*{0.15}} 
\put(2.2,4.9){\footnotesize $\beta$}
%Chains
\put(0,1){\line(1,1){3}}
%Other edges 
\put(1,0){\line(-1,1){1}} 
\put(3,4){\line(-1,1){1}} 
\end{picture}
}
\put(0,-0.25){
\begin{picture}(3,3.5)
\put(3,4){$P_{2}$}
%Vertex labels
\put(4,1){\circle*{0.15}} 
\put(4.2,0.9){\footnotesize $\alpha$}
\put(3,2){\circle*{0.15}} 
\put(3.2,1.9){\footnotesize $\beta$}
\put(4,3){\circle*{0.15}} 
\put(4.2,2.9){\footnotesize $\beta$}
\put(5,4){\circle*{0.15}} 
\put(5.2,3.9){\footnotesize $\beta$}
\put(4,5){\circle*{0.15}} 
\put(4.2,4.9){\footnotesize $\alpha$}
%Chains
\put(3,2){\line(1,1){2}} 
%Other edges 
\put(5,4){\line(-1,1){1}} 
\put(4,1){\line(-1,1){1}} 
\end{picture}
}
\put(7.5,0.25){
\begin{picture}(3,3.5)
\put(1,4){$P$}
%Vertex labels
\put(4,1){\circle*{0.15}} 
\put(4.2,0.9){\footnotesize $\alpha$}
\put(1,0){\circle*{0.15}} 
\put(1.2,-0.1){\footnotesize $\beta$}
\put(3,2){\circle*{0.15}} 
\put(3.2,1.9){\footnotesize $\beta$}
\put(4,3){\circle*{0.15}} 
\put(4.2,2.9){\footnotesize $\beta$}
\put(5,4){\circle*{0.15}} 
\put(5.2,3.9){\footnotesize $\beta$}
\put(0,1){\circle*{0.15}} 
\put(0.2,0.9){\footnotesize $\alpha$}
\put(1,2){\circle*{0.15}} 
\put(1.2,1.9){\footnotesize $\alpha$}
\put(2,3){\circle*{0.15}} 
\put(2.2,2.9){\footnotesize $\alpha$}
\put(3,4){\circle*{0.15}} 
\put(3.2,3.9){\footnotesize $\alpha$}
\put(4,5){\circle*{0.15}} 
\put(4.2,4.9){\footnotesize $\alpha$}
\put(2,5){\circle*{0.15}} 
\put(2.2,4.9){\footnotesize $\beta$}
%Chains
\put(0,1){\line(1,1){4}}
\put(1,0){\line(1,1){4}} 
%Other edges 
\put(1,0){\line(-1,1){1}} 
\put(5,4){\line(-1,1){1}} 
\put(4,1){\line(-1,1){1}} 
\put(3,4){\line(-1,1){1}} 
\end{picture}
}
\end{picture}
%\end{center}
%\end{figure}

\setlength{\unitlength}{0.8cm}
\begin{picture}(9,11.5)
\put(2,9){$L$}
% vertices
\put(3,0){\circle*{0.15}}
\put(2,1){\circle*{0.15}}
\put(4,1){\circle*{0.15}}
\put(3,2){\circle*{0.15}}
\put(5,2){\circle*{0.15}}
\put(2,3){\circle*{0.15}}
\put(4,3){\circle*{0.15}}
\put(6,3){\circle*{0.15}}
\put(1,4){\circle*{0.15}}
\put(3,4){\circle*{0.15}}
\put(5,4){\circle*{0.15}}
\put(7,4){\circle*{0.15}}
\put(0,5){\circle*{0.15}}
\put(2,5){\circle*{0.15}}
\put(4,5){\circle*{0.15}}
\put(6,5){\circle*{0.15}}
\put(8,5){\circle*{0.15}}
\put(1,6){\circle*{0.15}}
\put(3,6){\circle*{0.15}}
\put(5,6){\circle*{0.15}}
\put(7,6){\circle*{0.15}}
\put(9,6){\circle*{0.15}}
\put(2,7){\circle*{0.15}}
\put(4,7){\circle*{0.15}}
\put(6,7){\circle*{0.15}}
\put(8,7){\circle*{0.15}}
\put(3,8){\circle*{0.15}}
\put(5,8){\circle*{0.15}}
\put(7,8){\circle*{0.15}}
\put(4,9){\circle*{0.15}}
\put(6,9){\circle*{0.15}}
\put(3,10){\circle*{0.15}}
\put(5,10){\circle*{0.15}}
\put(4,11){\circle*{0.15}}
% colors, edges
\put(3,0){\line(-1,1){1}}
\put(3,0){\NWEdgeLabelForLatticeI{\alpha}}
\put(3,0){\line(1,1){1}}
\put(3,0){\NEEdgeLabelForLatticeIII{\beta}}
%=====
\put(2,1){\line(1,1){1}}
\put(2,1){\NEEdgeLabelForLatticeIII{\beta}}
%=====
\put(4,1){\line(-1,1){1}}
\put(4,1){\NWEdgeLabelForLatticeIII{\alpha}}
\put(4,1){\line(1,1){1}}
\put(4,1){\NEEdgeLabelForLatticeIII{\alpha}}
%=====
\put(3,2){\line(-1,1){1}}
\put(3,2){\NWEdgeLabelForLatticeIII{\beta}}
\put(3,2){\line(1,1){1}}
\put(3,2){\NEEdgeLabelForLatticeIII{\alpha}}
%=====
\put(5,2){\line(-1,1){1}}
\put(5,2){\NWEdgeLabelForLatticeIII{\alpha}}
\put(5,2){\line(1,1){1}}
\put(5,2){\NEEdgeLabelForLatticeIII{\alpha}}
%=====
\put(2,3){\line(-1,1){1}}
\put(2,3){\NWEdgeLabelForLatticeIII{\beta}}
\put(2,3){\line(1,1){1}}
\put(2,3){\NEEdgeLabelForLatticeIII{\alpha}}
%=====
\put(4,3){\line(-1,1){1}}
\put(4,3){\NWEdgeLabelForLatticeIII{\beta}}
\put(4,3){\line(1,1){1}}
\put(4,3){\NEEdgeLabelForLatticeIII{\alpha}}
%=====
\put(6,3){\line(-1,1){1}}
\put(6,3){\NWEdgeLabelForLatticeIII{\alpha}}
\put(6,3){\line(1,1){1}}
\put(6,3){\NEEdgeLabelForLatticeIII{\alpha}}
%=====
\put(1,4){\line(-1,1){1}}
\put(1,4){\NWEdgeLabelForLatticeIII{\beta}}
\put(1,4){\line(1,1){1}}
\put(1,4){\NEEdgeLabelForLatticeIII{\alpha}}
%=====
\put(3,4){\line(-1,1){1}}
\put(3,4){\NWEdgeLabelForLatticeIII{\beta}}
\put(3,4){\line(1,1){1}}
\put(3,4){\NEEdgeLabelForLatticeIII{\alpha}}
%=====
\put(5,4){\line(-1,1){1}}
\put(5,4){\NWEdgeLabelForLatticeIII{\beta}}
\put(5,4){\line(1,1){1}}
\put(5,4){\NEEdgeLabelForLatticeIII{\alpha}}
%=====
\put(7,4){\line(-1,1){1}}
\put(7,4){\NWEdgeLabelForLatticeIII{\alpha}}
\put(7,4){\line(1,1){1}}
\put(7,4){\NEEdgeLabelForLatticeIII{\alpha}}
%=====
\put(0,5){\line(1,1){1}}
\put(0,5){\NEEdgeLabelForLatticeIII{\alpha}}
%=====
\put(2,5){\line(-1,1){1}}
\put(2,5){\NWEdgeLabelForLatticeIII{\beta}}
\put(2,5){\line(1,1){1}}
\put(2,5){\NEEdgeLabelForLatticeIII{\alpha}}
%=====
\put(4,5){\line(-1,1){1}}
\put(4,5){\NWEdgeLabelForLatticeIII{\beta}}
\put(4,5){\line(1,1){1}}
\put(4,5){\NEEdgeLabelForLatticeIII{\alpha}}
%=====
\put(6,5){\line(-1,1){1}}
\put(6,5){\NWEdgeLabelForLatticeIII{\beta}}
\put(6,5){\line(1,1){1}}
\put(6,5){\NEEdgeLabelForLatticeIII{\alpha}}
%=====
\put(8,5){\line(-1,1){1}}
\put(8,5){\NWEdgeLabelForLatticeIII{\alpha}}
\put(8,5){\line(1,1){1}}
\put(8,5){\NEEdgeLabelForLatticeIII{\beta}}
%=====
\put(1,6){\line(1,1){1}}
\put(1,6){\NEEdgeLabelForLatticeIII{\alpha}}
%=====
\put(3,6){\line(-1,1){1}}
\put(3,6){\NWEdgeLabelForLatticeIII{\beta}}
\put(3,6){\line(1,1){1}}
\put(3,6){\NEEdgeLabelForLatticeIII{\alpha}}
%=====
\put(5,6){\line(-1,1){1}}
\put(5,6){\NWEdgeLabelForLatticeIII{\beta}}
\put(5,6){\line(1,1){1}}
\put(5,6){\NEEdgeLabelForLatticeIII{\alpha}}
%=====
\put(7,6){\line(-1,1){1}}
\put(7,6){\NWEdgeLabelForLatticeIII{\beta}}
\put(7,6){\line(1,1){1}}
\put(7,6){\NEEdgeLabelForLatticeIII{\beta}}
%=====
\put(9,6){\line(-1,1){1}}
\put(9,6){\NWEdgeLabelForLatticeIII{\alpha}}
%=====
\put(2,7){\line(1,1){1}}
\put(2,7){\NEEdgeLabelForLatticeIII{\alpha}}
%=====
\put(4,7){\line(-1,1){1}}
\put(4,7){\NWEdgeLabelForLatticeIII{\beta}}
\put(4,7){\line(1,1){1}}
\put(4,7){\NEEdgeLabelForLatticeIII{\alpha}}
%=====
\put(6,7){\line(-1,1){1}}
\put(6,7){\NWEdgeLabelForLatticeIII{\beta}}
\put(6,7){\line(1,1){1}}
\put(6,7){\NEEdgeLabelForLatticeIII{\beta}}
%=====
\put(8,7){\line(-1,1){1}}
\put(8,7){\NWEdgeLabelForLatticeIII{\beta}}
%=====
\put(3,8){\line(1,1){1}}
\put(3,8){\NEEdgeLabelForLatticeIII{\alpha}}
%=====
\put(5,8){\line(-1,1){1}}
\put(5,8){\NWEdgeLabelForLatticeIII{\beta}}
\put(5,8){\line(1,1){1}}
\put(5,8){\NEEdgeLabelForLatticeIII{\beta}}
%=====
\put(7,8){\line(-1,1){1}}
\put(7,8){\NWEdgeLabelForLatticeIII{\beta}}
%=====
\put(4,9){\line(-1,1){1}}
\put(4,9){\NWEdgeLabelForLatticeIII{\alpha}}
\put(4,9){\line(1,1){1}}
\put(4,9){\NEEdgeLabelForLatticeIII{\beta}}
%=====
\put(6,9){\line(-1,1){1}}
\put(6,9){\NWEdgeLabelForLatticeIII{\beta}}
%=====
\put(3,10){\line(1,1){1}}
\put(3,10){\NEEdgeLabelForLatticeIII{\beta}}
%=====
\put(5,10){\line(-1,1){1}}
\put(5,10){\NWEdgeLabelForLatticeIII{\alpha}}
%=====
\end{picture}
\end{center}

\vspace*{-0.25in}
\end{figure}

\noindent 
{\bf \ClassificationTheoremTwo\ \cite{DonTwoColor}}\ \ {\sl 
Let  $P$  be an indecomposable two-color grid poset.  
If  $L = J_{color}(P)$  has the $M$-structure property for some $2 \times 2$ integer 
matrix  $M$,  then  $L$  is a $\mathfrak{g}$-fundamental lattice,  
$\mathfrak{g} \in 
\{A_{1} \oplus A_{1}, A_{2}, C_{2}, G_{2}\}$.}

\noindent 
These two statements are combinatorial Dynkin diagram classification 
theorems:  No Lie theory or algebraic concepts of any kind appear in their 
hypotheses, but the short list of Dynkin diagram-indexed rank two Cartan 
matrices plays the central role in their conclusions.

Let  $P$  be a two-color grid poset and set  $L = J_{color}(P)$.  
To apply \CharacterCorollary\ to  $L$  
via \CharacterProposition, \WeightPresBijection\ was required: 
the elements of the $\mathfrak{g}$-semistandard lattices 
matched up with 
tableaux of Littelmann.  The precise match-up required here should make 
one pessimistic about obtaining the rank generating function identities of 
\CharacterCorollary\ for 
these more general  $L$. This pessimism is 
intuitively heightened by the classification results above, which 
emphasize how special the $\mathfrak{g}$-semistandard lattices are.  
(But it is possible to 
obtain the total count results mentioned at the end of \CharacterCorollary\ with 
elementary combinatorial reasoning \cite{ADLP}.)  
After representations for the cases listed in the 
introduction to this paper are constructed, Corollary 5.3 of \cite{ADLP} 
notes that the 
$\mathfrak{g}$-semistandard lattices in 
those cases are strongly Sperner.  Although this approach cannot be used 
for the rest of the 
$\mathfrak{g}$-semistandard lattices, it is natural to hope that 
those lattices have this property.  When addressing these 
extremal set theory issues, here it would now seem reasonable to attempt a 
combinatorial approach:  Can one find symmetric chain decompositions of  
$L = J_{color}(P)$  for certain two-color grid posets  $P$?  

Although the rank two cases in Lie theory are much simpler than the 
general rank cases, it is also true in Lie theory that the key aspect of a 
higher rank case often reduces to consideration of that aspect for just 
the rank two cases.  
Various aspects of this rank two paper will be used for many higher rank cases 
in \cite{DW}. 
The forms of the  
$\mathfrak{g}$-semistandard tableaux of Section \SemiNum\ 
may seem unmotivated to readers who 
are familiar only with \cite{Hum}.  Space and time permitting, much motivation 
could be supplied.  Strict columns of length two arise because the second 
fundamental representation in each simple case can be realized as 
the ``big 
piece'' of the second exterior power of the first fundamental 
representation.  Standard monomial theory (and earlier papers concerning 
algebras with straightening laws) explain how the restricted concatenation 
of columns corresponds to the multiplication of ``Pl\"{u}cker coordinates'' for 
flag manifolds.  Going further, it may be possible to ``explain'' the simple 
root colorings of the elements of the posets  $P$  in the spirit of the 
heaps of Stembridge, along the lines of Theorem 11.1 of \cite{PrEur}.

Our main result states that the $\mathfrak{g}$-semistandard lattices are 
splitting posets for their representations.  For any representation, 
the crystal graph (of Kashiwara) is a splitting poset 
(cf.\ Lemma 3.6 of  
\cite{DonSupp}).  More generally, this is true for 
Stembridge's overarching crystal graph-like ``admissible systems'' 
\cite{Stem}.  The second author has observed that any such ``general crystal graph'' 
for a given representation is ``edge minimal'' within the set of splitting posets 
for the representation:  It contains no splitting poset for the representation 
as a proper subgraph.  For all irreducible representations of types $A_{2}$ 
and $C_{2}$, it can be seen from \cite{KN} that Kashiwara's crystal graphs are 
subgraphs of the corresponding  $\mathfrak{g}$-semistandard lattices.  
It seems likely that all $\mathfrak{g}$-semistandard lattices give rise to 
admissible systems.  This would directly show that the elements of these lattices 
generate the associated Weyl characters.  In \cite{DW} we will consider most simple Lie 
algebras of arbitrary rank and uniformly define $\mathfrak{g}$-fundamental posets for 
their fundamental weights which have the following property:  the longest 
element in the associated Bruhat order is ``fully commutative''.  This definition is 
type-independent.  Using these fundamental posets, as in Section 
\SemiNum\ we build  $\mathfrak{g}$-semistandard posets and lattices for many 
representations.  Along with this paper, this should start a new program:  Find modular 
lattice splitting posets for all irreducible representations of all semisimple Lie 
algebras and show that they give rise to admissible systems.  If these hopes are 
realized, these modular lattices (including the  $\mathfrak{g}$-semistandard lattices) 
would in general contain ``extra'' edges with respect to the admissible system.  
But the lattices might be more combinatorially interesting than most or all admissible systems' 
directed graphs, and hopefully more accessible.  One consequence might be the formulation of 
analogs of the Littlewood-Richardson tensor product rule in terms of manipulations of the 
underlying  $\mathfrak{g}$-semistandard posets (or their analogs in the modular/non-distributive cases).

%=============================================
% Bibliography
%
%=============================================
%\newpage%
\vspace*{-0.185in}
\renewcommand{\baselinestretch}{1}


\small\normalsize
\begin{thebibliography}{999999}

{\small
\bibitem[Alv]{Alv} L. W. Alverson II, ``Distributive lattices and
representations of the rank two simple Lie algebras,''
Master's thesis, Murray State University, 2003.

\bibitem[ADLP]{ADLP} L.\ W.\ Alverson II, R.\ G.\ Donnelly, S.\ 
J.\ Lewis, and R.\ Pervine, ``Constructions of representations of 
rank two semisimple Lie algebras with distributive lattices,'' 
{\em Electronic J.\ Combin.\ }, 
{\bf 13} (2006), \#R109 (44 pp). 

\bibitem[Don1]{DonSupp} R.\ G.\ Donnelly, ``Extremal properties of
bases for representations
of semisimple Lie algebras,'' {\em J.\ Algebraic Comb.\ } {\bf 17}
(2003), 255--282.

\bibitem[Don2]{DonTwoColor} R.\ G.\ Donnelly, ``A Dynkin diagram 
classification of posets satisfying certain structural properties,'' 
preprint.

\bibitem[DLP]{DLP1} R.\ G.\ Donnelly, S.\ J.\ Lewis, and R. Pervine,
``Constructions of representations of $\mathfrak{o}(2n+1,\mathbb{C})$
that imply Molev and Reiner-Stanton lattices are strongly Sperner,''
{\em Discrete Math.\ } {\bf 263} (2003), 61--79.

\bibitem[DW]{DW} R.\ G.\ Donnelly and N.\ J.\ Wildberger,
``Distributive lattice models for certain families of
irreducible semisimple Lie algebra representations,'' in
preparation. 

\bibitem[Hum]{Hum} J.\ E.\ Humphreys, {\em Introduction to Lie
Algebras and Representation Theory}, Springer, New York, 1972.

\bibitem[KN]{KN} M.\ Kashiwara and T.\ Nakashima, ``Crystal
graphs for representations of the $q$-analogue of classical Lie
algebras,'' {\em J.\ Algebra\ } {\bf 165} (1994), 295--345.

\bibitem[Lit]{Lit1} P.\ Littelmann, ``A generalization of the 
Littlewood-Richardson rule,'' {\em J.\ Algebra\ } {\bf 130} (1990),
328--368.

\bibitem[Mc]{McClard} M.\ McClard, ``Picturing Representations of
Simple Lie
Algebras of Rank
Two,'' Master's thesis, Murray State University, 2000.

\bibitem[Pr1]{PrEur} R.\ A.\ Proctor, ``Bruhat lattices, plane
partition generating functions, and minuscule representations,'' {\em
European J.\ Combin.\ }{\bf 5} (1984), 331-350.

\bibitem[Pr2]{PrGZ} R.\ A.\ Proctor, ``Solution of a Sperner
conjecture of Stanley with a construction of Gelfand,'' {\em J.\ Comb.\
Th.\ A} {\bf 54} (1990), 225-234.

\bibitem[Sta1]{StanUnim} R.\ P.\ Stanley, ``Unimodal sequences arising
from Lie algebras,'' in: {\em Young Day Proceedings}, ed. T.\ V.\ Narayana
{\em et al}, Marcel Dekker, New York, 1980, 127-136.

\bibitem[Sta2]{Stanley} R.\ P.\ Stanley, {\em Enumerative Combinatorics,
Vol.\ 1}, Wadsworth and Brooks/Cole, Monterey, CA, 1986.

\bibitem[Sta3]{StanText2} R.\ P.\ Stanley, {\em Enumerative Combinatorics,
Vol.\ 2}, Cambridge University Press, 1999.

\bibitem[Stem]{Stem} J.\ Stembridge, ``Combinatorial models for Weyl
characters,'' {\em Advances in Math.\ } {\bf 168} (2002), 96--131.

}
\end{thebibliography}
\end{document}